\documentclass{amsart}
\usepackage[dvips]{graphicx}
\usepackage{amscd}
\usepackage{amsmath}
\usepackage{amsxtra}
\usepackage{amsfonts}
\usepackage{amssymb}
\newtheorem{theorem}{Theorem}[section]
\newtheorem{corollary}[theorem]{Corollary}
\newtheorem{lemma}[theorem]{Lemma}
\newtheorem{proposition}[theorem]{Proposition}
\theoremstyle{definition}
\newtheorem{definition}[theorem]{Definition}
\newtheorem{remark}[theorem]{Remark}

\newtheorem{example}[theorem]{Example}
\theoremstyle{remark}
\newtheorem{claim}[theorem]{Claim}
\renewcommand{\theclaim}{\textup{\theclaim}}
\newtheorem*{claims}{Claims}
\newtheorem*{acknowledgements}{Acknowledgements}

\renewcommand{\theenumi}{\roman{enumi}}

\numberwithin{equation}{section}
\newlength{\qedskip}

\def\hooklongrightarrow{\lhook\joinrel\longrightarrow} 
\def\longequal{=\joinrel =} 
\def\openone
{\mathchoice
{\hbox{\upshape \small1\kern-3.3pt\normalsize1}}
{\hbox{\upshape \small1\kern-3.3pt\normalsize1}}
{\hbox{\upshape \tiny1\kern-2.3pt\SMALL1}}
{\hbox{\upshape \Tiny1\kern-2pt\tiny1}}}
\makeatletter
\newbox\ipbox
\newcommand{\ip}[2]{\left\langle #1\mathrel{\mathchoice
{\setbox\ipbox=\hbox{$\displaystyle \left\langle\mathstrut #1#2\right\rangle$}
\vrule height\ht\ipbox width0.25pt depth\dp\ipbox}
{\setbox\ipbox=\hbox{$\textstyle \left\langle\mathstrut #1#2\right\rangle$}
\vrule height\ht\ipbox width0.25pt depth\dp\ipbox}
{\setbox\ipbox=\hbox{$\scriptstyle \left\langle\mathstrut #1#2\right\rangle$}
\vrule height\ht\ipbox width0.25pt depth\dp\ipbox}
{\setbox\ipbox=\hbox{$\scriptscriptstyle \left\langle\mathstrut #1#2\right\rangle$}
\vrule height\ht\ipbox width0.25pt depth\dp\ipbox}
} #2\right\rangle}
\newcommand{\diracb}[1]{\left\langle #1\mathrel{\mathchoice
{\setbox\ipbox=\hbox{$\displaystyle \left\langle\mathstrut #1\right.$}
\vrule height\ht\ipbox width0.25pt depth\dp\ipbox}
{\setbox\ipbox=\hbox{$\textstyle \left\langle\mathstrut #1\right.$}
\vrule height\ht\ipbox width0.25pt depth\dp\ipbox}
{\setbox\ipbox=\hbox{$\scriptstyle \left\langle\mathstrut #1\right.$}
\vrule height\ht\ipbox width0.25pt depth\dp\ipbox}
{\setbox\ipbox=\hbox{$\scriptscriptstyle \left\langle\mathstrut #1\right.$}
\vrule height\ht\ipbox width0.25pt depth\dp\ipbox}
}\right. }
\newcommand{\dirack}[1]{\left. \mathrel{\mathchoice
{\setbox\ipbox=\hbox{$\displaystyle \left.\mathstrut #1\right\rangle$}
\vrule height\ht\ipbox width0.25pt depth\dp\ipbox}
{\setbox\ipbox=\hbox{$\textstyle \left.\mathstrut #1\right\rangle$}
\vrule height\ht\ipbox width0.25pt depth\dp\ipbox}
{\setbox\ipbox=\hbox{$\scriptstyle \left.\mathstrut #1\right\rangle$}
\vrule height\ht\ipbox width0.25pt depth\dp\ipbox}
{\setbox\ipbox=\hbox{$\scriptscriptstyle \left.\mathstrut #1\right\rangle$}
\vrule height\ht\ipbox width0.25pt depth\dp\ipbox}
} #1\right\rangle}

\input cyracc.def

\begin{document}
\title[Ruelle operators]{Ruelle operators: Functions which are harmonic with respect to a transfer operator}
\author{Palle~E.~T.~ Jorgensen}
\address{Department of Mathematics\\
The University of Iowa\\
14 MacLean Hall\\
Iowa City, IA 52242-1419\\
U.S.A.}
\email{jorgen@math.uiowa.edu}
\thanks{Work supported in part by the U.S. National Science Foundation.}
\subjclass{Primary 46L60, 47D25, 42A16, 43A65; Secondary 46L45, 42A65, 41A15}
\keywords{$C^{\ast}$-algebra, endomorphism, wavelet, cascade algorithm, refinement
operator, representation, orthogonal expansion, quadrature mirror filter,
isometry in Hilbert space}

\begin{abstract}
Let $N\in\mathbb{N}$, $N\geq2$, be given. Motivated by wavelet analysis, we
consider a class of normal representations of the $C^{\ast}$-algebra
$\mathfrak{A}_{N}$ on two unitary generators $U$, $V$ subject to the relation
\newline
\begin{minipage}{\textwidth}
\hskip-36pt\begin{minipage}{\textwidth}
\begin{equation*}
UVU^{-1}=V^{N}.
\end{equation*}
\end{minipage}\hskip-36pt\vspace*{\abovedisplayskip}
\end{minipage}
The representations are in one-to-one correspondence with solutions $h\in
L^{1}\left(  \mathbb{T}\right)  $, $h\geq0$, to $R\left(  h\right)  =h$ where
$R$ is a certain transfer operator (positivity-preserving) which was studied
previously by D.~Ruelle. The representations of $\mathfrak{A}_{N}$ may also be
viewed as representations of a certain (discrete) $N$-adic $ax+b$ group which
was considered recently by J.-B. Bost and A.~Connes.
\end{abstract}\maketitle
\tableofcontents

\section{\label{Intr}Introduction}

In multiresolution wavelet theory, there is a fundamental interplay and
interconnection between the following two operators: $M$ and $R$, where $M$ is
the \emph{cascade refinement operator} and $R$ is the corresponding
\emph{transfer operator.} We also denote the second of these, $R$, the
\emph{Ruelle operator,} because of its close connection to an operator that
David Ruelle used first in his study of phase transitions in quantum
statistical mechanics lattice models;
see \cite{Rue68}, \cite{May80}, and \cite{Rue78a}.
We first recall the two operators
here in a simple
wavelet
context, but the scope will be widened later (Chapter
\ref{Kean}): Let $M$ be the operator in $L^{2}\left(  \mathbb{R}\right)  $
given by%
\begin{equation}
\left(  M\psi\right)  \left(  x\right)  =\sqrt{N}\sum_{k\in\mathbb{Z}}%
a_{k}\psi\left(  Nx-k\right)  , \label{eqIntr.1}%
\end{equation}
where $N\geq2$ is integral, and $a_{k}\in\mathbb{C}$, $k\in\mathbb{Z}$, are
given subject to $\sum_{k}\left|  a_{k}\right|  ^{2}=1$; $\psi\in L^{2}\left(
\mathbb{R}\right)  $, $x\in\mathbb{R}$. If%
\begin{equation}
m_{0}\left(  z\right)  :=\sum_{k\in\mathbb{Z}}a_{k}z^{k},\qquad z\in
\mathbb{T}, \label{eqIntr.2}%
\end{equation}
then the condition $m_{0}\in L^{\infty}\left(  \mathbb{T}\right)  $ implies
that $M$ is a bounded operator in $L^{2}\left(  \mathbb{R}\right)  $. With the
usual identification
\begin{equation}
\mathbb{R}\diagup2\pi\mathbb{Z}\ni\omega\longmapsto e^{-i\omega}%
=z\in\mathbb{T} \label{eqIntr.3}%
\end{equation}
we have the corresponding identification $m_{0}\left(  \omega\right)
=m_{0}\left(  z\right)  $ of $2\pi$-periodic functions on $\mathbb{R}$ with
functions on $\mathbb{T}$, and we shall use the same letter denoting the
function either way. Introducing the Fourier transform $\psi\mapsto\hat{\psi}$
in $L^{2}\left(  \mathbb{R}\right)  $, we get (\ref{eqIntr.1}) in the
equivalent form:%
\begin{equation}
\left(  M\psi\right)  \sphat\left(  \omega\right)  =\frac{m_{0}\left(
\frac{\omega}{N}\right)  }{\sqrt{N}\mathstrut}\hat{\psi}\left(  \frac{\omega
}{N}\right)  ,\qquad\omega\in\mathbb{R}. \label{eqIntr.4}%
\end{equation}
The Ruelle transfer operator is defined on $L^{1}\left(  \mathbb{T}\right)  $
by%
\begin{equation}
\left(  Rf\right)  \left(  z\right)  =\frac{1}{N}\sum_{w^{N}=z}\left|
m_{0}\left(  w\right)  \right|  ^{2}f\left(  w\right)  ,\qquad f\in
L^{1}\left(  \mathbb{T}\right)  ,\;z\in\mathbb{T}, \label{eqIntr.5}%
\end{equation}
where the summation is over the $N$ roots $w$, i.e., the $N$ solutions to
$w^{N}=z$.

For the \emph{quadrature wavelet filters,} there is the further restriction%
\begin{equation}
\sum_{w^{N}=z}\left|  m_{0}\left(  w\right)  \right|  ^{2}=N, \label{eqIntr.6}%
\end{equation}
or equivalently%
\begin{equation}
R\left(  \openone\right)  =\openone\label{eqIntr.7}%
\end{equation}
where $\openone$ denotes the constant function in $L^{2}\left(  \mathbb{T}%
\right)  $. But our analysis will \emph{not} be restricted to this special case.

When $m_{0}\in L^{\infty}\left(  \mathbb{T}\right)  $ is given, and $R$ is the
corresponding Ruelle operator, we study the eigenvalue problem%
\begin{equation}
h\in L^{1}\left(  \mathbb{T}\right)  ,\qquad h\geq0,\qquad R\left(  h\right)
=h. \label{eqIntr.8}%
\end{equation}
But without (\ref{eqIntr.6}), a nonzero solution $h$ to (\ref{eqIntr.8}) is
not then guaranteed. The problem (\ref{eqIntr.8}) is closely connected to the
problem%
\begin{equation}
\varphi\in L^{2}\left(  \mathbb{R}\right)  ,\qquad M\varphi=\varphi,
\label{eqIntr.9}%
\end{equation}
whose nonzero solutions (if any) are the \emph{scaling functions} (or
\emph{father functions\/}) in wavelet theory.

Suppose (\ref{eqIntr.6}) is given: then a famous argument of Mallat
\cite{Mal89} states that the $L^{2}\left(  \mathbb{R}\right)  $-norm of the
functions $F_{n}$,%
\begin{equation}
F_{n}\left(  \omega\right)  =\prod_{k=1}^{n}\frac{m_{0}\left(  \frac{\omega
}{N^{k\mathstrut}}\right)  }{\sqrt{N}\mathstrut}\chi_{\left[  -\pi
,\pi\right\rangle }^{{}}\left(  \frac{\omega}{N^{n\mathstrut}}\right)
\label{eqIntr.10}%
\end{equation}
is constant. In fact $\left\|  F_{n}\right\|  _{L^{2}\left(  \mathbb{R}%
\right)  }=1$, $n=1,2,\dots$. If moreover $\omega\mapsto m_{0}\left(
\omega\right)  $ is Lipschitz near $\omega=0$, and $m_{0}\left(  0\right)
=\sqrt{N}$ (the low-pass condition), then%
\begin{equation}
F\left(  \omega\right)  =\prod_{k=1}^{\infty}\frac{m_{0}\left(  \frac{\omega
}{N^{k\mathstrut}}\right)  }{\sqrt{N}\mathstrut} \label{eqIntr.11}%
\end{equation}
is pointwise convergent. If the eigenspace $\left\{  h\mid Rh=h\right\}  $ is
further given to be one-dimensional, then $F\in L^{2}\left(  \mathbb{R}%
\right)  $ and $\lim_{n\rightarrow\infty}\left\|  F-F_{n}\right\|
_{L^{2}\left(  \mathbb{R}\right)  }=0$. In view of (\ref{eqIntr.4}), the
inverse Fourier transform $\varphi=\check{F}$ will then solve the eigenvalue
problem (\ref{eqIntr.9}), and $\hat{\varphi}\left(  0\right)  =1$.

To motivate the more general problem (\ref{eqIntr.8}), we note that, if $h\in
L^{1}\left(  \mathbb{T}\right)  $ solves (\ref{eqIntr.8}), then%
\[
\int_{-\pi N^{n}}^{\pi N^{n}}h\left(  \frac{\omega}{N^{n\mathstrut}}\right)
\prod_{k=1}^{n}\frac{\left|  m_{0}\left(  \frac{\omega}{N^{k\mathstrut}%
}\right)  \right|  ^{2}}{N}\,d\omega=\int_{0}^{2\pi}R^{n}h\left(
\omega\right)  \,d\omega=\int_{0}^{2\pi}h\left(  \omega\right)  \,d\omega
\geq0.
\]
So again, if $h\neq0$, then the sequence $F_{n}$, now defined \emph{via} $h$
by%
\begin{equation}
F_{n}\left(  \omega\right)  :=\chi_{\left[  -\pi,\pi\right\rangle }^{{}%
}\left(  \frac{\omega}{N^{n\mathstrut}}\right)  
\left(  h\left(  \frac{\omega}{N^{n\mathstrut}}\right)  \right)  ^{\frac{1}{2}}
\prod_{k=1}^{n}\frac{m_{0}\left(  \frac{\omega
}{N^{k\mathstrut}}\right)  }{\sqrt{N}\mathstrut} \label{eqIntr.12}%
\end{equation}
has constant norm in $L^{2}\left(  \mathbb{R}\right)  $. If, for example,
$h\left(  0\right)  =1$, then it can be checked that $F_{n}\rightarrow F$ in
$L^{2}\left(  \mathbb{R}\right)  $ for $n\rightarrow\infty$, where $F$ has a
(\ref{eqIntr.11})-representation.

For more background references on
wavelets, filters, and scaling functions,
from the operator-theoretic viewpoint,
we give \cite{CoDa96}, \cite{CoRy95}, 
\cite{DaLa}, \cite{Dau92}, \cite{Hor95}, \cite{MePa93}, 
\cite{Mey98}, and \cite{Vil94}.
However,
a summary of some main ideas
is included for the convenience of
the reader, and to make the paper
more self-contained. Our terminology
is close to that of the listed references.

One of the aims of our paper is to widen the scope of the quadrature analysis
and to study the more general eigenvalue problem (\ref{eqIntr.8}). For that,
it is helpful to adopt a representation-theoretic viewpoint, and not to insist
on $L^{2}\left(  \mathbb{R}\right)  $ as the Hilbert space for the eigenvalue
problem (\ref{eqIntr.9}). We will look for abstract Hilbert spaces
$\mathcal{H}$ which admit solutions $\varphi\in\mathcal{H}$, $\varphi\neq0$,
to $M\varphi=\varphi$ in a way that naturally generalizes (\ref{eqIntr.9}).
This will also lead to results on multiresolutions which give solutions
\emph{up to unitary equivalence,} as well as conditions for equivalence. This
approach dictates another slight modification: if $\mathcal{H}$ is given only
abstractly, we must specify a unitary operator $U\colon\mathcal{H}%
\rightarrow\mathcal{H}$ which corresponds to the scaling operator%
\begin{equation}
U\colon\psi\longmapsto\frac{1}{\sqrt{N}\mathstrut}\psi\left(  \frac{x}%
{N}\right)  \label{eqIntr.13}%
\end{equation}
for the special case when $\mathcal{H}=L^{2}\left(  \mathbb{R}\right)  $.
Similarly, we must specify a representation $\pi$ of $L^{\infty}\left(
\mathbb{T}\right)  $ on $\mathcal{H}$ such that
\begin{equation}
U\pi\left(  f\right)  =\pi\left(  f\left(  z^{N}\right)  \right)  U,\qquad
f\in L^{\infty}\left(  \mathbb{T}\right)  \label{eqIntr.14}%
\end{equation}
as a commutation relation for operators on $\mathcal{H}$. In this wider
setting, the problem (\ref{eqIntr.9}) then takes the form:%
\begin{equation}
\varphi\in\mathcal{H}\colon\qquad U\varphi=\pi\left(  m_{0}\right)  \varphi.
\label{eqIntr.15}%
\end{equation}
This means that $\left(  U,\pi,\mathcal{H},\varphi\right)  $ are specified,
and satisfy (\ref{eqIntr.14})--(\ref{eqIntr.15}). The symbol $U$ denotes both
an element in $\mathfrak{A}_{N}$
(i.e., the $C^{\ast }$-algebra on two
unitary generators $U$ and $V$ subject
to $UVU^{-1}=V^{N}$),
and a unitary operator in $\mathcal{H}$. If
$\tilde{\pi}\in\operatorname*{Rep}\left(  \mathfrak{A}_{N},\mathcal{H}\right)  $
is the corresponding representation, then $\tilde{\pi}\left(  U\right)  =U$
(with the double meaning for $U$).

The restricting condition we place on $\left(  U,\pi,L^{2}\left(
\mathbb{R}\right)  \right)  $ is that there is some \linebreak $\varphi\in
L^{2}\left(  \mathbb{R}\right)  $, which satisfies $U\varphi=\pi\left(
m_{0}\right)  \varphi$, and $\hat{\varphi}\left(  0\right)  =1$, where
$\hat{\varphi}$ is the Fourier transform, i.e.,
\begin{equation}
\hat{\varphi}\left(  \omega\right)  :=\int_{\mathbb{R}}e^{-i\omega x}%
\varphi\left(  x\right)  \,dx.\label{eqax+b.8}%
\end{equation}
A system $\left(  U,\pi,L^{2}\left(  \mathbb{R}\right)  ,\varphi,m_{0}\right)
$ with these properties will be called a \emph{wavelet representation,} and
$\varphi$ scaling function \textup{(}or \emph{father function\/}\textup{).}
They are used in the construction of wavelets via \emph{multiresolutions;} see
\cite{BrJo97b} and \cite{Dau92}.

If a wavelet representation is given with scaling function, we form the
$L^{1}\left(  \mathbb{T}\right)  $-function $h_{\varphi}$ as follows:
$\mathbb{T}=\mathbb{R}\diagup2\pi\mathbb{Z}$, $z=e^{-i\omega}\in\mathbb{T}$ a
parametrization, and
\[
h_{\varphi}\left(  z\right)  =\sum_{n\in\mathbb{Z}}z^{n}\ip{\pi\left
( e_{n}\right) \varphi}{\varphi}_{L^{2}\left(  \mathbb{R}\right)  }=\frac
{1}{2\pi}\sum_{k\in\mathbb{Z}}\left|  \hat{\varphi}\left(  \omega+2\pi
k\right)  \right|  ^{2}.
\]
If the \emph{Ruelle transfer operator} is defined as%
\begin{equation}
\left(  Rf\right)  \left(  z\right)  :=\frac{1}{N}\sum_{\substack{w\in
\mathbb{T}\\w^{N}=z}}\left|  m_{0}\left(  w\right)  \right|  ^{2}f\left(
w\right)  , \label{eqax+b.9}%
\end{equation}
then%
\begin{equation}
R\left(  h_{\varphi}\right)  =h_{\varphi}. \label{eqax+b.10}%
\end{equation}
(See Lemma \ref{LemPoof.3}.)

In this paper, we prove a converse to this result: Every solution $h$ to
$Rh=h$ arises this way as $h=h_{\varphi_{\pi}}$ for some representation $\pi$.
(See Theorems \ref{Thmax+b.3} and \ref{ThmKean.6}.)

We further show in Chapter \ref{R-ha} that the solutions to (\ref{eqIntr.15})
may be represented by a family of finite (positive) Borel measures $\nu
=\nu\left(  m_{0},h\right)  $ on a certain Bohr compactification $K_{N}$ of
$\mathbb{R}$. It is simply the Bohr--Besicovitch compactification
(see \cite{Bes55})
corresponding to frequencies of the form
\begin{equation}
\sum_{i}\frac{n_{i}}{N^{i}},\qquad n_{i}\in\mathbb{Z}\text{, finite sums.}
\label{eqIntr.16}%
\end{equation}
Hence, in this case, $\mathcal{H}=L^{2}\left(  K_{N},\nu\right)  $. But, more
surprisingly, $\nu$ may be chosen (depending on $m_{0}$, $h$ as in
(\ref{eqIntr.8})) such that the solution $\varphi$ to (\ref{eqIntr.15}), in
the $L^{2}\left(  K_{N},\nu\right)  $ representation, is simply $\varphi
=\openone$, i.e., the constant unit function in $L^{2}\left(  K_{N}%
,\nu\right)  $.

A second advantage of the present wider scope is that it includes applications
outside wavelet theory. We note in Chapter \ref{Kean} that Ruelle operators of
the form (\ref{eqIntr.5}) arise naturally (in fact first!) in statistical mechanics
\cite{Rue68}
and ergodic theory
\cite{Kea72},
and our representation-theoretic results are
applied there.

Our motivation for the study of the interplay between the two operators $M$
and $R$ in (\ref{eqIntr.1}) and (\ref{eqIntr.5}) came from the now familiar
connection between the spectral analysis of $R$ and the convergence properties
of the iterations%
\begin{equation}
\psi^{\left(  0\right)  },\;M\psi^{\left(  0\right)  },\;M^{2}\psi^{\left(
0\right)  },\;\dots\label{eqIntr.17}%
\end{equation}
when $\psi^{\left(  0\right)  }$ is given. Similarly the spectral theory of
$R$ dictates directly the regularity properties of the solution $\varphi$ to
(\ref{eqIntr.9}). (See \cite{Dau92}, \cite{CoDa96}, \cite{BEJ97},
\cite{Vil94}, and Chapter \ref{Wave} below
for more
details on that point.)

Until recently, there was only very little in the literature on notions of
\emph{equivalence} for wavelets, or for the filter functions which are used to
generate them, or even \emph{isomorphism,} or \emph{invariants,} for these
objects; see, however, \cite{BrJo97b}.

A second motivation centers around the cocycle equivalence problems for
wavelet filters. Let $m_{0}^{{}}$ and $m_{0}^{\prime}$ be in $L^{\infty
}\left(  \mathbb{T}\right)  $. We say that they are \emph{cocycle equivalent}
if there is an $f\in L^{2}\left(  \mathbb{T}\right)  $ such that
\begin{enumerate}
\item \label{Intr(1)}$f\left(  z^{N}\right)  m_{0}^{\prime}\left(  z\right)
=m_{0}^{{}}\left(  z\right)  f\left(  z\right)  $ and

\item \label{Intr(2)}$\sum_{w^{N}=z}\left|  m_{0}^{\prime}\left(  w\right)
\right|  ^{2}=N$ for a.a.\ $z\in\mathbb{T}$.
\end{enumerate}
We stress that condition (\ref{Intr(2)}) is only assumed to hold for
$m_{0}^{\prime}$ and not for $m_{0}^{{}}$. Let $m_{0}^{{}}$ be fixed, and let
$R$ be the corresponding Ruelle operator of order $N$, i.e.,
\begin{equation}
\left(  Rf\right)  \left(  z\right)  :=\frac{1}{N}\sum_{w^{N}=z}\left|
m_{0}\left(  w\right)  \right|  ^{2}f\left(  w\right)  . \label{eqIntr.18}%
\end{equation}
In Chapter \ref{Cocy}, we show the following result: If $f\in L^{2}\left(
\mathbb{T}\right)  $ defines a cocycle equivalence, then $h\left(  z\right)
:=\left|  f\left(  z\right)  \right|  ^{2}$ solves $R\left(  h\right)  =h$. If
conversely $h\in L^{1}\left(  \mathbb{T}\right)  $, $h\geq0$, is given to
satisfy $R\left(  h\right)  =h$, then $f\left(  z\right)  :=h\left(  z\right)
^{\frac{1}{2}}$ defines a cocycle equivalence for some $m_{0}^{\prime}$, i.e.,
(\ref{Intr(1)})--(\ref{Intr(2)}) hold.

The relation $UfU^{-1}=f\left(  z^{N}\right)  $ for the abelian algebra
$L^{\infty}\left(  \mathbb{T}\right)  $ and a single unitary $U$ may be
rewritten as $Uf=f\left(  z^{N}\right)  U$, and then $U$ may possibly not be
unitary. A representation of the more general relation is then a pair $\left(
\pi,U\right)  $ where $\pi$ is a representation of $L^{\infty}\left(
\mathbb{T}\right)  $ in some Hilbert space $\mathcal{H}$, and $U$ is a bounded
operator in $\mathcal{H}$ such that%
\begin{equation}
U\pi\left(  f\right)  =\pi\left(  f\left(  z^{N}\right)  \right)  U
\label{eqIntr.19}%
\end{equation}
holds as an identity on $\mathcal{H}$ for all $f\in L^{\infty}\left(
\mathbb{T}\right)  $. Let $m_{0}\in L^{\infty}\left(  \mathbb{T}\right)  $,
and let $h\in L^{1}\left(  \mathbb{T}\right)  $ be given such that $h\geq0$
and $R_{m_{0}}\left(  h\right)  =h$. Then a representation $\left(
\pi,U\right)  $ results as follows: take $\mathcal{H}=L^{2}\left(
\mathbb{T},h\,d\mu\right)  $ and $U=S_{0}$ given by%
\begin{equation}
\left(  S_{0}\xi\right)  \left(  z\right)  =m_{0}\left(  z\right)  \xi\left(
z^{N}\right)  ,\qquad\xi\in L^{2}\left(  h\right)  :=L^{2}\left(
\mathbb{T},h\,d\mu\right)  , \label{eqIntr.20}%
\end{equation}
i.e., $\int_{\mathbb{T}}\left|  \xi\left(  z\right)  \right|  ^{2}h\left(
z\right)  \,d\mu\left(  z\right)  <\infty$ where $\mu$ is the normalized Haar
measure on $\mathbb{T}$, and%
\begin{equation}
\pi_{0}\left(  f\right)  \xi\left(  z\right)  =f\left(  z\right)  \xi\left(
z\right)  ,\qquad f\in L^{\infty}\left(  \mathbb{T}\right)  ,\;\xi\in
L^{2}\left(  h\right)  . \label{eqIntr.21}%
\end{equation}
We show in Theorem \ref{ThmCocy.5} that, if $m_{0}$ is a \emph{non-singular}
(defined below) wavelet filter, then $S_{0}$ (of (\ref{eqIntr.20})) is a
\emph{pure shift,} i.e., it is isometric in $L^{2}\left(  h\right)  $, and%
\begin{equation}
\bigcap_{n=1}^{\infty}S_{0}^{n}\left(  L^{2}\left(  h\right)  \right)
=\left\{  0\right\}  . \label{eqIntr.22}%
\end{equation}
In fact, for the wavelet filters, we show in Chapter \ref{Poof} that
$L^{2}\left(  h\right)  $ embeds isometrically into $L^{2}\left(
\mathbb{R}\right)  $; that is, there is an \emph{intertwining isometry}
$L^{2}\left(  h\right)  \overset{W}{\hooklongrightarrow}L^{2}\left(
\mathbb{R}\right)  $ such that%
\begin{align}
WS_{0}  &  =UW\label{eqIntr.23}\\%
\intertext{where}%
U\psi\left(  x\right)   &  :=N^{-\frac{1}{2}}\psi\left(  \frac{x}{N}\right)
,\qquad\psi\in L^{2}\left(  \mathbb{R}\right)  ,\label{eqIntr.24}\\%
\intertext{and}%
W\left(  \openone\right)   &  =\varphi. \label{eqIntr.25}%
\end{align}
In that case, $\varphi$ (in $L^{2}\left(  \mathbb{R}\right)  $) is the scaling
function, and%
\begin{equation}
\varphi=M\varphi\label{eqIntr.26}%
\end{equation}
with $M=U^{-1}\pi\left(  m_{0}\right)  $. Moreover, $M$ satisfies
\begin{equation}
M^{\ast}\pi\left(  f\right)  M=\pi\left(  \left|  m_{0}\right|  ^{2}f\left(
z^{N}\right)  \right)  ,\qquad f\in L^{\infty}\left(  \mathbb{T}\right)  .
\label{eqIntr.27}%
\end{equation}
More generally, if $\left(  M,\pi\right)  $ is given in a Hilbert space
$\mathcal{H}$, subject to (\ref{eqIntr.27}), we say that $M$ is a $\pi
$\emph{-isometry,} and the previous paper \cite{Jor98a} gives a complete
structure result for $\pi$-isometries, starting with the following Wold
decomposition:%
\begin{equation}
\mathcal{H}=\sideset{}{^{\smash{\oplus}}}{\sum}\limits_{n=0}^{\infty}\left[
M^{k}\mathcal{L}\right]  \oplus\bigcap_{1}^{\infty}\left[  M^{n}%
\mathcal{H}\right]  , \label{eqIntr.28}%
\end{equation}
where $\mathcal{L}:=\ker\left(  M^{\ast}\right)  $ and $\left[  M^{k}%
\mathcal{L}\right]  $ denotes the closure in $\mathcal{H}$ of $\left\{
M^{k}l\mid l\in\mathcal{L}\right\}  $. Moreover, the components in the
decomposition are mutually orthogonal, and each one reduces the representation
$\pi$.

For the case of non-singular wavelet filters, it follows from (\ref{eqIntr.22}%
) that $U$ is then unitarily equivalent to the \emph{bilateral shift} which
extends $S_{0}$. (Recall $S_{0}$ is then a \emph{unilateral} shift.) If the
wavelet filter is singular, i.e., if $m_{0}$ vanishes on a subset of
$\mathbb{T}$ of positive measure, then $\mathcal{L}=\ker M^{\ast}\neq\left\{
0\right\}  $, and the decomposition (\ref{eqIntr.28}) then has the shift part%
\[
\left[  M^{k}\mathcal{L}\right]  \overset{M}{\longrightarrow}\left[
M^{k+1}\mathcal{L}\right]  .
\]
But the scaling function $\varphi$ must be in $\bigcap_{1}^{\infty}\left[
M^{n}\mathcal{H}\right]  $. For representations more general than wavelets, we
show in Chapter \ref{R-ha} that there is an \emph{intertwining isometry}
$W_{B}$ (analogous to $W$ above):%
\begin{equation}
L^{2}\left(  h\right)  \overset{W_{B}}{\hooklongrightarrow}L^{2}\left(
K_{N},\nu\right)  \label{eqIntr.29}%
\end{equation}
where $K_{N}$ is the compact Bohr group, and $\nu$ is a Borel probablilty
measure on $K_{N}$ which depends on $m_{0}$ and $h$. The relations which
correspond to (\ref{eqIntr.23})--(\ref{eqIntr.26}) are:%
\begin{align}
W_{B}S_{0}  &  =UW_{B},\label{eqIntr.30}\\
\left(  U\psi\right)  \left(  \chi\right)   &  =m_{0}\left(  \chi\right)
\psi\left(  \chi^{N}\right)  ,\qquad\psi\in L^{2}\left(  K_{N},\nu\right)
,\;\chi\in K_{N},\label{eqIntr.31}\\%
\intertext{where $K_{N}=\left( \mathbb{Z}\left[ \frac{1}{N}\right
] \right) \sphat{}$, $\chi^{n}\left( \lambda
\right) =\chi\left( n\lambda\right) $, $\lambda\in\mathbb{Z}\left[ \frac{1}%
{N}\right] $, $n\in\mathbb{Z}$, and}%
m_{0}\left(  \chi\right)   &  :=\sum_{n\in\mathbb{Z}}a_{n}\chi^{n}
\label{eqIntr.32}%
\end{align}
is the natural extension of $m_{0}$ from $\mathbb{T}$ to $K_{N}$.

But, in this case, $W_{B}$ sends the constant function in $\mathbb{T}$ to that
of $K_{N}$, i.e.,
\begin{equation}
W_{B}\left(  \openone\right)  =\openone, \label{eqIntr.33}%
\end{equation}
so that $\varphi$ is now represented by the constant function $\openone$ in
$L^{2}\left(  K_{N},\nu\right)  $, i.e.,%
\begin{equation}
M\left(  \openone\right)  =\openone\label{eqIntr.34}%
\end{equation}
with $M=U^{-1}\pi_{0}\left(  m_{0}\right)  $ now taking the form%
\begin{equation}
M\psi\left(  \chi\right)  =\psi\left(  \chi^{-N}\right)  ,\qquad\psi\in
L^{2}\left(  K_{N},\nu\right)  ,\;\chi\in K_{N}. \label{eqIntr.35}%
\end{equation}

In the next two chapters, we introduce a certain discrete solvable group
$G_{N}$ and its $C^{\ast}$-algebra $\mathfrak{A}_{N}$; and we prove the
representation theorem alluded to above. The scaling function is formulated in
an abstract setting, and we identify a corresponding class of representations
which are defined from solutions $h$ to $Rh=h$, $h\in L^{1}$, $h\geq0$.
Conversely, we show that every such eigenfunction $h$ defines a representation
with an (abstract) scaling vector.

One of our motivations for the
analysis of (\ref{eqIntr.1}) or (\ref{eqIntr.9}) was recent
work on structural properties of the solutions 
$\varphi $, and aimed at giving new invariants
for them. Our own papers
\cite{BEJ97} and \cite{Jor98a} establish such
representation-theoretic invariants.
Analysis of (\ref{eqIntr.1}) in a variety of guises can also be found in
\cite{CDM91}, \cite{ChDe60}, \cite{CoDa96}, 
\cite{CoRa90}, \cite{Ho96}, \cite{Her95}, \cite{JLS98}, 
\cite{LMW96},  and \cite{LWC95}.

\section{\label{ax+b}A discrete $ax+b$ group}

Let $N\in\left\{  2,3,\dots\right\}  $ be given, and let $\Lambda
_{N}:=\mathbb{Z}\left[  \frac{1}{N}\right]  $ be the ring obtained from
$\mathbb{Z}$ by extending with the fraction $\frac{1}{N}$, i.e., $\Lambda_{N}$
contains $\mathbb{Z}$ and all powers $\left\{  \frac{1}{N^{k\mathstrut}}\mid
k=1,2,\dots\right\}  $. We will then consider the group $G=G_{N}$ of all
$2\times2$ matrices%
\begin{equation}
\left\{  \left(
\begin{smallmatrix}
N^{j} & \lambda\\
0 & 1
\end{smallmatrix}
\right)  \mid j\in\mathbb{Z},\;\lambda\in\Lambda_{N}\right\}  .
\label{eqax+b.1}%
\end{equation}
We showed in \cite{BreJor} that there is a one-to-one correspondence between
the unitary representations of $G_{N}$ and the $\ast$-representations of the
$C^{\ast}$-algebra $\mathfrak{A}$ on two unitary generators $U$, $V$, subject to
the relations
\begin{equation}
UVU^{-1}=V^{N}. \label{eqax+b.2}%
\end{equation}
By a representation of (\ref{eqax+b.2}), we mean a realization of $U$ and $V$
as unitary operators on some Hilbert space $\mathcal{H}$, say, such that
(\ref{eqax+b.2}) also holds for those operators. Let $f\in L^{\infty}\left(
\mathbb{T}\right)  $; then $f\left(  V\right)  $ is defined by the spectral
theorem, applied to $V$, and $\pi_{V}\left(  f\right)  :=f\left(  V\right)  $
is a representation of $L^{\infty}\left(  \mathbb{T}\right)  $ in the sense
that $\pi_{V}\left(  f_{1}f_{2}\right)  =\pi_{V}\left(  f_{1}\right)  \pi
_{V}\left(  f_{2}\right)  $, and $\pi_{V}\left(  f\right)  ^{\ast}=\pi
_{V}\left(  \bar{f}\right)  $ where $\bar{f}\left(  z\right)  :=\overline
{f\left(  z\right)  }$, $z\in\mathbb{T}$, $f_{1},f_{2},f\in L^{\infty}\left(
\mathbb{T}\right)  $. Moreover, (\ref{eqax+b.2}) then takes the form%
\begin{equation}
U\pi_{V}\left(  f\right)  U^{-1}=\pi_{V}\left(  f\left(  z^{N}\right)
\right)  , \label{eqax+b.3}%
\end{equation}
and conversely, every pair $\left(  \pi,U\right)  $ where $\pi$ is a
representation of $L^{\infty}\left(  \mathbb{T}\right)  $ on $\mathcal{H}$,
and $U$ is a unitary operator on $\mathcal{H}$ such that%
\begin{equation}
U\pi\left(  f\right)  U^{-1}=\pi\left(  f\left(  z^{N}\right)  \right)  ,
\label{eqax+b.4}%
\end{equation}
is of this form for some $V$. In fact, let $e_{n}\left(  z\right)  =z^{n}$,
$n\in\mathbb{Z}$, and set $V:=\pi\left(  e_{1}\right)  $.

Since $V$ is unitary, it has a spectral resolution $V=\int_{\mathbb{T}}%
\lambda\,E\left(  d\lambda\right)  $ with a projection-valued spectral measure
$E\left(  \,\cdot\,\right)  $ on $\mathbb{T}$. We will study cyclic vectors
$\varphi\in\mathcal{H}$. If the measure $\left\|  E\left(  \,\cdot\,\right)
\varphi\right\|  ^{2}$ on $\mathbb{T}$ is absolutely continuous with respect
to Haar measure on $\mathbb{T}$, we say that the corresponding representation
is \emph{normal.} The normal representations will be denoted
$\operatorname*{Rep}\left(  \mathfrak{A}_{N},\mathcal{H}\right)  $.

The $C^{\ast}$-algebra on the relations (\ref{eqax+b.2}),
introduced in \cite{BreJor},
will be denoted
$\mathfrak{A}_{N}$, and we shall always use the correspondence between the
representation $\operatorname*{Rep}\left(  \mathfrak{A}_{N},\mathcal{H}\right)  $
on some Hilbert space $\mathcal{H}$, and the corresponding unitary
representations of $G_{N}$.

We also showed in \cite{BreJor} that the
discrete (solvable) group $G_{N}$ has
representations which are not
predicted from Mackey's theory
of semidirect products of
\emph{continuous} groups. In fact,
the simplest discrete group
constructions lead to type $\mathrm{III}$
representations, even in cases
where the analogous continuous
groups have type $\mathrm{I}$ representations;
see also \cite{Bla77} and \cite{BoCo95}.

A state $\sigma$ on $\mathfrak{A}_{N}$ is said to be \emph{normal} if there is an
$L^{1}\left(  \mathbb{T}\right)  $-function $h$, $h\geq0$, such that%
\[
\sigma\left(  f\right)  =\int_{\mathbb{T}}fh\,d\mu,\qquad f\in L^{\infty
}\left(  \mathbb{T}\right)  ,
\]
where $\mu$ denotes the normalized Haar measure on $\mathbb{T}$, and we view
$\mathfrak{A}_{N}$ as the $C^{\ast}$-algebra on generators $f\in L^{\infty}\left(
\mathbb{T}\right)  $ and a single unitary element $U$, such that%
\[
UfU^{-1}=f\left(  z^{N}\right)  .
\]
We say that the state $\sigma$ is $U$\emph{-invariant} if%
\[
\sigma\left(  UAU^{-1}\right)  =\sigma\left(  A\right)  ,\qquad A\in
\mathfrak{A}_{N}.
\]
If $\sigma$ is normal, then $U$-invariance is equivalent to the condition%
\begin{equation}
\int_{\mathbb{T}}fh\,d\mu=\int_{\mathbb{T}}f\left(  z^{N}\right)  h\left(
z\right)  \,d\mu\left(  z\right)  ,\qquad f\in L^{\infty}\left(
\mathbb{T}\right)  . \label{eqax+b.4bis}%
\end{equation}

\begin{lemma}
\label{Lemax+b.2bis}A normal state $\sigma$ on $\mathfrak{A}_{N}$ with density $h$
\textup{(}$\in L^{1}\left(  \mathbb{T}\right)  $\textup{)} is $U$-invariant if
and only if $h$ is the constant function.
\end{lemma}

\begin{proof}
One direction is immediate, so let $h\in L^{1}\left(  \mathbb{T}\right)  $,
$h\geq0$, be a density for a fixed state $\sigma$ on $\mathfrak{A}_{N}$, and
assume $U$-invariance. Then apply (\ref{eqax+b.4bis}) to $f\left(  z\right)
=e_{n}\left(  z\right)  =z^{n}$, $z\in\mathbb{T}$, $n\in\mathbb{Z}$. Let
$\tilde{h}\left(  n\right)  =\int_{\mathbb{T}}\bar{e}_{n}h\,d\mu$ be the
Fourier coefficients of $h$. We get%
\[
\tilde{h}\left(  n\right)  =\tilde{h}\left(  Nn\right)  ,\qquad n\in
\mathbb{Z}.
\]
The operator%
\[
\left(  R_{N}f\right)  \left(  z\right)  =\frac{1}{N}\sum_{w^{N}=z}f\left(
w\right)
\]
satisfies%
\[
\left(  R_{N}f\right)  \sptilde\left(  n\right)  =\tilde{f}\left(  Nn\right)
,\qquad n\in\mathbb{Z},
\]
and a standard argument yields%
\begin{equation}
R_{N}^{n}\left(  f\right)  \underset{n\rightarrow\infty}{\longrightarrow}%
\int_{\mathbb{T}}f\,d\mu; \label{eqax+b.4ter}%
\end{equation}
see \cite{BrJo97b} for details. If the state $\sigma$ is $U$-invariant, $h$
must therefore satisfy $R_{N}\left(  h\right)  =h$, and by (\ref{eqax+b.4ter}%
), $h$ must be a constant function, which completes the proof.
\end{proof}

Our interest in $\operatorname*{Rep}\left(  \mathfrak{A}_{N},\mathcal{H}\right)  $
started with the following example:
\begin{equation}
\mathcal{H}=L^{2}\left(  \mathbb{R}\right)  ,\qquad\left(  U\psi\right)
\left(  x\right)  =\frac{1}{\sqrt{N}\mathstrut}\psi\left(  \frac{x}{N}\right)
,\label{eqax+b.5}%
\end{equation}
and%
\begin{equation}
\pi\left(  e_{n}\right)  \psi\left(  x\right)  :=\psi\left(  x-n\right)
,\qquad\psi\in L^{2}\left(  \mathbb{R}\right)  ,\;x\in\mathbb{R}%
,\;n\in\mathbb{Z}.\label{eqax+b.6}%
\end{equation}

\begin{definition}
\label{Defax+b.1}Let $m_{0}\in L^{\infty}\left(  \mathbb{T}\right)  $ be
given, and assume the following three properties:
\begin{enumerate}
\item \label{Defax+b.1(1)}%
\begin{equation}
\sum_{k=0}^{N-1}\left|  m_{0}\left(  e^{i\frac{2\pi k}{N}}z\right)  \right|
^{2}=N, \label{eqax+b.7}%
\end{equation}

\item \label{Defax+b.1(2)}$m_{0}$ is continuous on $\mathbb{T}$ near $z=1$, and

\item \label{Defax+b.1(3)}$m_{0}\left(  1\right)  =\sqrt{N}$.
\end{enumerate}
These functions are called \emph{low-pass filters} and are central
to the theory of orthogonal wavelets.
\end{definition}

The transfer operator $R$ of Ruelle
(\ref{eqIntr.18})
plays a crucial r\^{o}le in the study of
the regularity properties of the wavelets which are defined from some given
filter $m_{0}$. (See, e.g., \cite{CoDa96}.) When $m_{0}$ is given, subject to
(\ref{Defax+b.1(1)})--(\ref{Defax+b.1(3)}) above, there is an algorithm due to
S.G. Mallat \cite{Mal89,Dau92} which exhibits $\varphi$ as the inverse Fourier
transform of the infinite product%
\begin{equation}
\lim_{n\rightarrow\infty}\prod_{k=1}^{n}\frac{m_{0}\left(  \frac{\omega
}{N^{k\mathstrut}}\right)  }{\sqrt{N}\mathstrut}\chi_{\left[  -\pi N^{n},\pi
N^{n}\right\rangle }^{{}}\left(  \omega\right)  .\label{eqax+b.11}%
\end{equation}
The limit is known to be well defined, and specifying an $L^{2}\left(
\mathbb{R}\right)  $-function $F\left(  \omega\right)  $. Moreover
$\varphi\left(  x\right)  =\frac{1}{2\pi}\int_{\mathbb{R}}e^{i\omega
x}F\left(  \omega\right)  \,d\omega$ will then be in $L^{2}\left(
\mathbb{R}\right)  $ and satisfy%
\begin{equation}
U\varphi=\pi\left(  m_{0}\right)  \varphi,\label{eqax+b.12}%
\end{equation}
$\varphi$ depending on the representation. Let $m_{0}\left(  z\right)
=\sum_{n\in\mathbb{Z}}a_{n}z^{n}$. Then recall the right-hand side is%
\[
\sum_{n\in\mathbb{Z}}a_{n}\varphi\left(  x-n\right)  .
\]
Equivalently, the scaling identity (\ref{eqax+b.12}) may be rephrased as%
\begin{equation}
\varphi\left(  x\right)  =\sqrt{N}\sum_{n\in\mathbb{Z}}a_{n}\varphi\left(
Nx-n\right)  ,\qquad x\in\mathbb{R}.\label{eqax+b.13}%
\end{equation}

\begin{lemma}
\label{Lemax+b.2}The wavelet representation $\left(  L^{2}\left(
\mathbb{R}\right)  ,\pi,U\right)  $ which is given by a wavelet filter $m_{0}$
and a scaling function $\varphi$ is irreducible.
\end{lemma}

\begin{proof}
Recall $m_{0}$ satisfies conditions (\ref{Defax+b.1(1)})--(\ref{Defax+b.1(3)})
of Definition \ref{Defax+b.1}, and the corresponding scaling function
$\varphi\in L^{2}\left(  \mathbb{R}\right)  $ is then determined by the Mallat
algorithm; see (\ref{eqax+b.11}). But (\ref{eqax+b.11}) also shows that any
other function $\varphi_{1}$, say, which satisfies $U\varphi_{1}=\pi\left(
m_{0}\right)  \varphi_{1}$, or equivalently%
\[
\varphi_{1}\left(  x\right)  =\sqrt{N}\sum_{n}a_{n}\varphi_{1}\left(
Nx-n\right)  ,
\]
where $m_{0}\left(  z\right)  =\sum_{n}a_{n}z^{n}$, must be a constant times
$\varphi$, i.e., $\varphi_{1}=c\varphi$. If $P$ is an operator on
$L^{2}\left(  \mathbb{R}\right)  $ which commutes with both $U$ and
$\pi\left(  L^{\infty}\left(  \mathbb{T}\right)  \right)  $, then
$MP\varphi=P\varphi$ where $M=U^{-1}\pi\left(  m_{0}\right)  $. It follows
that $P\varphi=c\varphi$ for some constant $c$. But $\varphi$ is a cyclic
vector for the von Neumann algebra $\mathfrak{A}$ generated by $U$ and $\pi\left(
L^{\infty}\left(  \mathbb{T}\right)  \right)  $, so%
\[
PA\varphi=AP\varphi=Ac\varphi=cA\varphi
\text{\qquad for all }A\in\mathfrak{A},
\]
and we conclude that $P$ is a scalar times the identity operator in
$L^{2}\left(  \mathbb{R}\right)  $, concluding the proof of irreducibility.
\end{proof}

\begin{theorem}
\label{Thmax+b.3}\ 

\begin{enumerate}
\item \label{Thmax+b.3(1)}Let $m_{0}\in L^{\infty}\left(  \mathbb{T}\right)
$, and suppose $m_{0}$ does not vanish on a subset of $\mathbb{T}$ of positive
measure. Let%
\begin{equation}
\left(  Rf\right)  \left(  z\right)  =\frac{1}{N}\sum_{w^{N}=z}\left|
m_{0}\left(  w\right)  \right|  ^{2}f\left(  w\right)  ,\qquad f\in
L^{1}\left(  \mathbb{T}\right)  . \label{eqax+b.14}%
\end{equation}
Then there is a one-to-one correspondence between the data
\textup{(\ref{Thmax+b.3(1)(1)})} and \textup{(\ref{Thmax+b.3(1)(2)})} below,
where \textup{(\ref{Thmax+b.3(1)(2)})} is understood as equivalence classes
under unitary equivalence:

\begin{enumerate}
\renewcommand{\theenumi}{\relax}

\item \label{Thmax+b.3(1)(1)}$h\in L^{1}\left(  \mathbb{T}\right)  $, $h\geq
0$, and%
\begin{equation}
R\left(  h\right)  =h. \label{eqax+b.15}%
\end{equation}

\item \label{Thmax+b.3(1)(2)}$\tilde{\pi}\in\operatorname*{Rep}\left(
\mathfrak{A}_{N},\mathcal{H}\right)  $, $\varphi\in\mathcal{H}$, and the unitary
$U$ from $\tilde{\pi}$ satisfying%
\begin{equation}
U\varphi=\pi\left(  m_{0}\right)  \varphi. \label{eqax+b.16}%
\end{equation}
\end{enumerate}

\item \label{Thmax+b.3(2)}From \textup{(\ref{Thmax+b.3(1)(1)})}$\rightarrow
$\textup{(\ref{Thmax+b.3(1)(2)}),} the correspondence is given by%
\begin{equation}
\ip{\varphi}{\pi\left( f\right) \varphi}_{\mathcal{H}}=\int_{\mathbb{T}%
}fh\,d\mu, \label{eqax+b.17}%
\end{equation}
where $\mu$ denotes the normalized Haar measure on $\mathbb{T}$.

From \textup{(\ref{Thmax+b.3(1)(2)})}$\rightarrow$%
\textup{(\ref{Thmax+b.3(1)(1)}),} the correspondence is given by%
\begin{equation}
h\left(  z\right)  =h_{\varphi}\left(  z\right)  =\sum_{n\in\mathbb{Z}}%
z^{n}\ip{\pi\left( e_{n}\right) \varphi}{\varphi}_{\mathcal{H}}.
\label{eqax+b.18}%
\end{equation}

\item \label{Thmax+b.3(3)}When \textup{(\ref{Thmax+b.3(1)(1)})} is given to
hold for some $h$, and $\tilde{\pi}\in\operatorname*{Rep}\left(  \mathfrak{A}%
_{N},\mathcal{H}\right)  $ is the corresponding cyclic representation with
$U\varphi=\pi\left(  m_{0}\right)  \varphi$, then the representation is unique
from $h$ and \textup{(\ref{eqax+b.17})} up to unitary equivalence: that is, if
$\pi^{\prime}\in\operatorname*{Rep}\left(  \mathfrak{A}_{N},\mathcal{H}^{\prime
}\right)  $, $\varphi^{\prime}\in\mathcal{H}^{\prime}$ also cyclic and
satisfying%
\begin{align*}
\ip{\varphi^{\prime}}{\pi^{\prime}\left( f\right) \varphi^{\prime}} &
=\int_{\mathbb{T}}fh\,d\mu\\%
\intertext{and}%
U^{\prime}\varphi^{\prime} &  =\pi^{\prime}\left(  m_{0}\right)
\varphi^{\prime},
\end{align*}
then there is a unitary isomorphism $W$ of $\mathcal{H}$ onto $\mathcal{H}%
^{\prime}$ such that $W\pi\left(  A\right)  =\pi^{\prime}\left(  A\right)  W$,
$A\in\mathfrak{A}_{N}$, and $W\varphi=\varphi^{\prime}$.
\end{enumerate}
\end{theorem}

In the setup for the theorem, we are \emph{not} assuming that $\frac{1}{N}%
\sum_{w^{N}=z}\left|  m_{0}\left(  w\right)  \right|  ^{2}=1$, although this
will be the case for the applications to wavelets. Hence the existence of
solutions to the eigenvalue problem%
\begin{equation}
Rf=f,\qquad f\in L^{1}\left(  \mathbb{T}\right)  ,\;f\geq0,\;f\neq0,
\label{eqax+b.pound}%
\end{equation}
is not guaranteed.

We note further that the mapping $z\mapsto z^{N}$ of $\mathbb{T}$ into
$\mathbb{T\simeq}\mathbb{R}\diagup\mathbb{Z}$ is a special case of the
following more general setup. We will show in Chapter \ref{Kean} that Theorem
\ref{Thmax+b.3} carries over to the more general setting.

Let $I=\left[  0,1\right]  $, and let $T\colon I\rightarrow I$ be a
\emph{piecewise
expanding }$C^{2}$\emph{ surjective Markov map,}
i.e., there is%
\[
0=x_{0}<x_{1}<\dots<x_{N-1}<x_{N}=1
\]
such that the restriction of $T$ to each of the subintervals is monotone.
Further there exists $\beta>1$ such that $\inf_{x\in I}\left|  T^{\prime
}\left(  x\right)  \right|  \geq\beta$; and the implication
\begin{equation}
\text{if }T\left(  \left(  x_{i-1},x_{i}\right)  \right)  
\cap\left(  x_{j-1},x_{j}\right)  \neq\varnothing
\text{, then }\left(  x_{j-1},x_{j}\right)  
\subset T\left(  \left(  x_{i-1},x_{i}\right)  \right)  
  \label{eqax+b.poundbis}
\end{equation}
holds. Set
\[
Rf\left(  x\right)  =\sum_{Ty=x}\frac{f\left(  y\right)  }{\left|  T^{\prime
}\left(  y\right)  \right|  }.
\]
It is easy to see that each solution to the eigenvalue problem%
\[
Rf=f,\qquad f\in L^{1}\left(  I\right)  ,\;f\geq0,
\]
defines a measure $f\,dx$ on $I$ which is $T$-invariant.

\begin{proposition}
\label{Proax+b.poundpound}\textup{(Pollicott--Yuri)} Let $T$ and $R$ be as
described. Then there exists $f\in L^{1}\left(  I\right)  $, $f\geq0$,
$f\neq0$ with $Rf=f$.
\end{proposition}

\begin{proof}
\cite[p.\ 127]{PoYu98}.
\end{proof}

\begin{remark}
\label{Remax+b.4}\textup{(}\emph{Moments of Representations.}\textup{)} An
element $\tilde{\pi}\in\operatorname*{Rep}\left(  \mathfrak{A}_{N},\mathcal{H}%
\right)  $ is generated by the operators $\left\{  \pi\left(  f\right)  \mid
f\in L^{\infty}\left(  \mathbb{T}\right)  \right\}  $, and the unitary
operator $U\colon\mathcal{H}\rightarrow\mathcal{H}$, and the commutation
relation is%
\[
U\pi\left(  f\right)  U^{-1}=\pi\left(  f\left(  z^{N}\right)  \right)  .
\]
The theorem is concerned with solutions $\varphi\in\mathcal{H}$ to
$U\varphi=\pi\left(  m_{0}\right)  \varphi$ when $m_{0}$ is given. For a given
representation, we have a \emph{spectral measure} $\nu_{\varphi}$ on
$\mathbb{T}$ such that%
\[
\ip{\varphi}{\pi\left( f\right) \varphi}=\int_{\mathbb{T}}f\left(  z\right)
\,d\nu_{\varphi}\left(  z\right)  ,
\]
and we noted that $R\left(  \nu_{\varphi}\right)  =\nu_{\varphi}$. So, if
$\nu_{\varphi}$ is absolutely continuous, with Radon--Nikodym derivative
$\frac{d\nu_{\varphi}}{d\mu}=h$, then $R\left(  h\right)  =h$, and%
\[
\ip{\varphi}{\pi\left( f\right) \varphi}=\int_{\mathbb{T}}fh\,d\mu,
\]
and we say that the right-hand side represents \emph{the moments} of $\pi$ in
the state $\varphi$; specifically,%
\begin{equation}
\ip{\varphi}{\pi\left( e_{n}\right) \varphi}=\int_{\mathbb{T}}z^{n}h\left(
z\right)  \,d\mu\left(  z\right)  ,\qquad n\in\mathbb{Z}, \label{eqax+b.19}%
\end{equation}
where $\mu$ as usual denotes the normalized Haar measure on $\mathbb{T}$.
\emph{The other moments} are%
\begin{equation}
\ip{\varphi}{U^{n}\varphi}=\int_{\mathbb{T}}m_{0}^{\left(  n\right)  }\left(
z\right)  h\left(  z\right)  \,d\mu\left(  z\right)  ,\qquad n=0,1,2,\dots,
\label{eqax+b.20}%
\end{equation}
where%
\[
m_{0}^{\left(  n\right)  }\left(  z\right)  :=m_{0}\left(  z\right)
m_{0}\left(  z^{N}\right)  \cdots m_{0}\left(  z^{N^{n-1}}\right)  .
\]

The mixed moments
\[
\omega _{\varphi }\left( U^{-k}fU^{n}\right) 
=\ip{\varphi }{U^{\ast \, k}\pi\left( f\right) U^{n}\varphi }
=\ip{U^{k}\varphi }{\pi\left( f\right) U^{n}\varphi }
\text{\qquad for }0\leq k\leq n
\]
involve the Ruelle operator $R$ through the
formula
\begin{equation}
\omega _{\varphi }\left( U^{-k}fU^{n}\right) 
=\int _{\mathbb{T}}m_{0}^{\left( n-k\right) }R^{k}\left( fh\right) \,d\mu.
\label{eqax+b.poundter}
\end{equation}

We show
that a cyclic representation, in this case $\tilde{\pi}%
\in\operatorname*{Rep}\left(  \mathfrak{A}_{N},\mathcal{H}\right)  $, is uniquely
determined by its moments
\textup{(\ref{eqax+b.19})} and \textup{(\ref{eqax+b.poundter});}
\emph{cf.,} e.g., \textup{(\ref{eqax+b.5})} and \textup{(\ref{eqax+b.6}).}
\end{remark}

While, in the statement of Theorem \ref{Thmax+b.3}, we are not assuming that
$m_{0}$ satisfy%
\[
\sum_{w^{N}=z}\left|  m_{0}\left(  w\right)  \right|  ^{2}=N,
\]
we then cannot be guaranteed solutions $h\in L^{1}\left(  \mathbb{T}\right)
$, $h\geq0$, $h\neq0$, to $R_{m_{0}}\left(  h\right)  =h$, where%
\[
R_{m_{0}}\left(  f\right)  \left(  z\right)  :=\frac{1}{N}\sum_{w^{N}%
=z}\left|  m_{0}\left(  w\right)  \right|  ^{2}f\left(  w\right)  .
\]
However, there is a recent such existence theorem due to L.~Herv\'{e}
\cite{Her95} with a direct wavelet application. It is assumed in \cite{Her95}
that $m_{0}\in C^{\infty}\left(  \mathbb{T}\right)  $, $m_{0}\left(  1\right)
=\sqrt{N}$. This means that the infinite product
\[
\prod_{k=1}^{\infty}\frac{m_{0}\left(  e^{-i\frac{\omega}{N^{k\mathstrut}}%
}\right)  }{\sqrt{N}\mathstrut}%
\]
is well defined pointwise as a function $f$ of $\omega\in\mathbb{R}$, and, of
course,%
\[
\sqrt{N}f\left(  \omega\right)  =m_{0}\left(  e^{-i\frac{\omega}{N}}\right)
f\left(  \frac{\omega}{N}\right)  ,\qquad\omega\in\mathbb{R}.
\]
For the wavelet problem, we need $f\in L^{2}\left(  \mathbb{R}\right)  $ such
that the inverse Fourier transform $\varphi=\check{f}$ may serve as an
$L^{2}\left(  \mathbb{R}\right)  $-scaling function. The theorem of Herv\'{e}
states that, under the given conditions, $f\in L^{2}\left(  \mathbb{R}\right)
$ if and only if there is a solution $h\in C^{\infty}\left(  \mathbb{T}%
\right)  $, $h\geq0$, $h\neq0$, to $R_{m_{0}}\left(  h\right)  =h$, in fact
$h\left(  1\right)  >0$.

The Ruelle operator, also called the Perron--Frobenius--Ruelle
operator, or the transfer operator,
is based on a simple but powerful idea.
In addition to the diverse applications given in
\cite{Rue76}, \cite{Rue78b}, \cite{Rue79}, \cite{Rue88}, and \cite{Rue90},
it has also found applications in ergodic theory
\cite{Sin72}, \cite{Wal75}, 
and harmonic analysis \cite{JoPe98b}, \cite{Jor98a}, \cite{Sch74},
and in statistical mechanics \cite{Rue68}, \cite{Mey98}.

\section{\label{Poof}Proof of Theorem \textup{\ref{Thmax+b.3}}}

The present chapter is entirely devoted to the proof of Theorem
\textup{\ref{Thmax+b.3}, and we begin with five lemmas.}

An alternative proof of Theorem \textup{\ref{Thmax+b.3} would be to get the
cyclic }representation (which is asserted in the theorem) from the
Gelfand--Naimark--Segal (GNS) construction. But then we would have to show
first that the data which are given in the theorem either define a positive
definite function on the $N$-adic $ax+b$ group, or alternatively a positive
linear functional (state) on $\mathfrak{A}_{N}$, and there is \emph{not} a direct
approach to doing that. It turns out to be shorter to first directly construct
the representation, and then, \emph{a posteriori,} to conclude the positive
definite properties of an associated function on the group, or a functional on
$\mathfrak{A}_{N}$. This is also discussed in detail in Chapter \ref{Kean}, which
in fact provides representations in a context which is more general than that
of Theorem \textup{\ref{Thmax+b.3}.}
Our general reference for oparator algebras
and the GNS construction
is \cite{BrRoI},
but our particular application in Chapter \ref{R-ha}
is closer to the viewpoint taken in \cite{GlJa87}.

\begin{lemma}
\label{LemPoof.1}If $m_{0}\in L^{\infty}\left(  \mathbb{T}\right)  $ then the
Ruelle operator in \textup{(\ref{eqax+b.9})} maps $L^{1}\left(  \mathbb{T}%
\right)  $ into itself, and has $L^{1}\rightarrow L^{1}$ operator norm equal
to $\left\|  m_{0}\right\|  _{\infty}^{2}$.
\end{lemma}

\begin{proof}
Let $f\in L^{1}\left(  \mathbb{T}\right)  $. Then%
\begin{align*}
\int_{\mathbb{T}}\left|  \left(  Rf\right)  \left(  z\right)  \right|
\,d\mu\left(  z\right)   &  \leq\frac{1}{N}\int_{\mathbb{T}}\sum_{w^{N}%
=z}\left|  m_{0}\left(  w\right)  \right|  ^{2}\left|  f\left(  w\right)
\right|  \,d\mu\left(  z\right) \\
&  =\int_{\mathbb{T}}\left|  m_{0}\left(  z\right)  \right|  ^{2}\left|
f\left(  z\right)  \right|  \,d\mu\left(  z\right) \\
&  \leq\left\|  m_{0}\right\|  _{\infty}^{2}\cdot\int_{\mathbb{T}}\left|
f\left(  z\right)  \right|  \,d\mu\left(  z\right)  .
\end{align*}
The $L^{1}\rightarrow L^{1}$ norm is in fact
\emph{equal to} $\left\| m_{0}\right\| _{\infty }^{2}$, and this is
based on an argument in
\cite{BrJo98b} to which we refer.
\end{proof}

\begin{lemma}
\label{LemPoof.2}Let $h\in L^{1}\left(  \mathbb{T}\right)  $, $h\geq0$, be
given, and let $L^{2}\left(  h\right)  $ denote the $L^{2}$-space of functions
on $\mathbb{T}$ defined relative to the absolutely continuous measure
$h\,d\mu$ \textup{(}where $\mu$ is Haar measure on $\mathbb{T}$\textup{),}
i.e.,
\begin{equation}
\left\|  f\right\|  _{h}^{2}:=\int_{\mathbb{T}}\left|  f\right|  ^{2}h\,d\mu.
\label{eqPoof.1}%
\end{equation}
On $L^{2}\left(  h\right)  $, we have the representation $\pi_{0}$, and the
operator $S_{0}$, defined as follows:%
\begin{equation}
\left(  \pi_{0}\left(  f\right)  \xi\right)  \left(  z\right)   :=f\left(
z\right)  \xi\left(  z\right)  ,\qquad f\in L^{\infty}\left(  \mathbb{T}%
\right)  ,\;\xi\in L^{2}\left(  h\right)  ,\label{eqPoof.2}
\end{equation}
and
\begin{equation}
\left(  S_{0}\xi\right)  \left(  z\right)   :=m_{0}\left(  z\right)
\xi\left(  z^{N}\right)  . \label{eqPoof.3}%
\end{equation}
Then%
\begin{equation}
S_{0}\pi_{0}\left(  f\right)  =\pi_{0}\left(  f\left(  z^{N}\right)  \right)
S_{0}; \label{eqPoof.4}%
\end{equation}
and $S_{0}$ is isometric in $L^{2}\left(  h\right)  $ if and only if $R\left(
h\right)  =h$.
\end{lemma}

\begin{proof}
The properties in the lemma are clear except for the criteria for $S_{0}$ to
be isometric: Let $h\in L^{1}\left(  \mathbb{T}\right)  $, $h\geq0$ be given,
and let $f_{1},f_{2}\in L^{\infty}\left(  \mathbb{T}\right)  $. Then%
\begin{align*}
\ip{S_{0}f_{1}}{S_{0}f_{2}}_{L^{2}\left(  h\right)  }  &  =\int_{\mathbb{T}%
}\left|  m_{0}\left(  z\right)  \right|  ^{2}\,\overline{f_{1}\left(
z^{N}\right)  }\,f_{2}\left(  z^{N}\right)  h\left(  z\right)  \,d\mu\left(
z\right) \\
&  =\frac{1}{N}\int_{\mathbb{T}}\overline{f_{1}\left(  z\right)  }%
\,f_{2}\left(  z\right)  \sum_{w^{N}=z}\left|  m_{0}\left(  w\right)  \right|
^{2}h\left(  w\right)  \,d\mu\left(  z\right) \\
&  =\int_{\mathbb{T}}\overline{f_{1}\left(  z\right)  }\,f_{2}\left(
z\right)  \left(  Rh\right)  \left(  z\right)  \,d\mu\left(  z\right)  ,
\end{align*}
and it follows that $S_{0}$ is $L^{2}\left(  h\right)  $-isometric if $Rh=h$.
But taking $f_{1}=e_{n_{1}}=z^{n_{1}}$, $f_{2}=e_{n_{2}}=z^{n_{2}}$,
$n_{1},n_{2}\in\mathbb{Z}$, we can see that the identity%
\[
\ip{S_{0}e_{n_{1}}}{S_{0}e_{n_{2}}}_{L^{2}\left(  h\right)  }=\ip{e_{n_{1}}%
}{e_{n_{2}}}_{L^{2}\left(  h\right)  }%
\]
implies that the two functions $R\left(  h\right)  $ and $h$ must have
identical Fourier coefficients. Since we have Fourier uniqueness for
$L^{1}\left(  \mathbb{T}\right)  $, the result follows, i.e., $Rh=h$ must hold
when it is given that $S_{0}$ is isometric on $L^{2}\left(  h\right)  $.
\end{proof}

\begin{lemma}
\label{LemPoof.3}Let $h\in L^{1}\left(  \mathbb{T}\right)  $, $h\geq0$, be
given, and let $\pi\in\operatorname*{Rep}\left(  \mathfrak{A}_{N},\mathcal{H}%
\right)  $, $\varphi\in\mathcal{H}$, satisfy%
\begin{equation}
\ip{\varphi}{\pi\left( f\right) \varphi}_{\mathcal{H}}=\int_{\mathbb{T}%
}fh\,d\mu,\qquad f\in L^{\infty}\left(  \mathbb{T}\right)  . \label{eqPoof.5}%
\end{equation}
Then%
\begin{equation}
h\left(  z\right)  :=\sum_{n\in\mathbb{Z}}z^{n}\ip{\pi\left( e_{n}%
\right) \varphi}{\varphi} \label{eqPoof.6}%
\end{equation}
is in $L^{1}\left(  \mathbb{T}\right)  $; and if further $U\varphi=\pi\left(
m_{0}\right)  \varphi$, then $Rh=h$.
\end{lemma}

\begin{proof}
Clearly, the expression (\ref{eqPoof.6}) makes sense as a distribution on
$\mathbb{T}$, and its Fourier coefficients are $n\mapsto\ip{\pi\left
( e_{n}\right) \varphi}{\varphi}_{\mathcal{H}}$. But substituting $f=e_{-n}$,
$n\in\mathbb{Z}$, into (\ref{eqPoof.5}) shows that $\ip{\pi\left( e_{n}%
\right) \varphi}{\varphi}=\int_{\mathbb{T}}\overline{e_{n}}\,h\,d\mu$, which
are the Fourier coefficients for the given $L^{1}\left(  \mathbb{T}\right)
$-function $h$. Hence the right-hand side of (\ref{eqPoof.6}) must be $h$,
again by Fourier uniqueness. If $U\varphi=\pi\left(  m_{0}\right)  \varphi$,
then we calculate the $h$-Fourier coefficients as follows:%
\begin{align*}
\tilde{h}\left(  n\right)   &  =\int_{\mathbb{T}}\overline{e_{n}}\,h\,d\mu\\
&  =\ip{\pi\left( e_{n}\right) \varphi}{\varphi}\\
&  =\ip{U\pi\left( e_{n}\right) \varphi}{U\varphi}\\
&  =\ip{\pi\left( e_{Nn}\right) U\varphi}{\pi\left( m_{0}\right)\varphi}\\
&  =\ip{\pi\left( e_{Nn}m_{0}\right) \varphi}{\pi\left( m_{0}\right)\varphi}\\
&  =\ip{\varphi}{\pi\left( e_{-Nn}\left| m_{0}\right| ^{2}\right) \varphi}\\
&  =\int_{\mathbb{T}}e_{-Nn}\left|  m_{0}\right|  ^{2}h\,d\mu\\
&  =\frac{1}{N}\int_{\mathbb{T}}e_{-n}\left(  z\right)  \sum_{w^{N}=z}\left|
m_{0}\left(  w\right)  \right|  ^{2}h\left(  w\right)  \,d\mu\left(  z\right)
\\
&  =\int_{\mathbb{T}}z^{-n}\left(  Rh\right)  \left(  z\right)  \,d\mu\left(
z\right) \\
&  =\left(  Rh\right)  \sptilde\left(  n\right)  ,\qquad n\in\mathbb{Z}.
\end{align*}
Hence the two $L^{1}\left(  \mathbb{T}\right)  $-functions $h$ and $Rh$ have
the same Fourier coefficients, and therefore $h=R\left(  h\right)  $ holds as claimed
\end{proof}

The correspondence (\ref{Thmax+b.3(1)(2)})$\rightarrow$(\ref{Thmax+b.3(1)(1)})
in Theorem \ref{Thmax+b.3} follows now directly from the lemmas, and we turn
to (\ref{Thmax+b.3(1)(1)})$\rightarrow$(\ref{Thmax+b.3(1)(2)}).

Let $h\in L^{1}\left(  \mathbb{T}\right)  $ be given satisfying $h\geq0$ and
$Rh=h$. The conditions on $m_{0}$ are just $m_{0}\in L^{\infty}\left(
\mathbb{T}\right)  $ and that $m_{0}$ does not vanish on a subset of
$\mathbb{T}$ of positive measure. In Lemma \ref{LemPoof.2}, we already did the
first step in a recursive algorithm for constructing the desired
representation $\tilde{\pi}\in\operatorname*{Rep}\left(  \mathfrak{A}%
_{N},\mathcal{H}\right)  $ and cyclic vector $\varphi\in\mathcal{H}$ such that
$U\varphi=\pi\left(  m_{0}\right)  \varphi$. This construction will be a
\emph{unitary} dilation (or lifting) of the \emph{isometric} properties
(\ref{eqPoof.1})--(\ref{eqPoof.3}) in Lemma \ref{LemPoof.2}. The construction
is also a generalized multiresolution. In the case $\mathcal{H}=L^{2}\left(
\mathbb{R}\right)  $, it is directly related to the traditional
multiresolution from wavelet theory, and we shall follow up on this in
Chapters \ref{Wave}--\ref{Cocy} below.

Since $h$ is given, we have the Hilbert space $L^{2}\left(  h\right)
=L^{2}\left(  \mathbb{T},h\,d\mu\right)  $ from (\ref{eqPoof.1}) in Lemma
\ref{LemPoof.2}, and we set $\mathcal{H}_{0}=L^{2}\left(  h\right)  $. It is
of course the completion of $L^{\infty}\left(  \mathbb{T}\right)  $ in the
norm $\left\|  \,\cdot\,\right\|  _{h}$ from (\ref{eqPoof.1}), and we write
$\mathcal{H}_{0}=\widetilde{\mathcal{V}}_{0}$ with%
\[
\mathcal{V}_{0}:=\left\{  \left(  \xi,0\right)  \mid\xi\in L^{\infty}\left(
\mathbb{T}\right)  \right\}  .
\]
Then on%
\[
\mathcal{V}_{n}:=\left\{  \left(  \xi,n\right)  \mid\xi\in L^{\infty}\left(
\mathbb{T}\right)  \right\}
\]
we set%
\begin{align}
\left\|  \left(  \xi,n\right)  \right\|  _{\mathcal{H}}^{2}  &  :=\int
_{\mathbb{T}}R^{n}\left(  \left|  \xi\right|  ^{2}h\right)  \,d\mu
\label{eqPoof.7}\\%
\intertext{and}%
\ip{\left( \xi,n\right) }{\left( \eta,n\right) }_{\mathcal{H}}  &
=\int_{\mathbb{T}}R^{n}\left(  \bar{\xi}\eta h\right)  \,d\mu\text{\qquad for
}n=1,2,\dots,\nonumber
\end{align}
where $\mu$ is the usual normalized Haar measure on $\mathbb{T}$, i.e.,
$\frac{1}{2\pi}\int_{-\pi}^{\pi}\cdots\,d\omega$ relative to $z=e^{-i\omega}$,
$\omega\in\mathbb{R}$, when functions on $\mathbb{T}$ are identified with
$2\pi$-periodic functions on $\mathbb{R}$. We now let $\mathcal{H}_{n}$ be the
completion of $\mathcal{V}_{n}$ in this norm, and we construct $\mathcal{H}$
itself as an inductive limit of the Hilbert spaces $\mathcal{H}_{n}$,
$n=0,1,2,\dots$. To do this, we construct a system of isometries
\begin{equation}
\begin{minipage}{72pt}\mbox{\includegraphics
[bbllx=0bp,bblly=-36bp,bburx=72bp,bbury=72bp,width=72pt]
{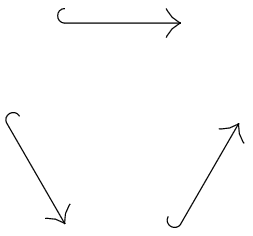}}\llap{\setlength{\unitlength}{32pt}\begin{picture}%
(2.25,3.375)(-1.125,-2.25)
\put(1.5,0.866){\makebox(0,0){\rlap
{$\displaystyle\mathcal{H}_{n+k+l}$}\phantom{$\displaystyle\mathcal{H}$}}}
\put(-1.5,0.866){\makebox(0,0){\rlap
{$\displaystyle\mathcal{H}_{n}$}\phantom{$\displaystyle\mathcal{H}$}}}
\put(0,-1.754){\makebox(0,0){\rlap
{$\displaystyle\mathcal{H}_{n+k}$}\phantom{$\displaystyle\mathcal{H}$}}}
\end{picture}}\end{minipage} \label{eqPoof.8}%
\end{equation}
and we get it from the completion of a corresponding isometry diagram
\[
\begin{minipage}{72pt}\mbox{\includegraphics
[bbllx=0bp,bblly=-36bp,bburx=72bp,bbury=72bp,width=72pt]
{funeq38.eps}%
}\llap{\setlength{\unitlength}{32pt}\begin{picture}(2.25,3.375)(-1.125,-2.25)
\put(1.5,0.866){\makebox(0,0){\rlap
{$\displaystyle\mathcal{V}_{n+k+l}$.}\phantom{$\displaystyle\mathcal{V}$}}}
\put(-1.5,0.866){\makebox(0,0){\rlap
{$\displaystyle\mathcal{V}_{n}$}\phantom{$\displaystyle\mathcal{V}$}}}
\put(0,-1.754){\makebox(0,0){\rlap
{$\displaystyle\mathcal{V}_{n+k}$}\phantom{$\displaystyle\mathcal{V}$}}}
\end{picture}}\end{minipage}%
\]
When $n$, $k$ are given, $n\geq0$, $k\geq1$, we construct the isometry
$\mathcal{V}_{n}\hookrightarrow\mathcal{V}_{n+k}$ by iteration of the one from
$\mathcal{V}_{n}$ to $\mathcal{V}_{n+1}$, i.e.,%
\[
\mathcal{V}_{n}\hooklongrightarrow\mathcal{V}_{n+1}\hooklongrightarrow
\mathcal{V}_{n+2}\hooklongrightarrow\cdots\hooklongrightarrow\mathcal{V}%
_{n+k},
\]
where $J\colon\mathcal{V}_{n}\rightarrow\mathcal{V}_{n+1}$ is defined by
\begin{equation}
J\left(  \left(  \xi,n\right)  \right)  :=\left(  \xi\left(  z^{N}\right)
,n+1\right)  . \label{eqPoof.9}%
\end{equation}
The isometric property of this operator is proved in the following lemma:

\begin{lemma}
\label{LemPoofIsomJ}The mapping $J$ defined in \textup{(\ref{eqPoof.9})} is
for each $n$ isometric from $\mathcal{H}_{n}$ into $\mathcal{H}_{n+1}$.
\end{lemma}

\begin{proof}
Let $\xi\in L^{\infty}\left(  \mathbb{T}\right)  $. Then
\begin{align*}
\left\|  \left(  \xi\left(  z^{N}\right)  ,n+1\right)  \right\|
_{\mathcal{H}}^{2}  &  =\int_{\mathbb{T}}R^{n+1}\left(  \left|  \xi\left(
z^{N}\right)  \right|  ^{2}h\left(  z\right)  \right)  \left(  z\right)
\,d\mu\left(  z\right) \\
& =\int_{\mathbb{T}}R^{n}\left(  R\left(  \left|  \xi\left(
z^{N}\right)  \right|  ^{2}h\left(  z\right)  \right)  \right)  
\,d\mu\left(  z\right) \\
&  =\int_{\mathbb{T}}R^{n}\left(  \left|  \xi\right|  ^{2}Rh\right)  \,d\mu\\
&  =\int_{\mathbb{T}}R^{n}\left(  \left|  \xi\right|  ^{2}h\right)  \,d\mu\\
&  =\left\|  \left(  \xi,n\right)  \right\|  _{\mathcal{H}}^{2},
\end{align*}
where we used the properties%
\begin{align*}
R\left(  f\left(  z^{N}\right)  g\left(  z\right)  \right)   &  =fR\left(
g\right) \\%
\intertext{and}%
Rh  &  =h
\end{align*}
for the Ruelle operator $R$ in (\ref{eqax+b.9}).
\end{proof}

We now define%
\begin{align}
U\left(  \xi,0\right)   &  :=\left(  S_{0}\xi,0\right)  =\left(  m_{0}\left(
z\right)  \xi\left(  z^{N}\right)  ,0\right)  ,\label{eqPoof.10}\\
U\left(  \xi,n+1\right)   &  :=\left(  m_{0}\left(  z^{N^{n}}\right)
\xi\left(  z\right)  ,n\right)  ,\label{eqPoof.11}\\%
\intertext{and}%
\pi\left(  f\right)  \left(  \xi,n\right)   &  :=\left(  f\left(  z^{N^{n}%
}\right)  \xi\left(  z\right)  ,n\right)  , \label{eqPoof.12}%
\end{align}
for $f,\xi\in L^{\infty}\left(  \mathbb{T}\right)  $, and $n=0,1,\dots$. A
direct substitution of the definitions then leads to commutativity of the
following three
commutative
diagrams (Figures 1a--1c):%
\[
\begin{minipage}{143pt}\mbox{\includegraphics
[bbllx=0bp,bblly=31bp,bburx=132bp,bbury=101bp,width=143pt]
{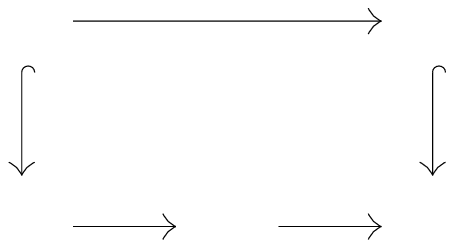}%
}\llap{\setlength{\unitlength}{16pt}\begin{picture}(8.8625,4.8)(-0.4,-0.4)
\put(0,0){\makebox(0,0){$\displaystyle\mathcal{V}_{1}$}}
\put(4,0){\makebox(0,0){$\displaystyle\mathcal{V}_{0}$}}
\put(8,0){\makebox(0,0){$\displaystyle\mathcal{V}_{1}$,}}
\put(0,4){\makebox(0,0){$\displaystyle\mathcal{V}_{0}$}}
\put(8,4){\makebox(0,0){$\displaystyle\mathcal{V}_{0}$}}
\put(2,-0.25){\makebox(0,0)[t]{$\scriptstyle U$}}
\put(6,-0.25){\makebox(0,0)[t]{$\scriptstyle J$}}
\put(0.25,2){\makebox(0,0)[l]{$\scriptstyle J$}}
\put(8.25,2){\makebox(0,0)[l]{$\scriptstyle J$}}
\put(4,3.75){\makebox(0,0)[t]{$\scriptstyle U$}}
\end{picture}}\\[6pt]\makebox[143pt]{\textsc{Figure }1\textup{a}%
}\end{minipage}%
\]%
\[
\begin{minipage}{143pt}\mbox{\includegraphics
[bbllx=0bp,bblly=31bp,bburx=132bp,bbury=101bp,width=143pt]
{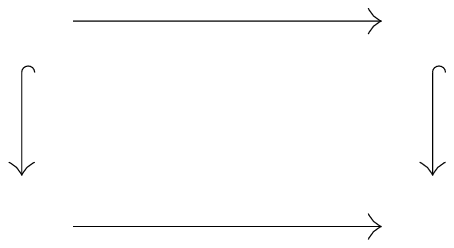}%
}\llap{\setlength{\unitlength}{16pt}\begin{picture}(8.8625,4.8)(-0.4,-0.4)
\put(0,0){\makebox(0,0){$\displaystyle\mathcal{V}_{n+1}$}}
\put(8,0){\makebox(0,0){$\displaystyle\mathcal{V}_{n}$,}}
\put(0,4){\makebox(0,0){$\displaystyle\mathcal{V}_{n}$}}
\put(8,4){\makebox(0,0){$\displaystyle\mathcal{V}_{n-1}$}}
\put(4,-0.25){\makebox(0,0)[t]{$\scriptstyle U$}}
\put(0.25,2){\makebox(0,0)[l]{$\scriptstyle J$}}
\put(8.25,2){\makebox(0,0)[l]{$\scriptstyle J$}}
\put(4,3.75){\makebox(0,0)[t]{$\scriptstyle U$}}
\end{picture}}\\[6pt]\makebox[143pt]{\textsc{Figure }1\textup{b}%
}\end{minipage}%
\]
for $n=1,2,\dots$, and%
\[
\begin{minipage}{143pt}\mbox{\includegraphics
[bbllx=0bp,bblly=31bp,bburx=132bp,bbury=101bp,width=143pt]
{funfig2a.eps}%
}\llap{\setlength{\unitlength}{16pt}\begin{picture}(8.8625,4.8)(-0.4,-0.4)
\put(0,0){\makebox(0,0){$\displaystyle\mathcal{V}_{n+1}$}}
\put(8,0){\makebox(0,0){$\displaystyle\mathcal{V}_{n+1}$.}}
\put(0,4){\makebox(0,0){$\displaystyle\mathcal{V}_{n}$}}
\put(8,4){\makebox(0,0){$\displaystyle\mathcal{V}_{n}$}}
\put(4,-0.25){\makebox(0,0)[t]{$\scriptstyle\pi\left( f\right) $}}
\put(0.25,2){\makebox(0,0)[l]{$\scriptstyle J$}}
\put(8.25,2){\makebox(0,0)[l]{$\scriptstyle J$}}
\put(4,3.75){\makebox(0,0)[t]{$\scriptstyle\pi\left( f\right) $}}
\end{picture}}\\[8pt]\makebox[143pt]{\textsc{Figure }1\textup{c}%
}\end{minipage}%
\]
In fact, substitution of (\ref{eqPoof.10})--(\ref{eqPoof.12}) into the
diagrams leads to the following, easily verified, identities (Figures 2a--2c):%
\[
\begin{minipage}{124pt}\mbox{\includegraphics
[bbllx=0bp,bblly=-42bp,bburx=124bp,bbury=121bp,width=124pt]
{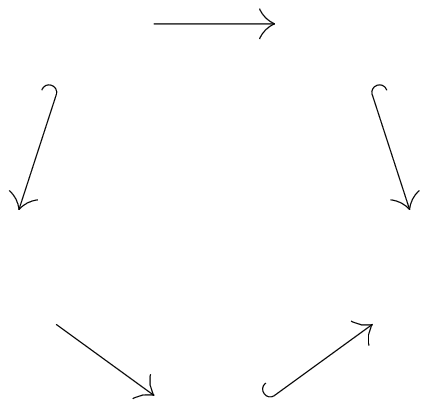}%
}\llap{\setlength{\unitlength}{56.5pt}\begin{picture}%
(2.195,2.888)(-1.097,-1.809)
\put(0.691,0.951){\makebox(0,0){\rlap
{$\displaystyle\left( S_{0}\xi,0\right) \in\mathcal{V}_{0}$}\phantom
{$\displaystyle\left( S_{0}\xi,0\right) $}}}
\put(-0.691,0.951){\makebox(0,0){\phantom{$\displaystyle\left( \xi,0\right
)$}\llap
{$\displaystyle\mathcal{V}_{0}\ni\left( \xi,0\right)$}}}
\put(-1.118,-0.363){\makebox(0,0){\phantom{$\displaystyle\left( \xi
\left( z^{N}\right) ,1\right)$}\llap
{$\displaystyle\mathcal{V}_{1}\ni\left( \xi\left( z^{N}\right) ,1\right)$}}}
\put(0,-1.123){\makebox(0,0)[t]{$\displaystyle\begin{matrix}
\displaystyle\left( m_{0}\left( z\right) \xi\left( z^{N}\right) ,0\right) \\
\makebox[0pt]{\hss$\displaystyle\cap$\hss}\makebox[0pt]
{\hss\rule[-0.15pt]{0.225pt}{6pt}\hss}  \\
\rlap{$\displaystyle\mathcal{V}_{0}$,}\phantom{\displaystyle\mathcal{V}}
\end{matrix}$}}
\put(1.118,-0.363){\makebox(0,0){\rlap
{$\displaystyle\left( \left( S_{0}\xi\right)  \left( z^{N}\right
) ,1\right) \in\mathcal{V}_{1}$}\phantom{$\displaystyle\left( \left( S_{0}%
\xi\right)  \left( z^{N}\right) ,1\right)$}}}
\put(0.968,0.315){\makebox(0,0)[l]{$\scriptstyle J$}}
\put(0,1.018){\makebox(0,0)[b]{$\scriptstyle U$}}
\put(-0.968,0.315){\makebox(0,0)[r]{$\scriptstyle J$}}
\put(-0.599,-0.824){\makebox(0,0)[tr]{$\scriptstyle U$}}
\put(0.599,-0.824){\makebox(0,0)[tl]{$\scriptstyle J$}}
\end{picture}}\linebreak\makebox[124pt]{\textsc{Figure }2\textup{a}%
}\end{minipage}%
\]%
\[
\begin{minipage}{191pt}\mbox{\includegraphics
[bbllx=0bp,bblly=31bp,bburx=132bp,bbury=101bp,width=143pt]
{funfig2a.eps}%
}\llap{\setlength{\unitlength}{16pt}\begin{picture}(8.8625,4.8)(-0.4,-0.4)
\put(0,0){\makebox(0,0){\llap{$\displaystyle\mathcal{V}_{n+1}\ni\left
( \xi\left( z^{N}\right) ,n+1\right) $\hskip-12pt}}}
\put(8,0){\makebox(0,0){\rlap{\hskip-12pt$\displaystyle\left( m_{0}%
\left( z^{N^{n}}%
\right) \xi\left( z^{N}\right) ,n\right) \in\mathcal{V}_{n}$,}}}
\put(0,4){\makebox(0,0){\llap{$\displaystyle\mathcal{V}_{n}\ni\left
( \xi,n\right) $\hskip-12pt}}}
\put(8,4){\makebox(0,0){\rlap{\hskip-12pt$\displaystyle\left( m_{0}%
\left( z^{N^{n-1}%
}\right) \xi\left( z\right) ,n-1\right) \in\mathcal{V}_{n-1}$}}}
\put(4,-0.25){\makebox(0,0)[t]{$\scriptstyle U$}}
\put(0.25,2){\makebox(0,0)[l]{$\scriptstyle J$}}
\put(8.25,2){\makebox(0,0)[l]{$\scriptstyle J$}}
\put(4,3.75){\makebox(0,0)[t]{$\scriptstyle U$}}
\end{picture}}\hspace*{48pt}\\[6pt]\makebox[191pt]{\textsc{Figure }2\textup
{b}%
}\end{minipage}%
\]
and the last simpler diagram, which does not involve a horizontal shift in the
$n$-index:%
\[
\begin{minipage}{191pt}\mbox{\includegraphics
[bbllx=0bp,bblly=31bp,bburx=132bp,bbury=101bp,width=143pt]
{funfig2a.eps}%
}\llap{\setlength{\unitlength}{16pt}\begin{picture}(8.8625,4.8)(-0.4,-0.4)
\put(0,0){\makebox(0,0){\llap{$\displaystyle\mathcal{V}_{n+1}\ni\left
( \xi\left( z^{N}\right) ,n+1\right) $\hskip-12pt}}}
\put(8,0){\makebox(0,0){\rlap{\hskip-12pt$\displaystyle\left( f\left
( z^{N^{n+1}}%
\right) \xi\left( z^{N}\right) ,n+1\right) \in\mathcal{V}_{n+1}$.}}}
\put(0,4){\makebox(0,0){\llap{$\displaystyle\mathcal{V}_{n}\ni\left
( \xi,n\right) $\hskip-12pt}}}
\put(8,4){\makebox(0,0){\rlap{\hskip-12pt$\displaystyle\left( f\left
( z^{N^{n}}%
\right) \xi\left( z\right) ,n\right) \in\mathcal{V}_{n}$}}}
\put(4,-0.25){\makebox(0,0)[t]{$\scriptstyle\pi\left( f\right) $}}
\put(0.25,2){\makebox(0,0)[l]{$\scriptstyle J$}}
\put(8.25,2){\makebox(0,0)[l]{$\scriptstyle J$}}
\put(4,3.75){\makebox(0,0)[t]{$\scriptstyle\pi\left( f\right) $}}
\end{picture}}\hspace*{48pt}\\[8pt]\makebox[191pt]{\textsc{Figure }2\textup
{c}%
}\end{minipage}%
\]
The purpose of the diagrams is to verify that the operator $U$ and the
representation $\pi$ (of $L^{\infty}\left(  \mathbb{T}\right)  $),
as defined from (\ref{eqPoof.10})--(\ref{eqPoof.12}),
pass to the
inductive limit construction which is obtained by the identification of
$\mathcal{H}_{n}$ with a closed subspace in $\mathcal{H}_{n+1}$ for each $n$,
and therefore, by iteration, in $\mathcal{H}_{n+k}$ for all $k=1,2,\dots$.
When the inductive limit
\begin{equation}
\mathcal{H}=\underset{n}{\varinjlim}\,\mathcal{H}_{n} \label{eqPoof.13}%
\end{equation}
is then formed, we get a well defined operator $U$ on $\mathcal{H}$, and a
representation $\pi$ of $L^{\infty}\left(  \mathbb{T}\right)  $ on
$\mathcal{H}$.
(Vectors $\xi $ in $\mathcal{H}$ may be characterized
by the following orthogonal expansion:
$\xi =\sum _{n=0}^{\infty }\xi _{n}$, where $\xi _{0}\in\mathcal{H}_{0}$, and
$\xi _{n}\in\mathcal{H}_{n}\ominus \mathcal{H}_{n-1}$, $n=1,2,\dots $; and
$\sum _{0}^{\infty }\left\|  \xi _{n}\right\|  ^{2}<\infty $. For an
alternative and purely function-theoretic characterization of the
Hilbert space $\mathcal{H}$, see also
Corollary \ref{CorPoof.pound} below.)
From Lemma \ref{LemPoof.2},
and
a direct verification, we also
get the identity
\[
U\pi\left(  f\right)  =\pi\left(  f\left(  z^{N}\right)  \right)  U
\]
for the corresponding operators on $\mathcal{H}$.

A final lemma now completes the proof of Theorem \ref{Thmax+b.3},
(\ref{Thmax+b.3(1)(1)})$\rightarrow$(\ref{Thmax+b.3(1)(2)}).

\begin{lemma}
\label{LemPoof.4}If $m_{0}\in L^{\infty}\left(  \mathbb{T}\right)  $ and
$Rh=h$, then $U$ is isometric, and if also $m_{0}$ vanishes on at most a
subset of $\mathbb{T}$ of measure zero, then $U$ is a unitary operator in
$\mathcal{H}$.
\end{lemma}

\begin{proof}
From the inductive construction
(\ref{eqPoof.10})--(\ref{eqPoof.11}),
we have $U\colon\mathcal{V}_{n+1}%
\rightarrow\mathcal{V}_{n}$, and we wish to pass $U$ to the completion
$U\colon\mathcal{H}_{n+1}\rightarrow\mathcal{H}_{n}$. That can be done if we
check first that $U$ is isometric from $\mathcal{V}_{n+1}$ to $\mathcal{V}%
_{n}$ for all $n=0,1,2,\dots$. We already checked, in fact, that $U$ is
isometric on $\mathcal{V}_{0}$, and therefore on the completion $L^{2}\left(
h\right)  $; that was Lemma \ref{LemPoof.2}.

Let $\xi\in L^{\infty}\left(  \mathbb{T}\right)  $. Then
\begin{align*}
\left\|  U\left(  \xi,n+1\right)  \right\|  _{\mathcal{H}}^{2}  &  =\left\|
\left(  m_{0}\left(  z^{N^{n}}\right)  \xi\left(  z\right)  ,n\right)
\right\|  _{\mathcal{H}}^{2}\\
&  =\int_{\mathbb{T}}R^{n}\left(  \left|  m_{0}\left(  z^{N^{n}}\right)
\xi\left(  z\right)  \right|  ^{2}h\left(  z\right)  \right)  \,d\mu\left(
z\right) \\
&  =\int_{\mathbb{T}}\left|  m_{0}\right|  ^{2}R^{n}\left(  \left|
\xi\right|  ^{2}h\right)  \,d\mu\\
&  =\frac{1}{N}\int_{\mathbb{T}}\sum_{w^{N}=z}\left|  m_{0}\left(  w\right)
\right|  ^{2}R^{n}\left(  \left|  \xi\right|  ^{2}h\right)  \left(  w\right)
\,d\mu\left(  z\right) \\
&  =\int_{\mathbb{T}}R^{n+1}\left(  \left|  \xi\right|  ^{2}h\right)  \left(
z\right)  \,d\mu\left(  z\right) \\
&  =\left\|  \left(  \xi,n+1\right)  \right\|  _{\mathcal{H}}^{2},
\end{align*}
which is the desired isometric property.

Using again the inductive limit construction of $\mathcal{H}$, we note that
$U$ will be unitary on $\mathcal{H}$, i.e., $U\left(  \mathcal{H}\right)
=\mathcal{H}$ if and only if%
\begin{equation}
U\left(  \mathcal{H}_{n+1}\right)  =\mathcal{H}_{n} \label{eqPoof.14}%
\end{equation}
for $n=0,1,2,\dots$. Equivalently, we must show that the spaces $\mathcal{H}%
_{n}\ominus U\left(  \mathcal{H}_{n+1}\right)  $ vanish for $n=0,1,\dots$. To
do this we need the following

\begin{claim}
\label{ClaPoof.5}The completion $\mathcal{H}_{n}=\widetilde{\mathcal{V}}_{n}$
in the norm on $\mathcal{V}_{n}$ consists of measurable functions $\xi$ on
$\mathbb{T}$ satisfying%
\begin{equation}
\int_{\mathbb{T}}R^{n}\left(  \left|  \xi\right|  ^{2}h\right)  \,d\mu<\infty.
\label{eqPoof.15}%
\end{equation}
\end{claim}

\begin{proof}
Let $\xi_{i}\in L^{\infty}\left(  \mathbb{T}\right)  $, $i=1,2,\dots$, and
suppose%
\[
\lim_{i,j\rightarrow\infty}\int_{\mathbb{T}}R^{n}\left(  \left|  \xi_{i}%
-\xi_{j}\right|  ^{2}h\right)  \,d\mu=0.
\]
Let%
\[
m_{0}^{\left(  n\right)  }\left(  z\right)  =m_{0}\left(  z\right)
m_{0}\left(  z^{N}\right)  \cdots m_{0}\left(  z^{N^{n-1}}\right)  .
\]
Then%
\[
\int_{\mathbb{T}}R^{n}\left(  \left|  \xi_{i}-\xi_{j}\right|  ^{2}h\right)
\,d\mu=\int_{\mathbb{T}}\left|  m_{0}^{\left(  n\right)  }\right|  ^{2}\left|
\xi_{i}-\xi_{j}\right|  ^{2}h\,d\mu.
\]
We conclude that there is a pointwise a.e.\ convergent subsequence $\xi
_{i_{1}},\xi_{i_{2}},\dots$ with limit $\xi$, say. We have%
\[
\int_{\mathbb{T}}\left|  m_{0}^{\left(  n\right)  }\right|  ^{2}\left|
\xi_{i_{k}}-\xi\right|  ^{2}h\,d\mu\underset{k\rightarrow\infty}%
{\longrightarrow}0,
\]
and%
\[
\int_{\mathbb{T}}\left|  m_{0}^{\left(  n\right)  }\right|  ^{2}\left|
\xi\right|  ^{2}h\,d\mu<\infty.
\]
Since%
\[
\int_{\mathbb{T}}\left|  m_{0}^{\left(  n\right)  }\right|  ^{2}\left|
\xi\right|  ^{2}h\,d\mu=\int_{\mathbb{T}}R^{n}\left(  \left|  \xi\right|
^{2}h\right)  \,d\mu,
\]
the claim follows.
\end{proof}

To prove the unitarity assertion of the lemma, we must show that if $\xi$
satisfies (\ref{eqPoof.15}) of Claim \ref{ClaPoof.5}, and if
\begin{equation}
\int_{\mathbb{T}}R^{n}\left(  \overline{\xi\left(  z\right)  }\,m_{0}\left(
z^{N^{n}}\right)  \eta\left(  z\right)  h\left(  z\right)  \right)
\,d\mu\left(  z\right)  =0 \label{eqPoof.16}%
\end{equation}
for all $\eta\in L^{\infty}\left(  \mathbb{T}\right)  $, then $\xi$ must
vanish a.e.\ on $\mathbb{T}$, and therefore
\[
\int_{\mathbb{T}}\left|  m_{0}^{\left(  n\right)  }\right|  ^{2}\left|
\xi\right|  ^{2}h\,d\mu=\int_{\mathbb{T}}R^{n}\left(  \left|  \xi\right|
^{2}h\right)  \,d\mu=0.
\]
Since $m_{0}\in L^{\infty}\left(  \mathbb{T}\right)  $, the function
$z\mapsto\overline{\xi\left(  z\right)  }\,m_{0}\left(  z^{N^{n}}\right)  $
also satisfies condition (\ref{eqPoof.15}) in the Claim. Since $\mathcal{H}%
_{n}=\widetilde{\mathcal{V}}_{n}$, we conclude that%
\[
\int_{\mathbb{T}}R^{n}\left(  \left|  \xi\left(  z\right)  m_{0}\left(
z^{N^{n}}\right)  \right|  ^{2}h\left(  z\right)  \right)  \,d\mu=0.
\]
But the integral is also%
\[
\int_{\mathbb{T}}\left|  m_{0}\right|  ^{2}R^{n}\left(  \left|  \xi\right|
^{2}h\right)  \,d\mu,
\]
and, if we now (finally!) use that $m_{0}$ does not vanish on a subset of
$\mathbb{T}$ of positive measure, we see that $R^{n}\left(  \left|
\xi\right|  ^{2}h\right)  $ must vanish pointwise a.e.\ on $\mathbb{T}$, and
therefore $\int_{\mathbb{T}}R^{n}\left(  \left|  \xi\right|  ^{2}h\right)
\,d\mu=0$, concluding the proof that%
\begin{equation}
\mathcal{H}_{n}\ominus U\left(  \mathcal{H}_{n+1}\right)  =0.\settowidth
{\qedskip}{$\displaystyle\mathcal{H}_{n}\ominus U\left( \mathcal{H}%
_{n+1}\right) =0.$}\addtolength{\qedskip}{-\textwidth}\rlap{\hbox
to-0.5\qedskip{\hfil\qed}} \label{eqPoof.17}%
\end{equation}
\renewcommand{\qed}{}
\end{proof}

To conclude the proof of Theorem \ref{Thmax+b.3}(\ref{Thmax+b.3(2)}), we only
need to identify $\varphi\in\mathcal{H}$ such that $U\varphi=\pi\left(
m_{0}\right)  \varphi$. Take $\varphi=\left(  \openone,0\right)  \sim\left(
\openone,1\right)  \sim\left(  \openone,2\right)  \sim\cdots$,
identification \emph{via} the isometry $J$ of (\ref{eqPoof.9}). Then%
\[
U\varphi=\left(  S_{0}\openone,0\right)  =\left(  m_{0},0\right)  =\pi\left(
m_{0}\right)  \left(  \openone,0\right)  =\pi\left(  m_{0}\right)  \varphi.
\]
It is clear from the construction that $\varphi$ is \emph{cyclic} for the
particular representation, i.e., $\tilde{\pi}\in\operatorname*{Rep}\left(
\mathfrak{A}_{N},\mathcal{H}\right)  $, which corresponds to the pair $\left(
U,\pi\right)  $ where $\pi$ is the $L^{\infty}\left(  \mathbb{T}\right)
$-representation which satisfies (\ref{eqax+b.4}).

The final assertion in Theorem \ref{Thmax+b.3}(\ref{Thmax+b.3(3)}) is that
$\tilde{\pi}$ is \emph{unique up to unitary equivalence.} The proof of this is
somewhat similar to the standard uniqueness part in the GNS construction: see,
e.g., \cite{BrRoI}. Since there are some differences as well, we sketch the details.

\begin{remark}
\label{RemPoof.6}Suppose $m_{0}$ does not vanish on a set of positive measure.
The function $\xi\left(  z\right)  :=\frac{1}{m_{0}\left(  z\right)  }$
represents an element in $\mathcal{H}_{1}$ \emph{via} $\left(  \xi,1\right)
$, even though generally $\frac{1}{m_{0}\left(  z\right)  }$ is not in
$L^{\infty}\left(  \mathbb{T}\right)  $; see, e.g., Examples
\textup{\ref{ExaWave.3(1)}--\ref{ExaWave.3(2)}.} We have%
\[
\left\|  \left(  \xi,1\right)  \right\|  _{\mathcal{H}_{1}}^{2}=\int
_{\mathbb{T}}R\left(  \frac{1}{\left|  m_{0}\left(  z\right)  \right|
^{2\mathstrut}}h\left(  z\right)  \right)  \,d\mu\left(  z\right)
=\int_{\mathbb{T}}h\left(  z\right)  \,d\mu\left(  z\right)  =\left\|
\varphi\right\|  _{\mathcal{H}}^{2},
\]
where $\varphi$ is the cyclic vector which corresponds to a given $h\in
L^{1}\left(  \mathbb{T}\right)  $, $h\geq0$, $Rh=h$. If $\left(
\mathcal{H}_{h},\pi,U\right)  $ denotes the representation of $\mathfrak{A}_{N}$
which is induced from $h$ \emph{via} Theorem \textup{\ref{Thmax+b.3},} then a
simple calculation shows that%
\begin{equation}
U^{\ast}\left(  \varphi\right)  =U^{-1}\left(  \varphi\right)  =\left(
\frac{1}{m_{0}},1\right)  \in\mathcal{H}_{1}. \label{eqPoof.18}%
\end{equation}
\end{remark}

\begin{proof}
[Proof of uniqueness in Theorem \textup{\ref{Thmax+b.3}}]The uniqueness up to
unitary equivalence is only asserted when $U$ is unitary, i.e., when the
representation $\tilde{\pi}\in\operatorname*{Rep}\left(  \mathfrak{A}%
_{N},\mathcal{H}\right)  $ is constructed from a given $m_{0}\in L^{\infty
}\left(  \mathbb{T}\right)  $, and an $h\in L^{1}\left(  \mathbb{T}\right)  $,
$h\geq0$, and $R_{m_{0}}\left(  h\right)  =h$. The determining conditions for
$\tilde{\pi}$ are:
\begin{enumerate}
\item \label{Proofofuniqueness(1)}\emph{cyclicity,}

\item \label{Proofofuniqueness(2)}$\ip{\varphi}{\pi\left( f\right) \varphi
}=\int_{\mathbb{T}}fh\,d\mu$, and

\item \label{Proofofuniqueness(3)}$U\varphi=\pi\left(  m_{0}\right)  \varphi$.
\end{enumerate}
But we just established that%
\[
U^{\ast}\left(  \mathcal{H}_{n}\right)  =\mathcal{H}_{n+1},
\]
$n=0,1,2,\dots$. Hence%
\begin{equation}
U^{\ast\,n}\left(  \mathcal{H}_{0}\right)  =\mathcal{H}_{n}, \label{eqPoof.19}%
\end{equation}
where $\mathcal{H}_{0}\subset\mathcal{H}_{1}\subset\mathcal{H}_{2}%
\subset\cdots$ is the resolution which defines $\mathcal{H}$ as an inductive
limit. In Lemma \ref{LemPoof.2}, we saw that $\mathcal{H}_{0}\simeq
L^{2}\left(  h\right)  $ with the isomorphism defined by $W\colon
\mathcal{H}_{0}\underset{\simeq}{\longrightarrow}L^{2}\left(  h\right)  $,
$W\pi\left(  f\right)  \varphi=f$, $f\in L^{\infty}\left(  \mathbb{T}\right)
$. This was based on the computation%
\begin{equation}
\left\|  \pi\left(  f\right)  \varphi\right\|  _{\mathcal{H}}^{2}=\int\left|
f\right|  ^{2}R_{m_{0}}\left(  h\right)  \,d\mu. \label{eqPoof.20}%
\end{equation}
If $R_{m_{0}}h=h$, then we get%
\[
\left\|  \pi\left(  f\right)  \varphi\right\|  _{\mathcal{H}}=\left\|
f\right\|  _{L^{2}\left(  h\right)  }\text{.}%
\]
But we also saw that $W$ intertwines $U$ and $S_{0}$, i.e., that
\begin{equation}
WU=S_{0}W. \label{eqPoof.21}%
\end{equation}
Using the identity (\ref{eqPoof.19}), we conclude that the intertwining
property on the $\mathcal{H}_{0}$ level extends to the $\mathcal{H}_{n}%
$-spaces for all $n\geq0$, and therefore that any two $\tilde{\pi}%
\in\operatorname*{Rep}\left(  \mathfrak{A}_{N},\mathcal{H}\right)  $ and
$\tilde{\pi}^{\prime}\in\operatorname*{Rep}\left(  \mathfrak{A}_{N},\mathcal{H}%
^{\prime}\right)  $ which both satisfy (\ref{Proofofuniqueness(1)}%
)--(\ref{Proofofuniqueness(3)}) must be unitarily equivalent. This completes
the last part of the proof of Theorem \ref{Thmax+b.3}.
\end{proof}

We conclude with a lemma which illustrates the above construction, and which
will be needed in the sequel.

\begin{lemma}
\label{LemPoof.7}The orthogonal projection of $\mathcal{H}_{1}$ onto
$\mathcal{H}_{0}$ is given by%
\begin{equation}
J^{\ast}\left(  \xi,1\right)  =\left(  \frac{R\left(  \xi h\right)  }%
{h},0\right)  . \label{eqPoof.22}%
\end{equation}
The condition on $\xi$ characterizing membership in $\mathcal{H}_{1}$ is%
\begin{equation}
\int_{\mathbb{T}}R\left(  \left|  \xi\right|  ^{2}h\right)  \,d\mu
=\int_{\mathbb{T}}\left|  m_{0}\right|  ^{2}\left|  \xi\right|  ^{2}%
h\,d\mu<\infty, \label{eqPoof.23}%
\end{equation}
and implied in \textup{(\ref{eqPoof.22})} is the assertion that the expression
$\frac{R\left(  \xi h\right)  \left(  z\right)  }{h\left(  z\right)  }$ is
then well defined, and that it is in $L^{2}\left(  h\right)  $, i.e., that%
\begin{equation}
\int_{\mathbb{T}}\frac{\left|  R\left(  \xi h\right)  \left(  z\right)
\right|  ^{2}}{h\left(  z\right)  }\,d\mu\left(  z\right)  <\infty.
\label{eqPoof.24}%
\end{equation}
\end{lemma}

\begin{proof}
To see that the expressions under the integral make sense pointwise, the
following estimate is needed (obtained by an iteration of Schwarz's
estimate!):%
\[
\left|  R\left(  \xi h\right)  \left(  z\right)  \right|  \leq R\left(
\left|  \xi\right|  ^{2}h\right)  \left(  z\right)  ^{\frac{1}{2}}h\left(
z\right)  ^{\frac{1}{2}}\leq\dots\leq\left(  R\left(  \left|  \xi\right|
^{2^{n}}h\right)  \left(  z\right)  \right)  ^{\frac{1}{2^{n}}%
}h\left(  z\right)  ^{\frac{1}{2}+\frac{1}{4}+\dots+\frac{1}{2^{n}}}.
\]
To prove (\ref{eqPoof.22}), we use the definition $J\colon\left(
\eta,0\right)  \mapsto\left(  \eta\left(  z^{N}\right)  ,1\right)  $, and
further compute that:%
\begin{align*}
\ip{J^{\ast}\left( \xi,1\right) }{\left( \eta,0\right) }  &  =\int
_{\mathbb{T}}R\left(  \bar{\xi}\left(  z\right)  \eta\left(  z^{N}\right)
h\left(  z\right)  \right)  \,d\mu\\
&  =\int_{\mathbb{T}}\eta\left(  z\right)  R\left(  \bar{\xi}h\right)  \left(
z\right)  \,d\mu\left(  z\right) \\
&  =\ip{\frac{R\left( \xi h\right) }{h}}{\eta}_{L^{2}\left(  h\right)  },
\end{align*}
and (\ref{eqPoof.22}) follows from this.
Since we saw that $J$ is isometric when $R\left(  h\right)  =h$ holds
(\emph{cf.\ }Lemma \ref{LemPoofIsomJ}), it
follows that $JJ^{\ast}$ is the desired projection of $\mathcal{H}_{1}$ onto
$\mathcal{H}_{0}$; but we may use $J$ in making an isometric identification.
\end{proof}

As a corollary, we get

\begin{corollary}
\label{CorPoof.pound}
The elements in the Hilbert space
$\mathcal{H}$, which is constructed from
the resolution
$\mathcal{H}_{0}\subset \mathcal{H}_{1}\subset \mathcal{H}_{2}\subset \cdots $,
may be given by precisely the sequences
of measurable functions $\left(  \xi _{n}\right)  $,
$n=0,1,\dots $, such that
\[
\sup _{n}\int _{\mathbb{T}}R^{n}\left(  \left|  \xi _{n}\right|  ^{2}h\right)  \,d\mu <\infty ,
\]
and
\[
R\left(  \xi _{n+1}h\right)  =\xi _{n}h,\qquad n=0,1,\dots .
\]
\end{corollary}

A system like this is also called a
martingale.

\section{\label{Wave}Wavelet filters}

In Chapter \ref{ax+b} we showed that if $m_{0}$ is a wavelet filter (see
Definition \ref{Defax+b.1}) then \emph{one of} the representations from
Theorem \ref{Thmax+b.3} may be realized in $\mathcal{H}=L^{2}\left(
\mathbb{R}\right)  $. Specifically, there is a cyclic vector $\varphi\in
L^{2}\left(  \mathbb{R}\right)  $ and a corresponding $h_{\varphi}\in
L^{1}\left(  \mathbb{T}\right)  $ such that $h_{\varphi}\geq0$ and $R\left(
h_{\varphi}\right)  =h_{\varphi}$. The assertion here is that, in this case,
the solution $\varphi$ to $U\varphi=\pi\left(  m_{0}\right)  \varphi$ may be
found in $L^{2}\left(  \mathbb{R}\right)  .$ For this, we need the
normalization (\ref{eqax+b.7}), i.e.,
\begin{equation}
\sum_{k=0}^{N-1}\left|  m_{0}\left(  e^{i\frac{2\pi k}{N}}z\right)  \right|
^{2}\equiv N\text{\qquad(a.e.\ on }\mathbb{T}\text{),} \label{eqWave.1}%
\end{equation}
which is one of the defining conditions on a wavelet filter. For $N=2$, it
reads:%
\begin{equation}
\left|  m_{0}\left(  z\right)  \right|  ^{2}+\left|  m_{0}\left(  -z\right)
\right|  ^{2}=2,\qquad z\in\mathbb{T}. \label{eqWave.2}%
\end{equation}
In that case, set%
\begin{equation}
m_{1}\left(  z\right)   :=z\,\overline{m_{0}\left(  -z\right)  },\qquad
z\in\mathbb{T},\label{eqWave.3}
\end{equation}
and
\begin{equation}
\left(  S_{j}f\right)  \left(  z\right)   =m_{j}\left(  z\right)  f\left(
z^{2}\right)  ,\qquad f\in L^{2}\left(  \mathbb{T}\right)  . \label{eqWave.4}%
\end{equation}

\begin{lemma}
\label{LemWave.1}The operators $S_{0}$ and $S_{1}$ are isometries in
$L^{2}\left(  \mathbb{T}\right)  $ and satisfy%
\begin{equation}
S_{i}^{\ast}S_{j}^{{}}=\delta_{ij}^{{}}I\text{\qquad and\qquad}\sum
_{i=0}^{1}S_{i}^{{}}S_{i}^{\ast}=I, \label{eqWave.4bis}%
\end{equation}
where $I$ denotes the identity operator.
\end{lemma}

\begin{proof}
See \cite{BrJo97b}.
The conclusion of the lemma
may be rephrased as the assertion
that the operators in (\ref{eqWave.4}) define
a representation of the Cuntz
$C^{\ast }$-algebra $\mathcal{O}_{2}$. There is a
similar conclusion for $\mathcal{O}_{N}$.
We use standard notation from
\cite{Cun77}
and \cite{BrRoI}
for the Cuntz algebras and
their representations.
\end{proof}

\begin{lemma}
\label{LemWave.2}Let $h\in L^{1}\left(  \mathbb{T}\right)  $, $h\geq0$,
satisfy $R\left(  h\right)  =h$. Then $S_{0}$ is isometric in $L^{2}\left(
h\right)  :=L^{2}\left(  \mathbb{T},h\,d\mu\right)  $, but $S_{1}$ is
generally not isometric in $L^{2}\left(  h\right)  $.
\end{lemma}

\begin{proof}
We already saw in Chapter \ref{ax+b} that, if $h\in L^{1}\left(
\mathbb{T}\right)  $, $h\geq0$, is given, then $S_{0}$ is isometric in
$L^{2}\left(  h\right)  $ if and only if $R\left(  h\right)  =h$. The same
argument shows that the condition for $S_{1}$ to be isometric in $L^{2}\left(
h\right)  $ is $R\left(  \check{h}\right)  =h$, where $\check{h}\left(
z\right)  :=h\left(  -z\right)  $, $z\in\mathbb{T}$, and there are easy
examples where this is not satisfied.
\renewcommand{\qed}{}
\end{proof}

\begin{example}
\label{ExaWave.3(1)}Let%
\begin{equation}
m_{0}\left(  z\right)  :=\frac{1}{\sqrt{2}\mathstrut}\left(  1+z^{3}\right)  ,
\label{eqWave.5}%
\end{equation}
and let $R$ be the corresponding Ruelle operator. It is easy to see that the
scaling function $\varphi$ for the wavelet representation is%
\begin{equation}
\varphi\left(  x\right)  =\frac{1}{3}\chi_{\left[  0,3\right\rangle }^{{}%
}\left(  x\right)  ,\qquad x\in\mathbb{R}. \label{eqWave.6}%
\end{equation}
The corresponding harmonic function $h_{\varphi}$ (i.e., $R\left(  h_{\varphi
}\right)  =h_{\varphi}$) is computed as follows:%
\begin{align}
h_{\varphi}\left(  z\right)   &  =\sum_{n}z^{n}\ip{\pi\left( e_{n}%
\right) \varphi}{\varphi}\label{eqWave.7}\\
&  =\frac{1}{9}\left(  z^{-2}+2z^{-1}+3+2z+z^{2}\right) \nonumber\\
&  =\frac{1}{9}\left(  1+2\cos\omega\right)  ^{2},\nonumber
\end{align}
where $z:=e^{-i\omega}$, $\omega\in\mathbb{R}$. Then $\check{h}_{\varphi
}\left(  z\right)  =h_{\varphi}\left(  -z\right)  \approx h_{\varphi}\left(
\omega+\pi\right)  =\frac{1}{9}\left(  1-2\cos\omega\right)  ^{2}$, and a
calculation shows that
\[
R\left(  \check{h}_{\varphi}\right)  \left(  e^{-i\omega}\right)  =h_{\varphi
}\left(  e^{-i\omega}\right)  -\frac{4}{9}\cos\left(  \frac{\omega}{2}\right)
\left(  1+\cos\frac{3\omega}{2}\right)  .
\]
In this case, therefore, $h=h_{\varphi}$ does not satisfy $R\left(  \check
{h}\right)  =h$, and so $S_{1}$ is not isometric in $L^{2}\left(  \frac{1}%
{9}\left(  1+2\cos\omega\right)  ^{2}\right)  $.
\hfill\qed
\end{example}

\begin{remark}
\label{RemWave.NoTrace}The structure of the states on $\mathfrak{A}_{2}$ which are
induced from solutions
\[
R\left(  h\right)  =h,\qquad h\in L^{1}\left(  \mathbb{T}\right)  ,\;h\geq0,
\]
is not yet well understood. But these states are clearly \emph{not tracial,}
not even for the simplest wavelet representations such as the one described
above in Example \textup{\ref{ExaWave.3(1)}. Recall, if }$\varphi$ is the
scaling function \textup{(\ref{eqWave.6})} in $L^{2}\left(  \mathbb{R}\right)
$, then the corresponding state $\sigma_{\varphi}$ is%
\[
\sigma_{\varphi}\left(  A\right)   :=\ip{\varphi}{\tilde{\pi}_{\varphi
}\left( A\right) \varphi}_{L^{2}\left(  \mathbb{R}\right)  },\qquad
A\in\mathfrak{A}_{2},
\]
and
\[
h_{\varphi}\left(  e^{-i\omega}\right)   =\frac{1}{9}\left(  1+2\cos
\omega\right)  ^{2}.
\]
The state $\sigma_{\varphi}$ is \emph{not tracial} because%
\[
\sigma_{\varphi}\left(  UVU^{-1}\right)  \neq\sigma_{\varphi}\left(
V\right)
\]
where%
\[
\left(  V\psi\right)  \left(  x\right)   =\psi\left(  x-1\right)  ,
\]
and
\[
\left(  U\psi\right)  \left(  x\right)   =\frac{1}{\sqrt{2}\mathstrut}%
\psi\left(  \frac{x}{2}\right)  ,\qquad\psi\in L^{2}\left(  \mathbb{R}\right)
.
\]
Since $UVU^{-1}=V^{2}$, we need only check that $\sigma_{\varphi}\left(
V^{2}\right)  \neq\sigma_{\varphi}\left(  V\right)  $; and a direct
calculation yields%
\begin{align*}
\sigma_{\varphi}\left(  V^{2}\right)   &  =\int_{\mathbb{T}}e_{2}h_{\varphi
}\,d\mu=\frac{1}{9},\\%
\intertext{while}%
\sigma_{\varphi}\left(  V\right)   &  =\int_{\mathbb{T}}e_{1}h_{\varphi}%
\,d\mu=\frac{2}{9}.
\end{align*}
\end{remark}

\begin{example}
\label{ExaWave.3(2)}We now turn to the representation associated with the
solution $R\openone=\openone$. \textup{(}Recall for the wavelet filters,
property \textup{(\ref{eqIntr.6})} or \textup{(\ref{eqax+b.7})} ensures that
the constant function $\openone$ is also an eigenfunction.\textup{)} For the
representation $\left(  \mathcal{H},\pi,U\right)  $ induced from $h=\openone$,
we may take%
\begin{equation}
\varphi=\frac{1}{\sqrt{3}\mathstrut}\left(  \chi_{\left[  0,1\right\rangle
}^{{}}\oplus\chi_{\left[  1,2\right\rangle }^{{}}\oplus\chi_{\left[
2,3\right\rangle }^{{}}\right)  \label{eqWave.8}%
\end{equation}
in $\mathcal{H}=L^{2}\left(  \mathbb{R}\right)  \oplus L^{2}\left(
\mathbb{R}\right)  \oplus L^{2}\left(  \mathbb{R}\right)  $. On the direct sum
$\mathcal{H}:=\sideset{}{^{\smash{\oplus}}}{\sum}L^{2}\left(  \mathbb{R}%
\right)  $ we introduce the usual representation $\tilde{\pi}=\left(
\pi,U\right)  $:%
\[%
\begin{cases}
\displaystyle\pi\left( f\right) \left( \psi_{i}\right) :=\left( \pi
\left( f\right) \psi
_{i}\right) , & f\in L^{\infty}\left( \mathbb{T}\right) ,  \\
\displaystyle U\left( \psi_{i}\right) :=\left( U\psi_{i}\right) , &
\end{cases}%
\]
where $\pi$ and $U$ on $L^{2}\left(  \mathbb{R}\right)  $ are given by the
usual formulas \textup{(\ref{eqax+b.5})} and \textup{(\ref{eqax+b.6}).} It is
clear that there is an isometry%
\[
\sideset{}{^{\smash{\oplus}}}{\sum}L^{2}\left(  \mathbb{R}\right)  \overset
{W}{\longrightarrow}L^{2}\left(  \mathbb{R}\right)
\]
which intertwines the respective representations of $L^{\infty}\left(
\mathbb{T}\right)  $. The vector $\varphi$ in $\mathcal{H}=L^{2}\left(
\mathbb{R}\right)  \oplus L^{2}\left(  \mathbb{R}\right)  \oplus L^{2}\left(
\mathbb{R}\right)  $ then satisfies%
\begin{align}
UW\varphi &  =W\pi\left(  m_{0}\right)  \varphi\label{eqWave.9}\\%
\intertext{or equivalently}%
W^{\ast}UW\varphi &  =\pi\left(  m_{0}\right)  \varphi. \label{eqWave.10}%
\end{align}
It is immediate from \textup{(\ref{eqWave.8})} that $h_{\varphi}=\openone$,
and we have identified the representation from Theorem \textup{\ref{Thmax+b.3}%
} induced by $h_{\varphi}=\openone$. \textup{(}We show later that $\tilde{\pi
}$ is the sum of \emph{two} mutually inequivalent irreducible representations.\textup{)}
\end{example}

Before getting to the structure theorem for the wavelet filters, we need a lemma:

\begin{lemma}
\label{LemWave.4}Consider the wavelet representation $\left(  L^{2}\left(
\mathbb{R}\right)  ,\pi,U\right)  $ defined from a given wavelet filter
$m_{0}$. Then for every $\varphi,\psi\in L^{2}\left(  \mathbb{R}\right)  $,
there is an $L^{1}\left(  \mathbb{T}\right)  $-function $h\left(  z\right)
=H\left(  \varphi,\psi\right)  \left(  z\right)  $ such that
\begin{equation}
\ip{\varphi}{\pi\left( f\right) \psi}=\int_{\mathbb{T}}f\left(  z\right)
H\left(  \varphi,\psi\right)  \left(  z\right)  \,d\mu\left(  z\right)  ,
\label{eqWave.11}%
\end{equation}
and it is represented by the Fourier expansion%
\begin{equation}
H\left(  \varphi,\psi\right)  \left(  z\right)  =\sum_{n\in\mathbb{Z}}z^{n}%
\ip{\pi\left( e_{n}\right) \varphi}{\psi}. \label{eqWave.12}%
\end{equation}
\end{lemma}

\begin{proof}
Let%
\[
Z\colon L^{2}\left(  \mathbb{R}\right)  \longrightarrow L^{2}\left(
\mathbb{T}\times\left[  0,1\right\rangle \right)
\]
denote the Zak transform; see \cite[p.\ 109]{Dau92} for details. We have%
\[
\left(  Z\psi\right)  \left(  z,x\right)  =\sum_{n\in\mathbb{Z}}z^{n}%
\psi\left(  x+n\right)  ,\qquad\psi\in L^{2}\left(  \mathbb{R}\right)
,\;x\in\mathbb{R},
\]
and%
\begin{equation}
\int_{\mathbb{T}}\int_{0}^{1}\left|  Z\psi\left(  z,x\right)  \right|
^{2}\,dx\,d\mu\left(  z\right)  =\int_{\mathbb{R}}\left|  \psi\left(
x\right)  \right|  ^{2}\,dx, \label{eqWave.13}%
\end{equation}
and $Z$ maps $L^{2}\left(  \mathbb{R}\right)  $ \emph{onto} $L^{2}\left(
\mathbb{T}\times\left[  0,1\right\rangle \right)  $. It follows, from
(\ref{eqWave.13}), polarization, and Fubini's theorem, that the function%
\begin{equation}
H\left(  \varphi,\psi\right)  \left(  z\right)  =\int_{0}^{1}\overline
{Z\varphi\left(  z,x\right)  }\,Z\psi\left(  z,x\right)  \,dx  \label{eqWave.13bis}
\end{equation}
is in $L^{1}\left(  \mathbb{T}\right)  $ for all $\varphi,\psi\in L^{2}\left(
\mathbb{R}\right)  $, and that this function satisfies the two desired
properties (\ref{eqWave.11}) and (\ref{eqWave.12}) stated in the lemma.
\end{proof}

\begin{remark}
The connection between the two
operators $M$ and $R$ of the
Introduction may be expressed with \textup{(\ref{eqWave.13bis})} as follows:
\[
R\left(  H\left(  \varphi ,\psi \right)  \right)  =H\left(  M\varphi ,M\psi \right)  .
\]
This identity, which is equivalent to \textup{(\ref{eqIntr.27})} above,
also provides the
direct link between spectral
theory of $R$ and approximation
properties of $\left\{  M^{n}\psi \mid n=1,2,\dots \right\}  $ as
$n\rightarrow \infty $. For details, see
\cite{BrJo98b}, \cite{Jor98a}, and \cite{BJR97}.
\end{remark}

Let $m_{0}$ be a wavelet filter with corresponding scaling function
$\varphi\in L^{2}\left(  \mathbb{R}\right)  $ and Ruelle operator $R$ in
$L^{1}\left(  \mathbb{T}\right)  $. Then%
\[
h_{\varphi}\left(  z\right)  =\sum_{n\in\mathbb{Z}}z^{n}\ip{\pi\left
( e_{n}\right) \varphi}{\varphi}_{L^{2}\left(  \mathbb{R}\right)  }%
\]
is the solution (in $L^{1}\left(  \mathbb{T}\right)  $) to $R\left(
h_{\varphi}\right)  =h_{\varphi}$ which generates the wavelet representation
$\tilde{\pi}_{\varphi}$ of $\mathfrak{A}_{N}$ on $L^{2}\left(  \mathbb{R}\right)
$.

Similarly, if $h\in L^{1}\left(  \mathbb{T}\right)  $, $h\geq0$, is given such
that $R\left(  h\right)  =h$, then, by Theorem \ref{Thmax+b.3}, there is a
(unique up to unitary equivalence) cyclic representation $\tilde{\pi}_{h}%
\in\operatorname*{Rep}\left(  \mathfrak{A}_{N},\mathcal{H}_{h}\right)  $ such that%
\begin{equation}
\ip{\Phi}{\pi_{h}\left( f\right) \Phi}=\int_{\mathbb{T}}fh\,d\mu
\label{eqWave.14}%
\end{equation}
and $U\Phi=\pi\left(  m_{0}\right)  \Phi$ where $\Phi\in\mathcal{H}_{h}$ is
the cyclic vector.

\begin{theorem}
\label{ThmWave.5}If $h\in L^{1}\left(  \mathbb{T}\right)  $ is given such that
$R\left(  h\right)  =h$, and if, for some $c\in\mathbb{R}_{+}$, we have
$h_{\varphi}\leq ch$, then the wavelet representation $\tilde{\pi}_{\varphi}$
is contained in $\tilde{\pi}_{h}$, i.e., we have $\mathcal{H}_{h}=L^{2}\left(
\mathbb{R}\right)  \oplus\mathcal{K}$ and $\tilde{\pi}_{h}=\tilde{\pi
}_{\varphi}\oplus\tilde{\pi}_{\mathcal{K}}$ for some $\tilde{\pi}%
_{\mathcal{K}}\in\operatorname*{Rep}\left(  \mathfrak{A}_{N},\mathcal{K}\right)  $.
\end{theorem}

\begin{remark}
\label{RemWave.TwoIrrReps}When the theorem is applied to Example
\textup{\ref{ExaWave.3(2)},} i.e., the representation $\tilde{\pi}_{1}$ of
$\mathfrak{A}_{2}$ which is induced by the pair $m_{0}=\frac{1}{\sqrt{2}%
\mathstrut}\left(  1+z^{3}\right)  $, $h_{1}=\openone$, we see that
$\tilde{\pi}_{1}$ is the sum of two mutually inequivalent irreducible
representations, $\tilde{\pi}_{\varphi}$ \textup{(}the wavelet
representation\textup{)} being one of them. The other one in $\mathcal{K}%
\subset\mathcal{H}_{1}$ is in fact irreducible.
\end{remark}

Let $\varphi\in L^{2}\left(  \mathbb{R}\right)  $ be given by
\textup{(\ref{eqWave.6}),} i.e., $\varphi=\frac{1}{3}\chi_{\left[
0,3\right\rangle }^{{}}$, and let
\[
h\left(  z\right)  =\frac{1}{9}\left(  3+2\left(  z+z^{-1}\right)
+z^{2}+z^{-2}\right)
\]
be the function \textup{(\ref{eqWave.7})} which defines the wavelet
representation $\tilde{\pi}_{\varphi}$ in $L^{2}\left(  \mathbb{R}\right)  $
for $m_{0}=\frac{1}{\sqrt{2}\mathstrut}\left(  1+z^{3}\right)  $. Recall
$\tilde{\pi}_{\varphi}=\left(  U,\pi_{0}\right)  $ is given by $U\psi\left(
x\right)  =\frac{1}{\sqrt{2}\mathstrut}\psi\left(  \frac{x}{2}\right)  $, and
$\pi_{0}\left(  e_{n}\right)  \psi=\psi\left(  x-n\right)  $, $\psi\in
L^{2}\left(  \mathbb{R}\right)  $. Let $\mathcal{K}=L^{2}\left(
\mathbb{R}\right)  \oplus L^{2}\left(  \mathbb{R}\right)  $, and
$\rho=e^{i\frac{2\pi}{3}}=-\frac{1}{2}+i\frac{\sqrt{3}}{2}$. Set
$U_{\mathcal{K}}=U\oplus U$ on $\mathcal{K}=L^{2}\left(  \mathbb{R}\right)
\oplus L^{2}\left(  \mathbb{R}\right)  $, $\alpha_{\rho}\left(  f\right)
\left(  z\right)  =f\left(  \rho z\right)  $, $f\in L^{\infty}\left(
\mathbb{T}\right)  $, and $\pi_{\mathcal{K}}\left(  f\right)  =\pi_{0}\left(
\alpha_{\rho}\left(  f\right)  \right)  \oplus\pi_{0}\left(  \alpha_{\rho^{2}%
}\left(  f\right)  \right)  $. Let $\tilde{\pi}_{\mathcal{K}}\in
\operatorname*{Rep}\left(  \mathfrak{A}_{2},\mathcal{K}\right)  $ be the
corresponding representation. We make the following claims:

\begin{claims}
If $\tilde{\pi}_{1}$ is the representation induced from $h_{1}\equiv\openone$,
i.e., the representation of Example \textup{\ref{ExaWave.3(2)},} then

\begin{enumerate}
\item \label{RemWave.TwoIrrRepsClaim(1)} $\tilde{\pi}_{1}=\tilde{\pi}%
_{\varphi}\oplus\tilde{\pi}_{\mathcal{K}},$ and

\item \label{RemWave.TwoIrrRepsClaim(2)}$\tilde{\pi}_{\mathcal{K}}$ is
irreducible on $\mathcal{K}=L^{2}\left(  \mathbb{R}\right)  \oplus
L^{2}\left(  \mathbb{R}\right)  $.
\end{enumerate}
\end{claims}

\begin{proof}
The formulas for $m_{0}$ and $h$ show that $m_{0}\left(  \rho z\right)
=m_{0}\left(  z\right)  $, $z\in\mathbb{T}$, and
\begin{equation}
h\left(  z\right)  +h\left(  \rho z\right)  +h\left(  \rho^{2}z\right)
\equiv\openone. \label{eqWave.sum}%
\end{equation}
Moreover, we have\renewcommand{\theenumi}{\alph{enumi}}
\begin{enumerate}
\item \label{RemWave.TwoIrrRepsPoof(1)}$\displaystyle R\left(  h\right)  =h$,

\item \label{RemWave.TwoIrrRepsPoof(2)}$\displaystyle R\left(  h\left(
\rho\,\cdot\,\right)  \right)  \left(  z\right)  =h\left(  \rho^{2}z\right)
,$ and

\item \label{RemWave.TwoIrrRepsPoof(3)}$\displaystyle R\left(  h\left(
\rho^{2}\,\cdot\,\right)  \right)  \left(  z\right)  =h\left(  \rho z\right)  .$
\end{enumerate}
The first assertions, (\ref{eqWave.sum}), and
(\ref{RemWave.TwoIrrRepsPoof(1)}) are immediate from substitution, so we check
only the last:%
\begin{align}
R\left(  h\left(  \rho\,\cdot\,\right)  \right)  \left(  z\right)   &
=\frac{1}{2}\sum_{w^{2}=z}\left|  m_{0}\left(  w\right)  \right|  ^{2}h\left(
\rho w\right) \tag{\emph{Ad} \textup{(\ref{RemWave.TwoIrrRepsPoof(2)})}}\\
&  =\frac{1}{2}\sum_{w^{2}=\rho^{2}z}\left|  m_{0}\left(  w\right)  \right|
^{2}h\left(  w\right) \nonumber\\
&  =R\left(  h\right)  \left(  \rho^{2}z\right) \nonumber\\
&  =h\left(  \rho^{2}z\right)  ,\nonumber
\end{align}
where we used $R\left(  h\right)  =h$ in the computation. The proof of
(\ref{RemWave.TwoIrrRepsPoof(3)}) is the same.

It follows that%
\begin{align*}
h_{\mathcal{K}}\left(  z\right)   &  :=1-h\left(  z\right)  =h\left(  \rho
z\right)  +h\left(  \rho^{2}z\right) \\%
\intertext{or}%
h_{\mathcal{K}}  &  =\alpha_{\rho}\left(  h\right)  +\alpha_{\rho^{2}}\left(
h\right)
\end{align*}
satisfies $R\left(  h_{\mathcal{K}}\right)  =h_{\mathcal{K}}$ and
$h_{\mathcal{K}}\geq0$. Let $\tilde{\pi}_{\mathcal{K}}$ be the corresponding
representation which is induced from $h_{\mathcal{K}}$ \emph{via} Theorem
\textup{\ref{Thmax+b.3}.} It then follows from \textup{(\ref{eqWave.sum})}
that Claim \textup{(\ref{RemWave.TwoIrrRepsClaim(1)})} holds. In fact,%
\begin{align*}
\ip{\varphi\oplus\varphi}{\pi_{\mathcal{K}}\left( f\right) \left
( \varphi\oplus\varphi\right) }  &  =\int_{\mathbb{T}}fh_{\mathcal{K}}\,d\mu\\
&  =\int_{\mathbb{T}}f\left(  \alpha_{\rho}\left(  h\right)  +\alpha_{\rho
^{2}}\left(  h\right)  \right)  \,d\mu\\
&  =\int_{\mathbb{T}}\left(  \alpha_{\rho^{-1}}\left(  f\right)  +\alpha
_{\rho^{-2}}\left(  f\right)  \right)  h\,d\mu\\
&  =\int_{\mathbb{T}}\left(  \alpha_{\rho^{2}}\left(  f\right)  +\alpha_{\rho
}\left(  f\right)  \right)  h\,d\mu\\
&  =\ip{\varphi}{\pi_{0}\left( \alpha_{\rho^{2}}\left( f\right) \right
) \varphi}+\ip{\varphi}{\pi_{0}\left( \alpha_{\rho}\left( f\right
) \right) \varphi},
\end{align*}
which is the assertion we made about the $\pi_{\mathcal{K}}$-representation of
$L^{\infty}\left(  \mathbb{T}\right)  $ on $\mathcal{K}=L^{2}\left(
\mathbb{R}\right)  \oplus L^{2}\left(  \mathbb{R}\right)  $.

It remains to prove the irreducibility assertion in
\textup{(\ref{RemWave.TwoIrrRepsClaim(2)}).} Each component in the sums
\[
U_{\mathcal{K}}  =U\oplus U
\]
and
\[
\pi_{\mathcal{K}}\left(  f\right)   =\pi_{0}\left(  \alpha_{\rho}\left(
f\right)  \right)  \oplus\pi_{0}\left(  \alpha_{\rho^{2}}\left(  f\right)
\right)
\]
may be viewed separately, but they are not representations of $\mathfrak{A}_{2}$
since $\alpha_{\rho}\left(  h\right)  $ is not an eigenfunction, i.e.,
$R\left(  \alpha_{\rho}\left(  h\right)  \right)  \neq\alpha_{\rho}\left(
h\right)  $; in fact, $R\left(  \alpha_{\rho}\left(  h\right)  \right)
=\alpha_{\rho^{2}}\left(  h\right)  $. But an operator in the commutant of
$\tilde{\pi}_{\mathcal{K}}=\left(  U_{\mathcal{K}},\pi_{\mathcal{K}}\right)
$, i.e., $Q\in\tilde{\pi}_{\mathcal{K}}\left(  \mathfrak{A}_{2}\right)  ^{\prime}%
$, $Q=\left(  Q_{ij}\right)  $, $1\leq i,j\leq2$, relative to the
decomposition $L^{2}\left(  \mathbb{R}\right)  \oplus L^{2}\left(
\mathbb{R}\right)  $ must have $Q_{ij}$ determined by the obvious matrix
identities, i.e., each $Q_{ij}$ commuting with $U$, $Q_{11}\pi_{0}\left(
\alpha_{\rho}\left(  f\right)  \right)  =\pi_{0}\left(  \alpha_{\rho}\left(
f\right)  \right)  Q_{11}$, $Q_{21}\pi_{0}\left(  \alpha_{\rho}\left(
f\right)  \right)  =\pi_{0}\left(  \alpha_{\rho^{2}}\left(  f\right)  \right)
Q_{21}$, etc. A direct check shows that $Q_{ij}=\delta_{ij}\lambda$ for
some $\lambda\in\mathbb{C}$. To see that $Q_{21}=0$, for example, use that
the stated representations of $L^{\infty}\left(  \mathbb{T}\right)  $ are
inequivalent. To see this, note that the functions $\alpha_{\rho}\left(
h\right)  $ and $\alpha_{\rho^{2}}\left(  h\right)  $ have different zero
sets: in fact,%
\[
\left\{  z\in\mathbb{T}\mid h\left(  \rho z\right)  =0\right\}  =\left\{
1,\rho\right\}  ,
\]
and
\[
\left\{  z\in\mathbb{T}\mid h\left(  \rho^{2}z\right)  =0\right\}  =\left\{
1,\rho^{2}\right\}  .\settowidth
{\qedskip}{$\displaystyle\left\{  z\in\mathbb{T}\mid h\left(  \rho^{2}%
z\right)  =0\right\}
=\left\{  1,\rho^{2}\right\}  .$}\addtolength{\qedskip}{-\textwidth}%
\rlap{\hbox
to-0.5\qedskip{\hfil\qed}}%
\]
\renewcommand{\qed}{}
\end{proof}

\begin{proof}
[Proof of Theorem \textup{\ref{ThmWave.5}}]We recall that $\mathfrak{A}_{N}$ is
the $C^{\ast}$-algebra generated by $L^{\infty}\left(  \mathbb{T}\right)  $
and a unitary symbol $U$ subject to the relation%
\begin{equation}
UfU^{-1}=f\left(  z^{N}\right)  ,\qquad f\in L^{\infty}\left(  \mathbb{T}%
\right)  , \label{eqWave.15}%
\end{equation}
and a representation $\tilde{\pi}\in\operatorname*{Rep}\left(  \mathfrak{A}%
_{N},\mathcal{H}\right)  $ amounts to a representation $\pi$ of $L^{\infty
}\left(  \mathbb{T}\right)  $ on $\mathcal{H}$, and a unitary operator, also
denoted $U$, on $\mathcal{H}$ such that%
\[
U\pi\left(  f\right)  U^{-1}=\pi\left(  f\left(  z^{N}\right)  \right)
,\qquad f\in L^{\infty}\left(  \mathbb{T}\right)  .
\]
In the computations for $\tilde{\pi}_{\varphi}$ and $\tilde{\pi}_{h}$ we use
the fact
that $\mathfrak{A}_{N}$ contains an increasing family of abelian subalgebras,
\emph{viz.,} $U^{-n}fU^{n}$, $n=0,1,\dots$. To see that the algebra for $n$ is
contained in the next one, use%
\[
U^{-n}fU^{n}=U^{-\left(  n+1\right)  }UfU^{-1}U^{n+1}=U^{-\left(  n+1\right)
}f\left(  z^{N}\right)  U^{n+1}.
\]
This structure also allows the interpretation of $\mathfrak{A}_{N}$ as an
inductive limit $C^{\ast}$-algebra; see \cite{BreJor} for more details.

\begin{lemma}
\label{LemWave.6}We have the estimate%
\begin{equation}
\left\|  \tilde{\pi}_{\varphi}\left(  A\right)  \varphi\right\|
_{L^{2}\left(  \mathbb{R}\right)  }\leq\sqrt{c}\left\|  \tilde{\pi}_{h}\left(
A\right)  \Phi\right\|  _{\mathcal{H}_{h}} \label{eqWave.16}%
\end{equation}
for all $A\in\mathfrak{A}_{N}$, where $\tilde{\pi}_{\varphi}\in\operatorname*{Rep}%
\left(  \mathfrak{A}_{N},L^{2}\left(  \mathbb{R}\right)  \right)  $ and
$\tilde{\pi}_{h}\in\operatorname*{Rep}\left(  \mathfrak{A}_{N},\mathcal{H}%
_{h}\right)  $ are the two cyclic representations described above, and where
$R\left(  h\right)  =h$, and $h_{\varphi}\leq ch$.
\end{lemma}

\begin{proof}
There is a dense $\ast$-subalgebra in $\mathfrak{A}_{N}$ consisting of finite sums%
\begin{equation}
\sum_{n\geq0}U^{-n}\alpha_{n}+\beta_{n}U^{n},\qquad\alpha_{0},\alpha_{1}%
,\dots,\beta_{0},\beta_{1},\dots\in L^{\infty}\left(  \mathbb{T}\right)  ,
\label{eqWave.17}%
\end{equation}
so we will estimate%
\[
\left\|  \tilde{\pi}_{h}\left(  A\right)  \Phi\right\|  _{\mathcal{H}_{h}}%
^{2}=\ip{\Phi}{\tilde{\pi}_{h}\left( A^{\ast}A\right) \Phi}%
\]
for general elements $A\in\mathfrak{A}_{N}$ of this form. In fact, we show that
there is a positive function $F_{A}$ in $L^{\infty}\left(  \mathbb{T}\right)
$ which only depends on $A$ such that%
\begin{equation}
\left\|  \tilde{\pi}_{h}\left(  A\right)  \Phi\right\|  _{\mathcal{H}_{h}}%
^{2}=\int_{\mathbb{T}}F_{A}h\,d\mu. \label{eqWave.18}%
\end{equation}
So this also applies to the wavelet representation $\tilde{\pi}_{\varphi}$,
and we get%
\[
\left\|  \tilde{\pi}_{\varphi}\left(  A\right)  \varphi\right\|
_{L^{2}\left(  \mathbb{R}\right)  }^{2}=\int_{\mathbb{T}}F_{A}h_{\varphi
}\,d\mu.
\]
Since $h_{\varphi}$ $\leq ch$, the estimate in the lemma will follow.

If $A\in\mathfrak{A}_{N}$ is given as above, we may compute $A^{\ast}A$ and find
the following expression (finite sum!):%
\begin{align*}
\ip{\Phi}{\tilde{\pi}_{h}\left( A^{\ast}A\right) \Phi}  &  =\left\|
\tilde{\pi}_{h}\left(  A\right)  \Phi\right\|  _{\mathcal{H}_{h}}^{2}\\
&  =\left\|  \pi\left(  \xi_{0}\right)  \Phi+U^{-1}\pi\left(  \xi_{1}\right)
\Phi+\dots+U^{-n}\pi\left(  \xi_{n}\right)  \Phi\right\|  _{\mathcal{H}_{h}%
}^{2},
\end{align*}
where $\xi_{0},\xi_{1},\dots,\xi_{n}\in L^{\infty}\left(  \mathbb{T}\right)
$. Introducing
\[
S_{0}\xi\left(  z\right)  =m_{0}\left(  z\right)  \xi\left(  z^{N}\right)
,\qquad\xi\in L^{\infty}\left(  \mathbb{T}\right)  ,
\]
the argument from the proof of Lemma \ref{LemPoof.7} yields%
\begin{align*}
\left\|  \sum_{k=0}^{n}U^{-k}\pi\left(  \xi_{k}\right)  \Phi\right\|
_{\mathcal{H}_{h}}^{2}  &  =\left\|  \sum_{k=0}^{n}\pi\left(  S_{0}^{k}%
\xi_{n-k}^{{}}\right)  \Phi\right\|  _{\mathcal{H}_{h}}^{2}\\
&  =\left\|  \sum_{k=0}^{n}S_{0}^{k}\xi_{n-k}^{{}}\right\|  _{L^{2}\left(
h\right)  }^{2}\\
&  =\int_{\mathbb{T}}\left|  \sum_{k=0}^{n}S_{0}^{k}\xi_{n-k}^{{}}\right|
^{2}h\,d\mu.
\end{align*}
The formula for $F_{A}$ may be read off from this, and the lemma follows.
\end{proof}

The lemma states that there is a well defined bounded operator $W\colon
\mathcal{H}_{h}\rightarrow L^{2}\left(  \mathbb{R}\right)  $ which is given by%
\[
W\tilde{\pi}_{h}\left(  A\right)  \Phi:=\tilde{\pi}_{\varphi}\left(  A\right)
\varphi,\qquad A\in\mathfrak{A}_{N};
\]
and it is clear that $W$ will intertwine the two representations, i.e., that%
\begin{equation}
W\tilde{\pi}_{h}\left(  A\right)  =\tilde{\pi}_{\varphi}\left(  A\right)
W,\qquad A\in\mathfrak{A}_{N}, \label{eqWave.20}%
\end{equation}
or $W\tilde{\pi}_{h}=\tilde{\pi}_{\varphi}W$ for short. But then $WW^{\ast}$
commutes with $\tilde{\pi}_{\varphi}$. Specifically,%
\[
WW^{\ast}\tilde{\pi}_{\varphi}=W\tilde{\pi}_{h}W^{\ast}=\tilde{\pi}_{\varphi
}WW^{\ast}.
\]
Since $\tilde{\pi}_{\varphi}$ is irreducible by Lemma \ref{Lemax+b.2}, we
conclude that $WW^{\ast}=\mathrm{const.}\,I_{L^{2}\left(  \mathbb{R}\right)
}$. It follows from the assumptions that the constant, $c$, say, is nonzero.
Hence $\left(  \sqrt{c}\right)  ^{-1}W^{\ast}$ is an isometry from
$L^{2}\left(  \mathbb{R}\right)  $ into $\mathcal{H}_{h}$ which intertwines,
i.e.,
\[
\tilde{\pi}_{h}\left(  \sqrt{c}\right)  ^{-1}W^{\ast}=\left(  \sqrt{c}\right)
^{-1}W^{\ast}\tilde{\pi}_{\varphi}.
\]
The assertion of the theorem follows from this.
\end{proof}

The argument from the theorem also gives the following corollary which
provides us with a numerical index for the convex cone of solutions $h\in
L^{1}\left(  \mathbb{T}\right)  $, $h\geq0$, $Rh=h$.

\begin{corollary}
\label{CorWave.7}Let $h$ and $h_{\varphi}$ be as described in Theorem
\textup{\ref{ThmWave.5};} then the intertwining operators%
\[
W\colon\mathcal{H}_{h}\longrightarrow L^{2}\left(  \mathbb{R}\right)
\]
form a Hilbert space, and the dimension of this Hilbert space is the
multiplicity of $\tilde{\pi}_{\varphi}$ in $\tilde{\pi}_{h}$.
\end{corollary}

\begin{proof}
The only argument in the proof of the corollary which is not already in the
theorem is the assertion of a Hilbert space structure on the intertwiners
$W\colon\mathcal{H}_{h}\rightarrow L^{2}\left(  \mathbb{R}\right)  $. If we
have two such, $W_{1}$, $W_{2}$, say, then $W_{1}^{{}}W_{2}^{\ast}$ commutes
with the wavelet representation $\tilde{\pi}_{\varphi}\in\operatorname*{Rep}%
\left(  \mathfrak{A}_{N},L^{2}\left(  \mathbb{R}\right)  \right)  $. Since the
latter is irreducible, $W_{1}^{{}}W_{2}^{\ast}$ must be a scalar times
$I_{L^{2}\left(  \mathbb{R}\right)  }$. Denoting this scalar $\left\langle
W_{1},W_{2}\right\rangle $, we have the desired inner product. We leave the
rest of the details to the reader.
\end{proof}

\begin{corollary}
\label{CorWave.8}Let $m_{0}=\frac{1}{\sqrt{2}\mathstrut}\left(  1+z^{3}%
\right)  $, and $h_{e}\left(  e^{-i\omega}\right)  =\frac{1}{9}\left(
1+2\cos\omega\right)  
^{2}
$. Since the wavelet representation in $L^{2}\left(
\mathbb{R}\right)  $ is induced by $h_{e}$ and irreducible, we conclude that
the only solutions $h$ to%
\[
0\leq h\leq\mathrm{const.}\,h_{e},\qquad R_{m_{0}}\left(  h\right)  =h,
\]
are of the form%
\[
h=\lambda h_{e},\qquad\lambda\geq0\text{.}%
\]
\end{corollary}

\begin{proof}
We saw in Example \ref{ExaWave.3(1)} that $h_{e}=\frac{1}{9}\left(
1+2\cos\omega\right)  ^{2}$ is the solution to $R\left(  h_{e}\right)  =h_{e}$
which corresponds to the wavelet representation for $m_{0}=\frac{1}{\sqrt
{2}\mathstrut}\left(  1+z^{3}\right)  $, and so it is irreducible by Lemma
\ref{Lemax+b.2}. But we saw that irreducibility of the representation
$\tilde{\pi}_{h}\in\operatorname*{Rep}\left(  \mathfrak{A}_{2},\mathcal{H}%
_{h}\right)  $ corresponds to extremality of the state $\omega_{h}$ on
$\mathfrak{A}_{2}$ which is given by%
\[
\omega_{h}\left(  fU^{n}\right)  =\int_{\mathbb{T}}fm_{0}^{\left(  n\right)
}h\,d\mu,
\]
where $m_{0}^{\left(  n\right)  }\left(  z\right)  =m_{0}^{{}}\left(
z\right)  m_{0}^{{}}\left(  z^{2}\right)  \cdots m_{0}^{{}}\left(  z^{2^{n-1}%
}\right)  $, $n=0,1,\dots$, $f\in L^{\infty}\left(  \mathbb{T}\right)  $. We
also saw that the estimate $\omega_{h}\leq\mathrm{const.}\,\omega_{h^{\prime}%
}$ on the positive elements in $\mathfrak{A}_{2}$ is equivalent to%
\[
\int_{\mathbb{T}}\left|  f\right|  ^{2}h\,d\mu\leq\mathrm{const.}%
\,\int_{\mathbb{T}}\left|  f\right|  ^{2}h^{\prime}\,d\mu,\qquad f\in
L^{\infty}\left(  \mathbb{T}\right)  .
\]
The conclusion of the corollary is immediate from this since we noted that the
wavelet representation is induced by $h_{e}$ and irreducible.
\end{proof}

The argument which we used in the proofs of Lemma \ref{LemWave.6} and
Corollary \ref{CorWave.7} yields a more general result about the commutant of
our representations of $\mathfrak{A}_{N}$. Let $m_{0}\in L^{\infty}\left(
\mathbb{T}\right)  $ and suppose $m_{0}$ satisfies identity (\ref{eqax+b.7}),
and moreover that it is non-singular. Our result will apply to wavelet
filters, but the other properties of wavelet filters will not be needed. Let
$R=R_{m_{0}}$ be the corresponding Ruelle operator, and consider a solution
$h$ to $R\left(  h\right)  =h$, $h\in L^{1}\left(  \mathbb{T}\right)  $,
$h\geq0$. Let $\tilde{\pi}\in\operatorname*{Rep}\left(  \mathfrak{A}%
_{N},\mathcal{H}_{h}\right)  $ be the cyclic representation which is induced
\emph{via} Theorem \ref{Thmax+b.3}, and let $\varphi\in\mathcal{H}_{h}$ denote
the cyclic vector. Recall that $\tilde{\pi}$ is determined by a representation
$\pi$ of $L^{\infty}\left(  \mathbb{T}\right)  $, and a unitary operator $U$,
on $\mathcal{H}_{h}$ such that%
\begin{align}
U\pi\left(  f\right)  U^{-1}  &  =\pi\left(  f\left(  z^{N}\right)  \right)
,\label{eqWave.21}\\
\ip{\varphi}{\pi\left( f\right) \varphi}  &  =\int_{\mathbb{T}}fh\,d\mu
,\label{eqWave.22}\\%
\intertext{and}%
U\varphi &  =\pi\left(  m_{0}\right)  \varphi. \label{eqWave.23}%
\end{align}

\begin{theorem}
\label{ThmWave.9}Let $\left(  m_{0},h\right)  $ be as described above. Then
there is a one-to-one correspondence between positive elements in the
commutant of $\tilde{\pi}\left(  \mathfrak{A}_{N}\right)  $ and solutions%
\begin{equation}
h_{Q}\in L^{1}\left(  \mathbb{T}\right)  ,\qquad R\left(  h_{Q}\right)
=h_{Q}, \label{eqWave.24}%
\end{equation}
satisfying the pointwise estimate%
\begin{equation}
0\leq h_{Q}\leq ch \label{eqWave.25}%
\end{equation}
for some constant $c$.
\end{theorem}

\begin{proof}
We introduce the notation $\mathcal{M}^{\prime}=\tilde{\pi}\left(
\mathfrak{A}_{N}\right)  ^{\prime}$ for the commutant, i.e., the bounded operators
$Q$ in $\mathcal{H}_{h}$ such that%
\begin{equation}
Q\tilde{\pi}\left(  A\right)  =\tilde{\pi}\left(  A\right)  Q,\qquad
A\in\mathfrak{A}_{N}. \label{eqWave.26}%
\end{equation}
Positivity of $Q$ means:
$\ip{\psi}{Q\psi}\geq0$, $\psi\in\mathcal{H}_{h}$. If
$Q$ is positive, the square root $Q^{\frac{1}{2}}$ is well defined by the
spectral theorem, and $Q^{\frac{1}{2}}$ is in $\mathcal{M}^{\prime}$ if $Q$ is.

Let $Q\in\mathcal{M}^{\prime}$ be given and positive. We then have%
\[
\ip{\varphi}{Q\tilde{\pi}\left( A^{\ast}A\right) \varphi}=\left\|  \tilde{\pi
}\left(  A\right)  Q^{\frac{1}{2}}\varphi\right\|  ^{2},\qquad A\in
\mathfrak{A}_{N}\text{,}%
\]
so $A\mapsto\ip{\varphi}{Q\tilde{\pi}\left( A\right) \varphi}$ is a positive
linear functional on $\mathfrak{A}_{N}$, and it will be denoted $\omega_{Q}$,
i.e.,%
\begin{equation}
\omega_{Q}\left(  A\right)  :=\ip{\varphi}{Q\tilde{\pi}\left( A\right
) \varphi},\qquad A\in\mathfrak{A}_{N}. \label{eqWave.27}%
\end{equation}
A standard estimate from operator theory (see, e.g., \cite{BrRoI}) yields the
estimate%
\begin{equation}
\omega_{Q}\left(  A^{\ast}A\right)  \leq\left\|  Q\right\|  \left\|
\tilde{\pi}\left(  A\right)  \varphi\right\|  ^{2} \label{eqWave.28}%
\end{equation}
Applying this to $A=\pi\left(  f\right)  $, $f\in L^{\infty}\left(
\mathbb{T}\right)  $, we see that the measure determined by $f\mapsto
\omega_{Q}\left(  \pi\left(  f\right)  \right)  $ is absolutely continuous.
Let $h_{Q}\in L^{1}\left(  \mathbb{T}\right)  $ be the Radon--Nikodym
derivative, i.e.,%
\[
\omega_{Q}\left(  \pi\left(  f\right)  \right)  =\int_{\mathbb{T}}fh_{Q}%
\,d\mu,
\]
where $\mu$ is Haar measure on $\mathbb{T}$. Setting $\varphi_{Q}=Q^{\frac
{1}{2}}\varphi$, we see that%
\begin{equation}
\ip{\varphi_{Q}}{\pi\left( f\right) \varphi_{Q}}=\int_{\mathbb{T}}fh_{Q}%
\,d\mu,\qquad f\in L^{\infty}\left(  \mathbb{T}\right)  . \label{eqWave.29}%
\end{equation}
Since $Q^{\frac{1}{2}}\in\mathcal{M}^{\prime}$, we also have%
\begin{equation}
U\varphi_{Q}=\pi\left(  m_{0}\right)  \varphi_{Q}. \label{eqWave.30}%
\end{equation}
Combining (\ref{eqWave.29}) and (\ref{eqWave.30}), we conclude that $R\left(
h_{Q}\right)  =h_{Q}$. Lemma \ref{LemPoof.3}, or Theorem \ref{Thmax+b.3},
yields that conclusion. In view of (\ref{eqWave.29}), a second application of
(\ref{eqWave.28}) yields the pointwise estimate%
\begin{equation}
0\leq h_{Q}\leq\left\|  Q\right\|  h, \label{eqWave.31}%
\end{equation}
concluding the proof in one direction.

Conversely, suppose $h_{Q}$ is given and satisfies (\ref{eqWave.24}%
)--(\ref{eqWave.25}) of the theorem. Since the representation $\tilde{\pi
}=\left(  \pi,U\right)  $ is cyclic, $\mathcal{H}_{h}$ is spanned (after
taking closure) by vectors of the form%
\[
U^{-n}\pi\left(  f\right)  \varphi,\qquad n=0,1,\dots,\qquad f\in L^{\infty}\left(
\mathbb{T}\right)  ,
\]
and we may define an operator $Q$ on $\mathcal{H}_{h}$ by the matrix entries%
\begin{multline}
\ip{U^{-n_{1}}\pi\left( f_{1}\right) \varphi}{QU^{-n_{2}}\pi\left( f_{2}%
\right) \varphi}\label{eqWave.31bis}\\
=%
\begin{cases}
\displaystyle\int_{\mathbb{T}}\overline{f_{1}\left( z^{N^{n_{2}-n_{1}}}%
\right) }f_{2}\left( z\right) \overline{m_{0}^{\left( n_{2}-n_{1}\right
) }\left( z\right
) }h_{Q}\left( z\right) \,d\mu\left( z\right) & \text{if }n_{2}\geq n_{1},  \\
\displaystyle
\vphantom{\raisebox{2pt}{$\displaystyle \int $}}
\int_{\mathbb{T}}\overline{f_{1}\left( z\right) }f_{2}%
\left( z^{N^{n_{1}-n_{2}}}\right) m_{0}^{\left( n_{1}-n_{2}\right) }%
\left( z\right) h_{Q}\left( z\right) \,d\mu\left( z\right)  & \text{if }%
n_{2}\leq n_{1},
\end{cases}%
\end{multline}
where $m_{0}^{\left(  k\right)  }\left(  z\right)  =m_{0}^{{}}\left(
z\right)  m_{0}^{{}}\left(  z^{N}\right)  \cdots m_{0}^{{}}\left(  z^{N^{k-1}%
}\right)  $. The argument from Lemma \ref{LemWave.6} shows that $Q$ is well
defined, and bounded. In fact, using (\ref{eqWave.31bis}) for $n_{1}=n_{2}=n$
and $f_{1}=f_{2}=f$, we get%
\[
\ip{U^{-n}\pi\left( f\right) \varphi}{QU^{-n}\pi\left( f\right) \varphi}%
=\int_{\mathbb{T}}\left|  f\right|  ^{2}h_{Q}\,d\mu\leq c\int_{\mathbb{T}%
}\left|  f\right|  ^{2}h\,d\mu=c\left\|  U^{-n}\pi\left(  f\right)
\varphi\right\|  _{\mathcal{H}_{h}}^{2},
\]
which is the desired boundedness. But the definition of the operator $Q$ also
entails that it is in the commutant, i.e., that it satisfies (\ref{eqWave.26}%
). The second implication is proved.
\end{proof}

\begin{corollary}
\label{CorWave.10}Let $m_{0}$ and $h$ be as described in Theorem
\textup{\ref{ThmWave.9},} and let $\tilde{\pi}_{h}$ be the corresponding
representation of $\mathfrak{A}_{N}$ on $\mathcal{H}_{h}$. Then the commutant
$\tilde{\pi}_{h}\left(  \mathfrak{A}_{N}\right)  ^{\prime}$ is abelian, and so
$\tilde{\pi}_{h}$ is the direct integral of a family of mutually inequivalent
irreducible representations of $\mathfrak{A}_{N}$.
\end{corollary}

\begin{proof}
The commutativity is a direct consequence of the formula $Q\mapsto h_{Q}$ in
Theorem \ref{ThmWave.9}, and, more specifically, of equation
(\ref{eqWave.31bis}), which expresses a fixed $Q\in\tilde{\pi}_{h}\left(
\mathfrak{A}_{N}\right)  ^{\prime}$ in terms of $h_{Q}$, where $R\left(
h_{Q}\right)  =h_{Q}$ and $h_{Q}\leq ch$. The assertion about $\tilde{\pi
}_{h}$ being a direct integral of a family of mutually inequivalent
irreducible representations follows from a standard fact in representation
theory: If $\tilde{\pi}_{h}$ contains two equivalent irreducibles, then the
commutant $\tilde{\pi}_{h}\left(  \mathfrak{A}_{N}\right)  ^{\prime}$ would
contain a copy of the $2$-by-$2$ complex matrices, or have such a copy of
$M_{2}\left(  \mathbb{C}\right)  $ in a direct integral.

If $m_{0}$ is continuous, then the eigenspace $\left\{  h\in L^{1}\left(
\mathbb{T}\right)  \mid R_{m_{0}}\left(  h\right)  =h\right\}  $ is
finite-dimensional by a theorem in \cite{CoRa90}, so in that case, $\tilde
{\pi}_{h}$ is a finite direct sum of a finite number of mutually inequivalent
irreducible representations. By a result in the next chapter, these
irreducible representations must be wavelet representations in the case when
$m_{0}$ is a wavelet filter which is non-singular.
\end{proof}

\section{\label{Cocy}Cocycle equivalence of filter functions}

Let $h\in L^{1}\left(  \mathbb{T}\right)  $, $h\ge 0$, be given, and form the
Hilbert space $L^{2}\left(  h\right)  :=L^{2}\left(  \mathbb{T},h\,d\mu
\right)  $ as usual. We saw in Lemma \ref{LemWave.2} that the operator
\begin{equation}
\left(  S_{0}\xi\right)  \left(  z\right)  :=m_{0}\left(  z\right)  \xi\left(
z^{N}\right)  \label{eqCocy.1}%
\end{equation}
is isometric in $L^{2}\left(  h\right)  $ if and only if $R\left(  h\right)
=h$ where $R=R_{m_{0}}$ is the Ruelle operator formed from $m_{0}$. We will
make the standing assumption that the function $m_{0}$ is a wavelet filter. It
follows that $R\left(  \openone\right)  =\openone$, so $S_{0}$ is an isometry
in $L^{2}\left(  \mathbb{T}\right)  $, which is the special case $h=\openone$.

\begin{definition}
\label{DefCocy.0}
We say that two wavelet filters $m_{0}^{{}}$ and $m_{0}^{\prime}$ are
\emph{cocycle equivalent} if there is a
nonzero
measurable function $f$ on
$\mathbb{T}$ such that%
\begin{equation}
f\left(  z^{N}\right)  m_{0}^{\prime}\left(  z\right)  =f\left(  z\right)
m_{0}\left(  z\right)  ,\qquad z\in\mathbb{T}. \label{eqCocy.2}%
\end{equation}
\end{definition}

It is immediate that $m_{0}^{{}}$ and $m_{0}^{\prime}$ satisfy (\ref{eqCocy.2})
for some $f$ if and only if the multiplication operator $\left(  M_{f}%
\xi\right)  \left(  z\right)  =f\left(  z\right)  \xi\left(  z\right)  $
intertwines the two isometries%
\begin{align}
S_{0}\xi\left(  z\right)   &  =m_{0}\left(  z\right)  \xi\left(  z^{N}\right)
\label{eqCocy.3}\\%
\intertext{and}%
S_{0}^{\prime}\xi\left(  z\right)   &  =m_{0}^{\prime}\left(  z\right)
\xi\left(  z^{N}\right)  \label{eqCocy.4}%
\end{align}
on $L^{2}\left(  \mathbb{T}\right)  $, i.e., if and only if%
\begin{equation}
M_{f}^{{}}S_{0}^{{}}=S_{0}^{\prime}M_{f}^{{}}. \label{eqCocy.5}%
\end{equation}

\begin{lemma}
\label{LemCocy.1}\ 

\begin{enumerate}
\item \label{LemCocy.1(1)}Let $f$ be a function on $\mathbb{T}$ which defines
a cocycle equivalence for wavelet filters $m_{0}^{{}}$ and $m_{0}^{\prime}$,
and set $h=\left|  f\right|  ^{2}$. Then%
\begin{equation}
R_{m_{0}}\left(  h\right)  =h. \label{eqCocy.6}%
\end{equation}

\item \label{LemCocy.1(2)}Conversely, if $m_{0}$ is given and $R_{m_{0}%
}\left(  h\right)  =h$, $h\geq0$,
nonzero,
then $f=\sqrt{h}$ defines a cocycle equivalence.
\end{enumerate}
\end{lemma}

\begin{proof}
\emph{Ad} (\ref{LemCocy.1(1)}): Let $h=\left|  f\right|  ^{2}$ be as described
in (\ref{LemCocy.1(1)}). Then
\begin{align*}
\left(  Rh\right)  \left(  z\right)   &  =\frac{1}{N}\sum_{w^{N}=z}\left|
m_{0}^{{}}\left(  w\right)  \right|  ^{2}h\left(  w\right) \\
&  =\frac{1}{N}\sum_{w^{N}=z}\left|  m_{0}^{\prime}\left(  w\right)  \right|
^{2}h\left(  w_{{}}^{N}\right) \\
&  =h\left(  z\right)  \underset{\equiv1}{\underbrace{\frac{1}{N}\sum
_{w^{N}=z}\left|  m_{0}^{\prime}\left(  w\right)  \right|  ^{2}}}=h\left(
z\right)  .
\end{align*}
Note that the normalization condition
$\frac{1}{N}\sum
_{w^{N}=z}\left|  m_{0}^{\prime}\left(  w\right)  \right|  ^{2}=1$
was needed only for the function
$m_{0}^{\prime}$,
not for the second filter function $m_{0}$ of Definition \ref{DefCocy.0}.

\emph{Ad} (\ref{LemCocy.1(2)}): Let $h$ be as in (\ref{LemCocy.1(2)}), and set
$f=\sqrt{h}$, and%
\begin{equation}
v\left(  z\right)  :=\left(  h\left(  z\right)  \right)  ^{\frac{1}{2}}%
m_{0}\left(  z\right)  . \label{eqCocy.7}%
\end{equation}
Then%
\[
\frac{1}{N}\sum_{w^{N}=z}\left|  v\left(  w\right)  \right|  ^{2}=\frac{1}%
{N}\sum_{w^{N}=z}h\left(  w\right)  \left|  m_{0}\left(  w\right)  \right|
^{2}=\left(  Rh\right)  \left(  z\right)  =h\left(  z\right)  .
\]
In particular,
\[
h\left(  z^{N}\right)  =\frac{1}{N}\sum_{k=0}^{N-1}\left|  v\left(
e^{i\frac{2\pi k}{N}}z\right)  \right|  ^{2};
\]
so if $h\left(  z^{N}\right)  =0$, then all $N$ terms on the right-hand side
must vanish. In particular, $v\left(  z\right)  =0$. It follows that the
function%
\begin{equation}
m_{0}^{\prime}\left(  z\right)  :=\frac{v\left(  z\right)  }{\left(  h\left(
z^{N}\right)  \right)  ^{\frac{1\mathstrut}{2}}} \label{eqCocy.8}%
\end{equation}
is well defined, and satisfies%
\[
\frac{1}{N}\sum_{w^{N}=z}\left|  m_{0}^{\prime}\left(  w\right)  \right|
^{2}=1,
\]
which is to say that $m_{0}^{\prime}$ is a wavelet filter. It is clear from
(\ref{eqCocy.7}) and (\ref{eqCocy.8}) that $f=\sqrt{h}$ does define a cocycle
equivalence between the two wavelet filters $m_{0}^{{}}$ and $m_{0}^{\prime}$
as claimed.
\end{proof}

\begin{definition}
\label{DefCocy.2}We say that a measurable function on $\mathbb{T}$ is
\emph{non-singular} if it does not vanish on a subset of $\mathbb{T}$ of
positive Lebesgue measure.
\end{definition}

The next result (Lemma \ref{LemCocy.3}) establishes
a general property of positive
harmonic functions defined from
a Ruelle operator, in the
non-singular case: and it is
an analogue of a classical
result for positive harmonic
functions defined in the usual way
from a Laplace operator; see,
e.g., \cite{Rud90} or \cite{SzFo70}.

\begin{lemma}
\label{LemCocy.3}If $m_{0}$ is a wavelet filter, and if $R_{m_{0}}h=h$, $h\in
L^{1}\left(  \mathbb{T}\right)  $, $h\geq0$, and $h\neq0$. Then $h$ is
non-singular whenever $m_{0}$ is.
\end{lemma}

\begin{proof}
Let $m_{0}$ and $h$ be as stated. Let $Z\left(  h\right)  $ be the complement
in $\mathbb{T}$ of the support of $h$.
We show that
$\mu \left( Z\left( h\right) \right) >0$ leads
to a contradiction.
Since $m_{0}$ is non-singular, and
$h\left(  z\right)  =\frac{1}{N}\sum_{w^{N}=z}\left|  m_{0}\left(  w\right)
\right|  ^{2}h\left(  w\right)  $, it follows that the set $Z^{\left(
1\right)  }\left(  h\right)  =\left\{  z\in\mathbb{T}\mid z^{N}\in
Z\left( h\right) \right\}  $
is contained in $Z\left(  h\right)  $ except possibly for a
zero-measure subset in $\mathbb{T}$. Moreover, $\mu\left(  Z^{\left(
1\right)  }\left(  h\right)  \right)  =\mu\left(  Z\left(  h\right)  \right)
$. Suppose now that $\mu\left(  Z\left(  h\right)  \right)  >0$, i.e., that
$h$ is not non-singular. Then continue recursively, defining%
\begin{equation}
Z^{\left(  n\right)  }\left(  h\right)  =\left\{  z\in\mathbb{T}\mid z^{N^{n}%
}\in Z\left(  h\right)  \right\}  . \label{eqCocy.9}%
\end{equation}
We get the inclusions
\[
Z^{\left( n+1\right) }\left( h\right) \subset 
Z^{\left( n\right) }\left( h\right) \subset 
\dots \subset 
Z\left( h\right) \subset 
\mathbb{T},
\]
except possibly for points of a subset of zero measure,
and $\mu\left(  Z^{\left(  n\right)  }\left(  h\right)  \right)  =\mu\left(
Z\left(  h\right)  \right)  >0,\qquad n=1,2,\dots$. Defining $Z^{\left(
\infty\right)  }\left(  h\right)  :=\bigcap_{n}Z^{\left(  n\right)  }\left(
h\right)  $, we note that $\mu\left(  Z^{\left(  \infty\right)  }\left(
h\right)  \right)  =\mu\left(  Z\left(  h\right)  \right)  $. Since
$Z^{\left(  \infty\right)  }\left(  h\right)  $ is invariant under $z\mapsto
z^{N}$, we conclude from ergodicity that $\mu\left(  Z^{\left(  \infty\right)
}\left(  h\right)  \right)  =\mu\left(  Z\left(  h\right)  \right)  =1$
relative to normalized Haar measure $\mu$ on $\mathbb{T}$. But this
contradicts the assumption $h\neq0$ in $L^{1}\left(  \mathbb{T}\right)  $, and
the proof is concluded.
\end{proof}

Assume that $m_{0}$ is non-singular. If $h\in L^{1}\left(  \mathbb{T}\right)
$, $h\geq0$, and $h\neq0$ solves $R_{m_{0}}h=h$, then $h$ will also be
non-singular by the lemma. It follows that the multiplication operator
$\Delta:=M_{h^{\frac{1}{2}}}$ is then an isometric isomorphism between the
Hilbert spaces $L^{2}\left(  h\right)  =L^{2}\left(  \mathbb{T},h\,d\mu
\right)  $ and $L^{2}\left(  \mathbb{T}\right)  $, i.e., we have, for $\xi\in
L^{2}\left(  h\right)  $,
\begin{equation}
\left\|  \Delta\xi\right\|  _{L^{2}\left(  \mathbb{T}\right)  }=\left\|
\xi\right\|  _{L^{2}\left(  h\right)  },\quad L^{2}\left(  h\right)
\overset{\Delta}{\longrightarrow}L^{2}\left(  \mathbb{T}\right)  ,\text{\quad
and\quad}\Delta\left(  L^{2}\left(  h\right)  \right)  =L^{2}\left(
\mathbb{T}\right)  . \label{eqCocy.10}%
\end{equation}

\begin{lemma}
\label{LemCocy.4}Let $\left(  S_{0}\xi\right)  \left(  z\right)  =m_{0}\left(
z\right)  \xi\left(  z^{N}\right)  $ and suppose $S_{0}$ is isometric in
$L^{2}\left(  h\right)  $, and further that $m_{0}$ is non-singular. Then
\begin{equation}
\left(  \Delta S_{0}\Delta^{-1}\xi\right)  \left(  z\right)  =m_{0}\left(
z\right)  \left(  \frac{h\left(  z\right)  }{h\left(  z^{N}\right)  }\right)
^{\frac{1}{2}}\xi\left(  z^{N}\right)  \label{eqCocy.11}%
\end{equation}
is an isometry in $L^{2}\left(  \mathbb{T}\right)  $. Moreover, the two
conditions%
\begin{equation}
\sum_{w^{N}=z}\left|  m_{0}\left(  w\right)  \right|  ^{2}=N \label{eqCocy.12}%
\end{equation}
and%
\[
m_{0}\left(  1\right)  =\sqrt{N}%
\]
are preserved under the cocycle operation%
\[
m_{0}\longmapsto m_{0}\left(  z\right)  \left(  \frac{h\left(  z\right)
}{h\left(  z^{N}\right)  }\right)  ^{\frac{1}{2}}.
\]
\end{lemma}

\begin{proof}
The proof follows essentially from the previous lemmas, but we recall that
$m_{0}$ is assumed to be continuous near $z=1$, and if $R_{m_{0}}\left(
h\right)  =h$, $h\in L^{1}\left(  \mathbb{T}\right)  $, $h\geq0$, it can be
checked that $h$ will then also be continuous near $z=1$ so that the
evaluation of
\begin{equation}
m_{0}^{\left(  h\right)  }\left(  z\right)  =m_{0}\left(  z\right)  \left(
\frac{h\left(  z\right)  }{h\left(  z^{N}\right)  }\right)  ^{\frac{1}{2}}
\label{eqCocy.13}%
\end{equation}
at $z=1$ is well defined.

It remains to check that $m_{0}^{\left(  h\right)  }$ satisfies
(\ref{eqCocy.12}). But
\begin{align*}
\frac{1}{N}\sum_{w^{N}=z}\left|  m_{0}^{\left(  h\right)  }\left(  w\right)
\right|  ^{2}  &  =\frac{1}{N}\sum_{w^{N}=z}\left|  m_{0}\left(  w\right)
\right|  ^{2}\frac{h\left(  w\right)  }{h\left(  w^{N}\right)  }\\
&  =\frac{1}{h\left(  z\right)  }\frac{1}{N}\sum_{w^{N}=z}\left|  m_{0}\left(
w\right)  \right|  ^{2}h\left(  w\right) \\
&  =\frac{1}{h\left(  z\right)  }\left(  Rh\right)  \left(  z\right)
=\frac{h\left(  z\right)  }{h\left(  z\right)  }=1,
\end{align*}
which completes the proof. Recall $R\left(  h\right)  =h$ holds since $S_{0}$
is isometric in $L^{2}\left(  h\right)  $.
\end{proof}

If $h\in L^{1}\left(  \mathbb{T}\right)  $, $h\geq0$,
$h\neq 0$,
is given and if $m_{0}$
is a non-singular wavelet filter, we showed that $S_{0}$ is an isometry in
$L^{2}\left(  h\right)  $ if and only if $R_{m_{0}}\left(  h\right)  =h$, and
we will now show that this isometry is pure, i.e., that $\bigcap_{n=1}%
^{\infty}S_{0}^{n}\left(  L^{2}\left(  h\right)  \right)  =\left\{  0\right\}
$. The proof is based on the previous lemmas and one of the main results in
\cite{BrJo97b}.

\begin{theorem}
\label{ThmCocy.5}Let $m_{0}$ and $h$ be as described above, i.e., $m_{0}$ a
non-singular wavelet filter and $R_{m_{0}}\left(  h\right)  =h$, $h\in
L^{1}\left(  \mathbb{T}\right)  $, $h\geq0$,
$h\neq 0$. Then $S_{0}\xi\left(  z\right)
=m_{0}\left(  z\right)  \xi\left(  z^{N}\right)  $, $\xi\in L^{2}\left(
h\right)  $, satisfies%
\begin{equation}
\bigcap_{n=1}^{\infty}S_{0}^{n}\left(  L^{2}\left(  h\right)  \right)
=\left\{  0\right\}  . \label{eqCocy.14}%
\end{equation}
If $U$ is the unitary operator in $\mathcal{H}_{h}$ with%
\[
U\pi\left(  f\right)  U^{-1}=\pi\left(  f\left(  z^{N}\right)  \right)
\]
which defines the $h$-induced cyclic representation, i.e., $U\varphi=\pi\left(
m_{0}\right)  \varphi$, and $\varphi$ denoting the cyclic vector, then%
\begin{equation}
\bigcap_{n=1}^{\infty}U^{n}\left(  V_{0}\left(  \varphi\right)  \right)
=\left\{  0\right\}  , \label{eqCocy.15}%
\end{equation}
where%
\begin{equation}
V_{0}\left(  \varphi\right)  :=\overline{\operatorname*{span}}\left\{
\pi\left(  f\right)  \varphi\mid f\in L^{\infty}\left(  \mathbb{T}\right)
\right\}  . \label{eqCocy.16}%
\end{equation}
\end{theorem}

\begin{proof}
We already showed that $\left\|  U^{n}\pi\left(  f\right)  \varphi\right\|
_{\mathcal{H}_{h}}=\left\|  S_{0}^{n}f\right\|  _{L^{2}\left(  h\right)  }
=\left\|  f\right\|  _{L^{2}\left(  h\right)  }$,
$n=0,1,\dots$, so it follows that the two intersection properties
(\ref{eqCocy.14}) and (\ref{eqCocy.15}) are equivalent. But, in view of Lemma
\ref{LemCocy.4}, we may check equivalently that
\[
\left(  S_{0}^{\left(  h\right)  }\xi\right)  \left(  z\right)  =m_{0}\left(
z\right)  \left(  \frac{h\left(  z\right)  }{h\left(  z^{N}\right)  }\right)
^{\frac{1}{2}}\xi\left(  z^{N}\right)  ,
\]
as an isometry in $L^{2}\left(  \mathbb{T}\right)  $, satisfies%
\begin{equation}
\bigcap_{n=1}^{\infty}\left(  S_{0}^{\left(  h\right)  }\right)  ^{n}%
L^{2}\left(  \mathbb{T}\right)  =\left\{  0\right\}  . \label{eqCocy.17}%
\end{equation}
Because of Lemma \ref{LemCocy.4}, this in turn is immediate from \cite[Theorem
3.1]{BrJo97b}, and the proof is concluded.
\end{proof}

The significance of the result is that $\left(  S_{0},L^{2}\left(  h\right)
\right)  $ will then be unitarily equivalent to $\left(  S,H^{2}\left(
\mathbb{T},\mathcal{K}\right)  \right)  $ with
\begin{equation}
\left(  SF\right)  \left(  z\right)  =zF\left(  z\right)  , \label{eqCocy.18}%
\end{equation}
$F\colon\mathbb{T}\rightarrow\mathcal{K}$ satisfying $F\left(  z\right)
=\sum_{n=0}^{\infty}z^{n}k_{n}$, and
\[
\int_{\mathbb{T}}\left\|  F\left(  z\right)  \right\|  _{\mathcal{K}}%
^{2}\,d\mu\left(  z\right)  =\sum_{n=0}^{\infty}\left\|  k_{n}\right\|
_{\mathcal{K}}^{2},
\]
where%
\[
k_{n}\in\mathcal{K}=\ker\left(  S_{0}^{\ast}\right)  =\left\{  \xi\in
L^{2}\left(  h\right)  \mid S_{0}^{\ast}\xi=0\right\}  .
\]
This is a simple application of the standard Wold decomposition
(see \cite{SzFo70})
for the
isometry $\left(  S_{0},L^{2}\left(  h\right)  \right)  $.

If $N=2$, set%
\begin{equation}
\left(  S_{1}\xi\right)  \left(  z\right)  :=z\,\overline{m_{0}\left(
-z\right)  }\left(  \frac{h\left(  -z\right)  }{h\left(  z\right)  }\right)
^{\frac{1}{2}}\xi\left(  z^{2}\right)  \text{,\qquad for }\xi\in L^{2}\left(
h\right)  . \label{eqCocy.19}%
\end{equation}
Then we have%
\begin{equation}
S_{i}^{\ast}S_{j}^{{}}=\delta_{ij}^{{}}I_{L^{2}\left(  h\right)  }^{{}%
}\text{,\qquad and\qquad}\sum_{i=0}^{1}S_{i}^{{}}S_{i}^{\ast}=I_{L^{2}\left(
h\right)  }^{{}} \label{eqCocy.20}%
\end{equation}
in view of Lemma \ref{LemCocy.1}. We therefore get $\mathcal{K}=\ker\left(
S_{0}^{\ast}\right)  =S_{1}^{{}}\left(  L^{2}\left(  h\right)  \right)  $.

The relations (\ref{eqCocy.20}) are called the Cuntz relations and correspond
to representations of the corresponding $C^{\ast}$-algebra $\mathcal{O}_{2}$
(see \cite{Cun77}) which is known to be simple. They are the representations
which act on $L^{2}\left(  h\right)  $, and \emph{via} $Wf=\pi\left(  f\right)
\varphi$, they intertwine with operators on $\mathcal{H}_{h}$, in particular
$WS_{0}f=UWf$, $f\in L^{2}\left(  h\right)  $. \emph{Via} Lemma
\ref{LemCocy.4}, they correspond to representations of $\mathcal{O}_{2}$ on
$L^{2}\left(  \mathbb{T}\right)  $, but these representations are different
from those of Lemma \ref{LemWave.1} and Lemma \ref{LemCocy.1}.

In fact, the representation on $L^{2}\left(  \mathbb{T}\right)  $ which
intertwines with (\ref{eqCocy.20}) is given by $\left(  S_{i}^{\left(
h\right)  }\xi\right)  \left(  z\right)  =m_{i}^{\left(  h\right)  }\left(
z\right)  \xi\left(  z^{2}\right)  $, $\xi\in L^{2}\left(  \mathbb{T}\right)
$, $i=0,1$, where
\begin{align*}
m_{0}^{\left(  h\right)  }\left(  z\right)   &  =m_{0}^{{}}\left(  z\right)
\left(  \frac{h\left(  z\right)  }{h\left(  z^{2}\right)  }\right)  ^{\frac
{1}{2}}\\%
\intertext{and}%
m_{1}^{\left(  h\right)  }\left(  z\right)   &  =z\,\overline{m_{0}^{{}%
}\left(  -z\right)  }\left(  \frac{h\left(  -z\right)  }{h\left(
z^{2}\right)  }\right)  ^{\frac{1}{2}},\qquad z\in\mathbb{T}.
\end{align*}

\begin{corollary}
\label{CorCocy.6}Let $m_{0}$ and $h$ be as above, and let $\left(
\pi,\mathcal{H}_{h},\varphi\right)  $ be the corresponding cyclic
representation. Then the operator given by $f\mapsto WS_{1}f$, from
$L^{2}\left(  h\right)  $ into $\mathcal{H}_{h}$ \textup{(}see
\textup{(\ref{eqCocy.19})),} maps $L^{2}\left(  h\right)  $ onto the space
$V_{0}\left(  \varphi\right)  \ominus U\left(  V_{0}\left(  \varphi\right)
\right)  $ where $V_{0}\left(  \varphi\right)  $ is the closed cyclic subspace
in $\mathcal{H}_{h}$ which is generated by $\varphi$ under the representation
$\pi\left(  L^{\infty}\left(  \mathbb{T}\right)  \right)  $.
\end{corollary}

\begin{proof}
The details are contained in the discussion above.
\end{proof}

\section{\label{Kean}The transfer operator of Keane}

Let $m_{0}$ be a
wavelet filter,
and let $R$ be the corresponding Ruelle operator. Then
we showed in \cite{BEJ97} and \cite{Jor98a} that, for $\lambda\in
\mathbb{T\,\diagdown}\left\{  1\right\}  $, i.e., $\lambda\in\mathbb{C}$,
$\left|  \lambda\right|  =1$, the eigenvalue problem $R\left(
h\right)  =\lambda h$ does not have nonzero solutions
$h$
in $L^{1}\left(
\mathbb{T}\right)  $ if the scaling function has orthogonal $\mathbb{Z}%
$-translates. In this chapter, we consider a more general framework which
admits such solutions, and which also includes problems in the theory of
iteration, other than the wavelet problems, e.g., iteration of conformal
transformations. The non-trivial solutions to $R\left(  h\right)  =\lambda h$
will be interpreted as functionals on a $C^{\ast}$-algebra analogous to
$\mathfrak{A}_{N}$, but they will not be positive if $\lambda\neq1$.

Let $\left(  X,\mathcal{B},\mu\right)  $ be a finite measure space, i.e.,
$\mu$ will be assumed to be \emph{finite} positive measure which is defined on
a $\sigma$-algebra of subsets of $X$. If $X$ is a topological space, we assume
that $\mathcal{B}$ includes the Borel subsets of $X$. Let $T\colon
X\rightarrow X$ be a measurable mapping of $X$ \emph{onto} $X$ which is
$N$-to-one, i.e., for $\mu$-a.a.\ $x$ in $X$, $T^{-1}\left(  x\right)
=\left\{  y\in X\mid Ty=x\right\}  $ is of cardinality $N$, and assume further
that $\mu$ is $T$-invariant, i.e., that%
\begin{equation}
\mu\left(  T^{-1}\left(  E\right)  \right)  =\mu\left(  E\right)  ,\qquad
E\in\mathcal{B}. \label{eqKean.1}%
\end{equation}
We shall further need a selection of measurable inverses $\sigma_{i}\colon
X\rightarrow X$, $i=1,\dots,N$, such that $T\left(  \sigma_{i}\left(
x\right)  \right)  =x$, $\mu$-a.a.\ $x$ in $X$, $i=1,\dots,N$, and such that%
\begin{equation}
\mu\left(  \sigma_{i}\left(  X\right)  \cap\sigma_{j}\left(  X\right)
\right)  =0 \label{eqKean.2}%
\end{equation}
for all $i\neq j$. The following invariance condition (which is slightly
stronger than (\ref{eqKean.1})) will be needed throughout the chapter. It may
be stated in the following three equivalent forms:%
\begin{align}
\mu &  =\frac{1}{N}\sum_{i=1}^{N}\mu\circ\sigma_{i}^{-1},\label{eqKean.3}\\
\mu\left(  E\right)   &  =\frac{1}{N}\sum_{i=1}^{N}\mu\left(  \sigma_{i}%
^{-1}\left(  E\right)  \right)  ,\qquad E\in\mathcal{B},\label{eqKean.4}\\
\int_{X}f\,d\mu &  =\frac{1}{N}\sum_{i=1}^{N}\int_{X}f\circ\sigma_{i}%
\,d\mu\text{\qquad for all }\mathcal{B}\text{-measurable }f\text{ on }X.
\label{eqKean.5}%
\end{align}

\begin{example}
\label{ExaKean.1}Let $X=\mathbb{T}$, $T\left(  z\right)  =z^{N}$, and let
$\left\{  \sigma_{i}\right\}  $ be the choice of $N$-th roots specified by,
e.g.,%
\begin{equation}
\sigma_{k}\left(  e^{-i\omega}\right)  =e^{-\frac{i\omega k}{N}},\qquad
k=0,1,\dots,N-1; \label{eqKean.6}%
\end{equation}
and let $\mu$ be Haar measure on $\mathbb{T}$, i.e., $\frac{1}{2\pi}\int
_{-\pi}^{\pi}\cdots\,d\omega$. Then it is easy to check that conditions
\textup{(\ref{eqKean.1})--(\ref{eqKean.5})} hold.
\end{example}

\setcounter{figure}{2}\begin{figure}[ptb]\begin{minipage}{285pt}\mbox{\includegraphics
[bbllx=49bp,bblly=0bp,bburx=239bp,bbury=288bp,height=432pt]
{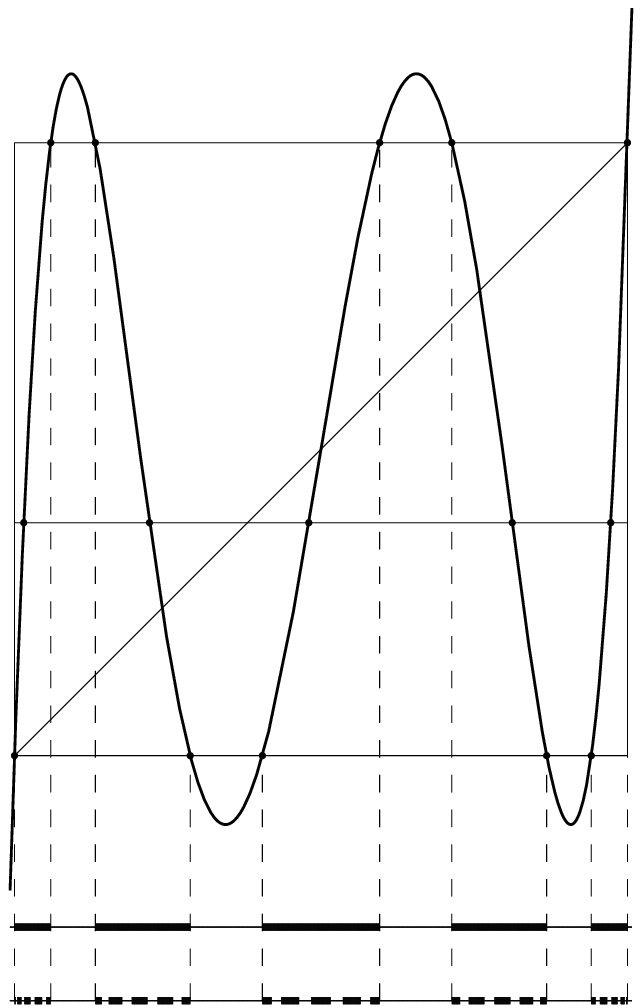}\llap{\setlength{\unitlength}{63.53pt}\begin{picture}(4.48,6.8)(-2.24,-3.797)
\put(-2.02125,-0.5){\makebox(0,0)[tr]{\raisebox
{-10pt}{$\sigma_{5}(x)$}$\nearrow$\hskip1pt\vphantom{\raisebox{1pt}{$\nearrow$}}}}
\put(-1.16642,-0.5){\makebox(0,0)[bl]{\vphantom{\raisebox
{-2pt}{$\swarrow$}}\hskip2pt$\swarrow$\raisebox{10pt}{$\sigma_{4}(x)$}}}
\put(-0.08278,-0.5){\makebox(0,0)[tr]{\raisebox
{-10pt}{$\sigma_{3}(x)$}$\nearrow$\hskip2pt\vphantom{\raisebox{2pt}{$\nearrow$}}}}
\put(1.30036,-0.5){\makebox(0,0)[bl]{\vphantom{\raisebox
{-2pt}{$\swarrow$}}\hskip2pt$\swarrow$\raisebox{10pt}{$\sigma_{2}(x)$}}}
\put(1.97008,-0.5){\makebox(0,0)[tr]{\raisebox
{-10pt}{$\sigma_{1}(x)$}$\nearrow$\hskip2pt\vphantom{\raisebox{2pt}{$\nearrow$}}}}
\put(-1.96088,-3.2){\makebox(0,0)[b]{$\sigma_{5}I$}}
\put(-1.21189,-3.2){\makebox(0,0)[b]{$\sigma_{4}I$}}
\put(0,-3.2){\makebox(0,0)[b]{$\sigma_{3}I$}}
\put(1.21189,-3.2){\makebox(0,0)[b]{$\sigma_{2}I$}}
\put(1.96088,-3.2){\makebox(0,0)[b]{$\sigma_{1}I$}}
\put(-0.61388,1.41665){\makebox(0,0){$I\times I$}}
\put(2.17776,-0.5){\makebox(0,0)[l]{$x\in I$}}
\put(2.17776,-3.25){\makebox(0,0)[l]{$\bigcup\limits_{i}^{5}\sigma_{i}I$}}
\put(2.17776,-3.75){\makebox(0,0)[l]{$\bigcup\limits_{i}^{5}\bigcup\limits_{j}^{5}
\sigma_{i}\sigma_{j}I$}}
\put(2.17776,-4){\makebox(0,0)[l]{\hskip3pt$\vdots$}}
\put(-1.96088,-4){\makebox(0,0){$\vdots$}}
\put(-1.21189,-4){\makebox(0,0){$\vdots$}}
\put(0,-4){\makebox(0,0){$\vdots$}}
\put(1.21189,-4){\makebox(0,0){$\vdots$}}
\put(1.96088,-4){\makebox(0,0){$\vdots$}}
\end{picture}}}\end{minipage}\vspace{14pt}%
\caption{{Example \ref{ExaKean.2}, case 1, $p(a)=a$:
$p(t)=t^{5}-5\alpha^{2}t^{3}+5\alpha^{4}t$, $\alpha =1.05$,
$I:=[a,b]$, $a=-b$, $b\approx 2.08411$.
(Figure designed by Brian Treadway.)}}%
\label{FigExactlyN}%
\end{figure}\begin{figure}[ptb]\begin{minipage}{285pt}\mbox{\includegraphics
[bbllx=49bp,bblly=0bp,bburx=239bp,bbury=288bp,height=432pt]
{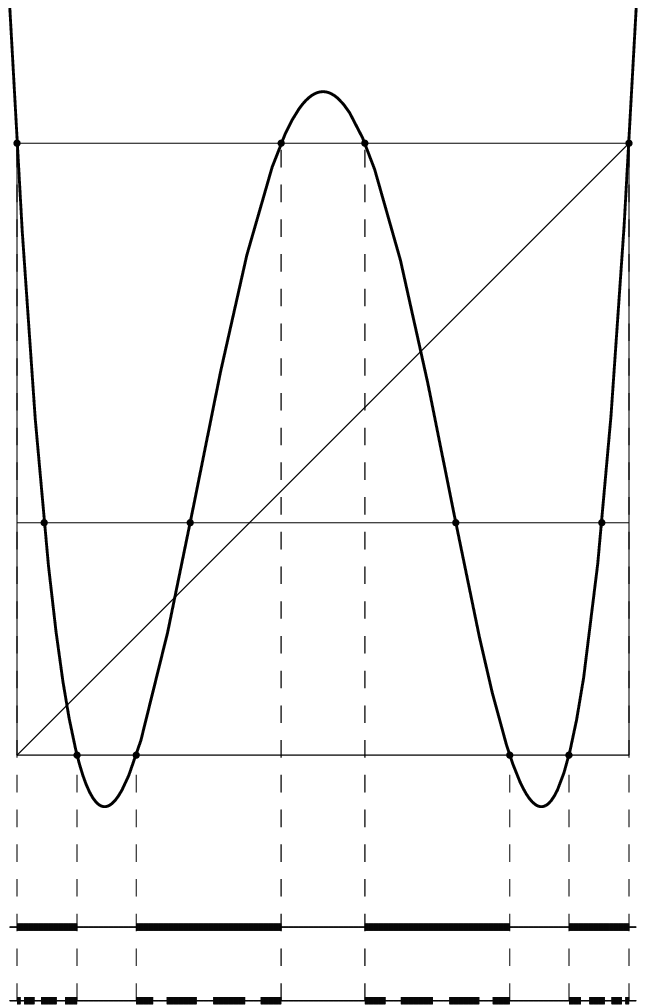}\llap{\setlength{\unitlength}{63.53pt}\begin{picture}(4.48,6.8)(-2.24,-3.797)
\put(-1.89595,-0.5){\makebox(0,0)[bl]{\vphantom{\raisebox
{-2pt}{$\swarrow$}}\hskip2pt$\swarrow$\raisebox{10pt}{$\sigma_{4}(x)$}}}
\put(-0.90299,-0.5){\makebox(0,0)[tr]{\raisebox
{-10pt}{$\sigma_{3}(x)$}$\nearrow$\hskip2pt\vphantom{\raisebox{2pt}{$\nearrow$}}}}
\put(0.90299,-0.5){\makebox(0,0)[bl]{\vphantom{\raisebox
{-2pt}{$\swarrow$}}\hskip2pt$\swarrow$\raisebox{10pt}{$\sigma_{2}(x)$}}}
\put(1.89595,-0.5){\makebox(0,0)[tr]{\raisebox
{-10pt}{$\sigma_{1}(x)$}$\nearrow$\hskip2pt\vphantom{\raisebox{2pt}{$\nearrow$}}}}
\put(-1.87652,-3.2){\makebox(0,0)[b]{$\sigma_{4}I$}}
\put(-0.77728,-3.2){\makebox(0,0)[b]{$\sigma_{3}I$}}
\put(0.77728,-3.2){\makebox(0,0)[b]{$\sigma_{2}I$}}
\put(1.87652,-3.2){\makebox(0,0)[b]{$\sigma_{1}I$}}
\put(-1.28286,1.41665){\makebox(0,0){$I\times I$}}
\put(2.17776,-0.5){\makebox(0,0)[l]{$x\in I$}}
\put(2.17776,-3.25){\makebox(0,0)[l]{$\bigcup\limits_{i}^{4}\sigma_{i}I$}}
\put(2.17776,-3.75){\makebox(0,0)[l]{$\bigcup\limits_{i}^{4}\bigcup\limits_{j}^{4}
\sigma_{i}\sigma_{j}I$}}
\put(2.17776,-4){\makebox(0,0)[l]{\hskip3pt$\vdots$}}
\put(-1.87652,-4){\makebox(0,0){$\vdots$}}
\put(-0.77728,-4){\makebox(0,0){$\vdots$}}
\put(0.77728,-4){\makebox(0,0){$\vdots$}}
\put(1.87652,-4){\makebox(0,0){$\vdots$}}
\end{picture}}}\end{minipage}\vspace{14pt}%
\caption{{Example \ref{ExaKean.2}, case 2, $p(a)=b$:
$p(t)=t^{4}-4\alpha^{2}t^{2}+2\alpha^{4}$, $\alpha =1.05$,
$I:=[a,b]$, $a=-b$, $b\approx 2.08064$.
(Figure designed by Brian Treadway.)}}%
\label{FigExactlyNeven}%
\end{figure}

\begin{example}
\label{ExaKean.2}\textup{(The }\emph{Julia set}\textup{ of a polynomial)} Let
$p\left(  z\right)  =z^{N}+p_{1}z^{N-1}+\dots+p_{N}$ be a polynomial with real
coefficients $p_{i}$, and leading coefficient $1$. When $z\in\overline
{\mathbb{C}}=\mathbb{C}\cup\left\{  \infty\right\}  $ is given, define
$z_{0}=z$, and $z_{n+1}=p\left(  z_{n}\right)  $. The fixed points are the
solutions $\zeta$ to $p\left(  \zeta\right)  =\zeta$, and of course
$\zeta=\infty$ is a fixed point. If $\zeta$ is a fixed point, the \emph{domain
of attraction} $\Omega\left(  \zeta\right)  $ is $\Omega\left(  \zeta\right)
=\left\{  z_{0}\in\overline{\mathbb{C}}\mid z_{n}\rightarrow\zeta\right\}  $,
and we say that $\zeta$ is \emph{attracting} if $\Omega\left(  \zeta\right)  $
contains a neighborhood of $\zeta$. Clearly $\zeta=\infty$ is attracting. The
Julia set $X=X\left(  p\right)  $ is the complement of the union of the
$\Omega\left(  \zeta\right)  $'s over all the fixed points; and it follows
that $X\left(  p\right)  $ is compact. For spectral theory, it is enough to
restrict to the case when $X\left(  p\right)  $ is contained in $\mathbb{R}$,
and it is known that the convex hull of $X\left(  p\right)  $ is of the form
$\left[  a,b\right]  $ where $a\in X\left(  p\right)  $, and $b$ is an
unstable fixed point. The equations $p\left(  y\right)  =a$, or $p\left(
y\right)  =b$, have exactly $N$ solutions in $\left[  a,b\right]  $
\textup{(}see Figures \textup{\ref{FigExactlyN}}
and \textup{\ref{FigExactlyNeven}),}
and so $p$ must have $N-1$
critical points in $\left[  a,b\right]  $. So if $x\in\left[  a,b\right]  $,
the equation $p\left(  y\right)  =x$ has $N$ solutions $y
_{i}
\in\left[
a,b\right]  $, $y
_{i}
=\sigma_{i}\left(  x\right)  $, $i=1,\dots,N$, $a\leq
\sigma_{N}\left(  x\right)  <\sigma_{N-1}\left(  x\right)  <\dots<\sigma
_{1}\left(  x\right)  \leq b$, and
\[
X\left(  p\right)  =\bigcap_{k=1}^{\infty}\bigcup\sigma_{i_{1}}\left(
\sigma_{i_{2}}\left(  \cdots\left(  \sigma_{i_{k}}\left(  \left[  a,b\right]
\right)  \right)  \cdots\right)  \right)  .
\]

Reasoning from the definitions, and insisting on
$X\left( p\right) \subset \mathbb{R}$, we see
that necessarily $p\left( b\right) =b$ with $b$ an unstable
fixed point, while for $a$ there
are only two possible cases, $p\left( a\right) =a$, or
$p\left( a\right) =b$. In the first case, $a$
is then also an unstable fixed
point.
In the real case, only $\infty $
is attracting. (Of course, the possibilities
are more varied in the case of complex
Julia sets.)

Defining $\mu$ by%
\begin{equation}
\mu\left(  \sigma_{i_{1}}\circ\sigma_{i_{2}}\circ\dots\circ\sigma_{i_{k}%
}\left(  \left[  a,b\right]  \right)  \right)  =N^{-k} \label{eqKean.7}%
\end{equation}
and extending to the Borel $\sigma$-algebra, we get a measure $\mu$ on
$X\left(  p\right)  $ which satisfies conditions \textup{(\ref{eqKean.1}%
)--(\ref{eqKean.5});} see, e.g., \cite{Bel92} and \cite{Bro65} for details.

Note that condition \textup{(\ref{eqax+b.poundbis})} is not satisfied in
this example and that $\mu $ is not absolutely
continuous.
\end{example}

Let $\left(  X,\mathcal{B},\mu\right)  $ be as described above, and let
$g\colon X\rightarrow\left[  0,1\right]  $ be a fixed measurable function
satisfying
\begin{align}
\sum_{i=1}^{N}g\left(  \sigma_{i}\left(  x\right)  \right)   &  =1\text{\qquad
a.a.\ }x\in X,\label{eqKean.8}\\%
\intertext{or, equivalently,}%
\sum_{Ty=x}g\left(  y\right)   &  =1\text{\qquad a.a.\ }x\in X.
\label{eqKean.9}%
\end{align}
Following Keane \cite{Kea72}, we define a more general transfer operator as
follows:%
\begin{equation}
\left(  Rf\right)  \left(  x\right)  =\sum_{Ty=x}g\left(  y\right)  f\left(
y\right)  , \label{eqKean.10}%
\end{equation}
and we note that $R\left(  \openone\right)  =\openone$ when $\openone$ is the
constant function $1$ on $X$. It will be assumed that a measure $\mu$ is
chosen (and that it exists, see \cite{Kea72}) satisfying condition
(\ref{eqKean.3}), and we will normalize $\mu$ such that $\mu\left(  X\right)
=1$.

The paper \cite{JoPe98b} considers
the orthonormal basis problem
in $L^{2}\left(  \mu\right)  $ for the measure $\mu $
described above in (\ref{eqKean.7}). The question
of when a particular
orthogonal family of functions $\mathcal{F}$
in $L^{2}\left(  \mu\right)  $ is total is shown
to reduce to when the
eigenspace $\left\{  h\mid R_{g}\left(  h\right)  =h\right\}  $ is
one-dimensional. In that
application, $g=g_{\mathcal{F}}$ will depend on
the particular function family $\mathcal{F}$ to
be tested.

Let $\left(  X,\mathcal{B},T,\mu\right)  $ be as described above. Following
\cite{Kea72} we place some additional restrictions on the system:

\begin{itemize}
\item $\left(  X,d\right)  $ is a compact metric space;

\item $T$ is a local homeomorphism;

\item  there are $\rho,\delta\in\mathbb{R}_{+}$ such that $d\left(
Tx,Ty\right)  \geq\rho d\left(  x,y\right)  $ if $d\left(  x,y\right)
\leq\delta$;

\item  for each $\varepsilon\in\mathbb{R}_{+}$ there is $n_{\varepsilon}$ such
that $\bigcup_{n\geq n_{\varepsilon}}T^{-n}\left(  x_{0}\right)  $ is $\varepsilon
$-dense in $X$ for all $x_{0}\in X$.
\end{itemize}

A measure $\nu=\nu_{g}$ on $X$ is said to be a $g$-measure if
\[
\int_{X}\left(  R_{g}f\right)  \left(  x\right)  \,d\nu_{g}\left(
x\right)  =\int_{X}f\,d\nu_{g}%
\]
for all $f\in C\left(  X\right)  $. Keane showed that, under the above
restrictions, $g$-measures exist, but they are generally not unique. The
examples of $g$-measures that we shall need here are constructed as follows: A
finite set of points $x_{1},x_{2},\dots,x_{k}$ in $X$ is called a cycle if
$x_{i+1}=Tx_{i}$, $x_{1}=Tx_{k}$, and if $g\left(  x_{i}\right)  =1$,
$i=1,\dots,k$.

\begin{lemma}
\label{LemKean.gMeas}Let the $g$ system be given as above, and let
$x_{1},\dots,x_{k}$ be a cycle. Then $\nu:=\frac{1}{k}\sum_{i=1}^{k}%
\delta_{x_{i}}$ is a $g$-measure, where $\delta_{x_{i}}$ denotes the point
mass at $x_{i}$.
\end{lemma}

\begin{proof}
Let $f\in C\left(  X\right)  $, and let $\nu$ be as stated. Then
\[
\int_{X}R\left(  f\right)  \,d\nu=\frac{1}{k}\sum_{i}R\left(  f\right)
\left(  x_{i}\right)  =\frac{1}{k}\sum_{i}\sum_{Ty=x_{i}}g\left(  y\right)
f\left(  y\right)  .
\]
But $T^{-1}\left(  x_{i}\right)  $ includes $x_{i-1}$, and $g$ vanishes on all
other points in $T^{-1}\left(  x_{i}\right)  $. Hence%
\[
\int_{X}R\left(  f\right)  \,d\nu=\frac{1}{k}\sum_{i}g\left(  x_{i-1}\right)
f\left(  x_{i-1}\right)  =\frac{1}{k}\sum_{i}f\left(  x_{i-1}\right)
=\int_{X}f\,d\nu.\settowidth
{\qedskip}{$\displaystyle\int_{X}R\left(  f\right)  \,d\nu=\frac{1}{k}\sum
_{i}g\left(  x_{i-1}\right)
f\left(  x_{i-1}\right)  =\frac{1}{k}\sum_{i}f\left(  x_{i-1}\right)
=\int_{X}f\,d\nu.$}\addtolength{\qedskip}{-\textwidth}\rlap{\hbox
to-0.5\qedskip{\hfil\qed}}
\]
\renewcommand{\qed}{}
\end{proof}

\begin{lemma}
\label{LemKean.3}Returning to the general setup, let the operator $R$ be as in
\textup{(\ref{eqKean.10}).} Then $R$ leaves invariant all the spaces
$L^{p}\left(  X,\mu\right)  $, $1\leq p\leq\infty$, and we have%
\begin{equation}
N\int_{X}g\left(  x\right)  \xi\left(  Tx\right)  f\left(  x\right)
\,d\mu\left(  x\right)  =\int_{X}\xi\left(  x\right)  \left(  Rf\right)
\left(  x\right)  \,d\mu\left(  x\right)  \label{eqKean.11}%
\end{equation}
for all $\xi\in L^{\infty}\left(  X\right)  $, and $f\in L^{1}\left(
X,\mu\right)  $.
\end{lemma}

\begin{proof}
The argument from Chapter \ref{Poof} shows that $L^{p}\left(  X,\mu\right)  $
is invariant under $R$ for $p=1,2$ and $p=\infty$. The other cases follow from
a standard interpolation argument for $L^{p}$-norms. The proof of
(\ref{eqKean.11}) is a simple application of (\ref{eqKean.10}) and
(\ref{eqKean.5}). Specifically, we have:%
\begin{align*}
N\int_{X}g\left(  x\right)  \xi\left(  Tx\right)  f\left(  x\right)
\,d\mu\left(  x\right)   &  \underset{\text{(by (\ref{eqKean.5}))}}%
{\longequal}\int_{X}\sum_{Ty=x}g\left(  y\right)  \xi\left(  x\right)
f\left(  y\right)  \,d\mu\left(  x\right) \\
&  \underset{\text{(by (\ref{eqKean.10}))}}{\longequal}\int_{X}\xi\left(
x\right)  \left(  Rf\right)  \left(  x\right)  \,d\mu\left(  x\right)  ,
\end{align*}
which is the desired identity (\ref{eqKean.11}) of the lemma.
\end{proof}

It turns out that almost all the results of Chapters \ref{ax+b}--\ref{Cocy}
carry over, \emph{mutatis mutandis,} to the more general transfer operators of
Keane. But now, instead of the $C^{\ast}$-algebra $\mathfrak{A}_{N}$ from Chapter
\ref{ax+b}, we need the following one, $\mathfrak{A}\left(  X,T\right)  $: It is
generated by $L^{\infty}\left(  X\right)  $ and a single (abstract) unitary
element $U$ subject to the relation%
\begin{equation}
UfU^{-1}=f\circ T. \label{eqKean.12}%
\end{equation}
Note that a main difference between $L^{\infty}\left(  \mathbb{T}\right)  $,
which was used in $\mathfrak{A}_{N}$ from Chapter \ref{ax+b}, and the general
case, is that $L^{\infty}\left(  X\right)  $ is generally not singly
generated, so $\mathfrak{A}\left(  X,T\right)  $ does not have an equivalent
formulation as $UVU^{-1}=V^{N}$, and similarly, there is not a discrete group
formulation which is parallel to the one we used for the $ax+b$ group in
Chapter \ref{ax+b}.

\begin{definition}
\label{DefKean.4}A representation $\tilde{\pi}$ of $\mathfrak{A}\left(
X,T\right)  $ in a Hilbert space $\mathcal{H}$ is a pair $\left(
\pi,U\right)  $ where $\pi$ is a representation of $L^{\infty}\left(
X\right)  $ on $\mathcal{H}$, and $U$ a unitary operator $U\colon
\mathcal{H}\rightarrow\mathcal{H}$, i.e., $U^{\ast}=U^{-1}$, such that%
\begin{equation}
U\pi\left(  f\right)  U^{-1}=\pi\left(  f\circ T\right)  ,\qquad f\in
L^{\infty}\left(  X\right)  . \label{eqKean.13}%
\end{equation}
Note that $\mathfrak{A}\left(  X,T\right)  $ is the norm-closure of the linear
span of elements of the form $U^{-n}fU^{k}$ where $n,k\in\left\{
0,1,\dots\right\}  $, and $f\in L^{\infty}\left(  X\right)  $, and we set%
\[
\tilde{\pi}\left(  U^{-n}fU^{k}\right)  =U^{\ast\,n}\pi\left(  f\right)  U^{k}%
\]
with the slight abuse of notation, denoting by $U$ both an abstract element in
$\mathfrak{A}\left(  X,T\right)  $ and an operator in $\mathcal{H}$. The
representation is \emph{cyclic} if there is a vector $\varphi\in\mathcal{H}$
such that $\left\{  \tilde{\pi}\left(  A\right)  \varphi\mid A\in
\mathfrak{A}\left(  X,T\right)  \right\}  $ is norm-dense in $\mathcal{H}$.
\end{definition}

\begin{lemma}
\label{LemKean.5}Let $\left(  X,\mathcal{B},\mu,g\right)  $ be as described
above, and set $m_{0}=\sqrt{Ng}$, i.e.,%
\begin{equation}
m_{0}\left(  x\right)  =\sqrt{N}\left(  g\left(  x\right)  \right)  ^{\frac
{1}{2}},\qquad x\in X. \label{eqKean.14}%
\end{equation}
Let $\tilde{\pi}=\left(  \pi,U\right)  $ be a \emph{normal} representation of
$\mathfrak{A}\left(  X,T\right)  $ in $\mathcal{H}$, and let $\varphi
\in\mathcal{H}$ satisfy%
\begin{equation}
U\varphi=\pi\left(  m_{0}\right)  \varphi, \label{eqKean.15}%
\end{equation}
i.e., $U\varphi=\sqrt{N}\pi\left(  g^{\frac{1}{2}}\right)  \varphi$; then the
spectral density $h\in L^{1}\left(  X,\mu\right)  $ satisfies%
\begin{equation}
R\left(  h\right)  =h. \label{eqKean.16}%
\end{equation}
\textup{(}Recall the spectral density $h$ is the Radon--Nikodym derivative of
the measure $\nu_{\varphi}$ defined by $f\mapsto\ip{\varphi}{\pi\left
( f\right) \varphi}$ with respect to $\mu$, i.e., $h=\frac{d\nu_{\varphi
\mathstrut}}{d\mu}$. The absolute continuity of $\nu_{\varphi}$ with respect
to $\mu$ is part of the definition of a normal representation.\textup{)}
\end{lemma}

\begin{proof}
For arbitrary $f\in L^{\infty}\left(  X\right)  $, we have
\begin{align*}
\int_{X}fh\,d\mu &  =\ip{\varphi}{\pi\left( f\right) \varphi}_{\mathcal{H}}\\
&  =\ip{U\varphi}{U\pi\left( f\right) \varphi}_{\mathcal{H}}\\
&  =\ip{\pi\left( m_{0}\right) \varphi}{\pi\left( f\circ T\right) U\varphi
}_{\mathcal{H}}\\
&  =\ip{\pi\left( m_{0}\right) \varphi}{\pi\left( f\circ T\right) \pi
\left( m_{0}\right) \varphi}_{\mathcal{H}}\\
&  =N\ip{\varphi}{\pi\left( g\left( f\circ T\right) \right) \varphi}_{\mathcal{H}}\\
&  =N\int_{X}g\left( f\circ T\right) h\,d\mu\\
&  =\int_{X}fR\left(  h\right)  \,d\mu,
\end{align*}
where Lemma \ref{LemKean.3} was used in the last step. Since $L^{\infty
}\left(  X\right)  $ separates points in $L^{1}\left(  X,\mu\right)  $ the
desired identity $R\left(  h\right)  =h$ follows from this, and the proof is completed.
\end{proof}

Our main theorem below is a converse of this lemma.

\begin{theorem}
\label{ThmKean.6}Let $\left(  X,\mathcal{B},\mu,g\right)  $ be as described
above, i.e., we have a measure space $\left(  X,\mathcal{B},\mu\right)  $
where $\mu$ satisfies \textup{(\ref{eqKean.5})} and a function $g\colon
X\rightarrow\left[  0,1\right]  $ which satisfies the normalization
\textup{(\ref{eqKean.8})} of Keane for some $N\geq2$, and let $R$ be the
corresponding transfer operator, see \textup{(\ref{eqKean.10}).}

Then there is a one-to-one correspondence between \textup{(\ref{ThmKean.6(1)}%
)} and \textup{(\ref{ThmKean.6(2)})} below:\renewcommand{\theenumi}%
{\alph{enumi}}

\begin{enumerate}
\item \label{ThmKean.6(1)}$h\in L^{1}\left(  X,\mu\right)  $, $h\geq0$,
$R\left(  h\right)  =h$, and

\item \label{ThmKean.6(2)}positive linear functionals $\omega_{h}$ on
$\mathfrak{A}\left(  X,T\right)  $ such that%
\begin{subequations}
\label{eqKean.17}
\renewcommand{\theequation}{\theparentequation\roman{equation}}
\begin{align}
\omega_{h}\left(  fU^{n}\right)    & =N^{\frac{n}{2}}\int_{X}f\sqrt{g^{\left(
n\right)  }}\,h\,d\mu\label{eqKean.17(1)}\\
\intertext{and}
\omega_{h}\left(  U^{-1}fU^{n}\right)    & =N^{\frac{n-1}{2}}\int
_{X}\sqrt{g^{\left(  n-1\right)  }}\,R\left(  fh\right)  \,d\mu,\qquad n\geq
1.\label{eqKean.17(2)}\end{align}
\end{subequations}%
\end{enumerate}

If $g$ is non-singular and if $\mu$ is $T$-ergodic, then the isometry $U$
induced from the Gelfand--Naimark--Segal \textup{(}GNS\/\textup{)} construction
applied to \textup{(\ref{eqKean.17})} is unitary as an operator in the Hilbert
space $\mathcal{H}_{h}$ of the GNS representation. The representation is
unique up to unitary equivalence.
\end{theorem}

\begin{proof}
The structure of the present proof is very close to that of Theorem
\ref{Thmax+b.3} and we refer to Chapter \ref{Poof} for details. Here we will
only sketch some points which are specific to the present more general
situation. By the terminology in (\ref{eqKean.17}), $g^{\left(  n\right)  }$ is%
\begin{equation}
g^{\left(  n\right)  }\left(  x\right)  =g\left(  x\right)  g\left(
Tx\right)  \cdots g\left(  T^{n-1}x\right)  , \label{eqKean.18}%
\end{equation}
and we have%
\[
N^{n}\int_{X}g^{\left(  n\right)  }\left(  \xi\circ T^{n}\right)  f\,d\mu
=\int_{X}\xi R^{n}\left(  f\right)  \,d\mu
\]
as a simple generalization of Lemma \ref{LemKean.3} ($\xi\in L^{\infty}\left(
X\right)  $, $f\in L^{1}\left(  X,\mu\right)  $).

In passing from (\ref{ThmKean.6(1)})$\rightarrow$(\ref{ThmKean.6(2)}),we adopt
the inductive limit construction for the Hilbert space $\mathcal{H}_{h}$ of
the representation. Once the representation of $\mathfrak{A}\left(  X,T\right)  $
is identified, it is clear that the linear functional $\omega_{h}$ in
(\ref{eqKean.17}) will be positive; indeed%
\[
\omega_{h}\left(  A^{\ast}A\right)  =\ip{\varphi}{\tilde{\pi}\left( A^{\ast
}A\right) \varphi}=\left\|  \tilde{\pi}\left(  A\right)  \varphi\right\|
^{2}\geq0,\qquad A\in\mathfrak{A}\left(  X,T\right)  .
\]

The inductive limit construction from Figures 1--2 in Chapter \ref{Poof} will
be briefly reviewed as it applies to the present case: Let%
\begin{equation}
\left(  L^{\infty}\left(  X\right)  ,n\right)  =\left\{  \left(  \xi,n\right)
\mid\xi\in L^{\infty}\left(  X\right)  \right\}  \label{eqKean.19}%
\end{equation}
for $n=0,1,\dots$, and let%
\begin{equation}
\left\|  \left(  \xi,n\right)  \right\|  _{\mathcal{H}_{h}}^{2}:=\int_{X}%
R^{n}\left(  \left|  \xi\right|  ^{2}h\right)  \,d\mu. \label{eqKean.20}%
\end{equation}
Then
as in (\ref{eqPoof.9})--(\ref{eqPoof.12}), we have that%
\begin{equation}
J\colon \left(  \xi,n\right)  \longmapsto\left(  \xi\circ T,n+1\right)
\label{eqKean.21}%
\end{equation}
is isometric from $\left(  L^{\infty}\left(  X\right)  ,n\right)  $ into
$\left(  L^{\infty}\left(  X\right)  ,n+1\right)  $, and if we define $U$ by%
\begin{align*}
U\left(  \xi,0\right)   &  =\sqrt{N}\left(  g^{\frac{1}{2}}\left(  \xi\circ
T\right)  ,0\right) \\%
\intertext{and}%
U\left(  \xi,n+1\right)   &  =\sqrt{N}\left(  \left(  g\circ T^{n}\right)
^{\frac{1}{2}}\xi,n\right)  ,
\end{align*}
then $U$ is isometric in the Hilbert space which results from the inductive
limit construction applied to (\ref{eqKean.21}). The argument from Chapter
\ref{Poof}
and Figures 1--2
also shows that $U$ is unitary, if $g$ is non-singular (i.e., does
not vanish on a subset $E\subset X$, $\mu\left(  E\right)  >0$), and if $T$ is
ergodic relative to $\mu$. The ergodicity is needed to guarantee that a
nonzero solution to (\ref{ThmKean.6(1)}) must automatically be non-singular if
$g$ is; see Chapter \ref{Poof}
and Lemma \ref{LemCocy.3}
for details.

In any case, even if $U$ is only isometric we have the formula%
\[
U\pi\left(  f\right)  =\pi\left(  f\circ T\right)  U
\]
when the representation $\pi$ of $L^{\infty}\left(  X\right)  $ is defined by%
\[
\pi\left(  f\right)  \left(  \xi,n\right)  =\left(  \left(  f\circ
T^{n}\right)  \xi,n\right)
\]
for $n=0,1,\dots$ and $f,\xi\in L^{\infty}\left(  X\right)  $.

Finally, the cyclic vector $\varphi\in\mathcal{H}_{h}$ (the inductive limit
Hilbert space) will be given by $\varphi=\left(  \openone,0\right)  \in\left(
L^{\infty}\left(  X\right)  ,0\right)  $, and by (\ref{eqKean.21}),
\[
\varphi=\left(  \openone,0\right)  \sim\left(  \openone,1\right)  \sim\left(
\openone,2\right)  \sim\cdots,
\]
so%
\begin{align*}
U\varphi &  =\sqrt{N}\pi\left(  g^{\frac{1}{2}}\right)  \varphi\\%
\intertext{and}%
\ip{\varphi}{\pi\left( f\right) \varphi}  &  =\int_{X}fh\,d\mu,
\end{align*}
which concludes the construction of the representation from $h$ as in
(\ref{ThmKean.6(1)}), i.e.,
\linebreak
$R\left(  h\right)  =h$.

It remains to prove that if a representation $\tilde{\pi}=\left(
\pi,U\right)  $ of $\mathfrak{A}\left(  X,T\right)  $ results from a positive
functional $\omega_{h}$ in (\ref{eqKean.17}), then
\begin{equation}
U\varphi=\pi\left(  m_{0}\right)  \varphi ,  \label{eqKean.22}%
\end{equation}
where $m_{0}:=\sqrt{N}g^{\frac{1}{2}}$ and $\varphi$ is the cyclic vector of
this GNS construction. Introducing%
\begin{align*}
m_{0}^{\left(  n\right)  }\left(  x\right)   &  =m_{0}\left(  x\right)
m_{0}\left(  Tx\right)  \cdots m_{0}\left(  T^{n-1}x\right) \\
&  =N^{\frac{n}{2}}\left(  g\left(  x\right)  g\left(  Tx\right)  \cdots
g\left(  T^{n-1}x\right)  \right)  ^{\frac{1}{2}}\\
&  =N^{\frac{n}{2}}\left(  g^{\left(  n\right)  }\left(  x\right)  \right)
^{\frac{1}{2}},
\end{align*}
we claim that
\begin{equation}
\ip{\tilde{\pi}\left( A\right) \varphi}{U\varphi}_{\mathcal{H}_{h}}=\ip
{\tilde{\pi}\left( A\right) \varphi}{\pi\left( m_{0}\right) \varphi
}_{\mathcal{H}_{h}} \label{eqKean.23}%
\end{equation}
for all $A\in\mathfrak{A}\left(  X,T\right)  $ if (\ref{eqKean.17}) is given. The
result follows from this and cyclicity, i.e., (\ref{eqKean.22}) must hold. In
checking (\ref{eqKean.23}), it is enough to consider $A=fU^{n}$, $n=0,1,\dots
$, $f\in L^{\infty}\left(  X\right)  $, and then
\[
\ip{\tilde{\pi}\left( fU^{n}\right) \varphi}{U\varphi}  
  =\int_{X}\overline{m_{0}^{\left(  n-1\right)  }}\,R\left(  \bar{f}h\right)
\,d\mu\text{,\qquad by (\ref{eqKean.17(2)}),}%
\]
while%
\begin{align*}
\ip{\tilde{\pi}\left( fU^{n}\right) \varphi}{\pi\left( m_{0}\right) \varphi}
&  =\int_{X}\bar{f}\,\overline{m_{0}^{\left(  n\right)  }}\,m_{0}h\,d\mu\\
&  =\int_{X}\left|  m_{0}\right|  ^{2}\,\overline{m_{0}^{\left(  n-1\right)
}\circ T}\,\bar{f}h\,d\mu\\
&  =\int_{X}\left|  m_{0}\left(  x\right)  \right|  ^{2}\bar{m}_{0}\left(
Tx\right)  \cdots\bar{m}_{0}\left(  T^{n-1}x\right)  \bar{f}\left(  x\right)
h\left(  x\right)  \,d\mu\left(  x\right) \\
&  =\int_{X}R^{\ast}\left(  \overline{m_{0}^{\left(  n-1\right)  }}\right)
\bar{f}h\,d\mu\\
&  =\int_{X}\overline{m_{0}^{\left(  n-1\right)  }}\,R\left(  \bar{f}h\right)
\,d\mu\text{,\qquad by (\ref{eqKean.17(1)}),}%
\end{align*}
which proves that the two sides in (\ref{eqKean.23}) are identical, and so the
desired (\ref{eqKean.22}) must hold. The proof is completed.
\end{proof}

\section{\label{R-ha}A representation theorem for $R$-harmonic functions}

While the result in this chapter may be formulated for the representations
which correspond to the general $N$-to-$1$ transformations $T\colon
X\rightarrow X$ of the previous chapter, we shall restrict attention here (for
simplicity) to the case from Chapters \ref{ax+b}--\ref{Poof} above, i.e.,
$X=\mathbb{T}$, and $T\colon z\mapsto z^{N}$ when $N\geq2$ is fixed. As in the
previous chapters we consider a fixed wavelet filter $m_{0}$ of order $N$ and
the corresponding Ruelle operator $R=R_{m_{0}}$. From Theorem \ref{Thmax+b.3},
we know that each solution, $h\in L^{1}\left(  \mathbb{T}\right)  $, $h\geq0$,
$Rh=h$, defines a representation $\left(  \pi,U\right)  $ on a Hilbert space
$\mathcal{H}$ with cyclic vector $\varphi$ such that $U\varphi=\pi\left(
m_{0}\right)  \varphi$. We show in this chapter that this $\mathcal{H}$ may be
taken to be the $L^{2}$-space $L^{2}\left(  K_{N},\nu\right)  $, where
$K_{N}=\left(  \Lambda_{N}\right)  \sphat{}$ (the Pontryagin compact dual of
$\Lambda_{N}=\mathbb{Z}\left[  \frac{1}{N}\right]  $), and where $\nu
=\nu\left(  m_{0},h\right)  $ is a measure on $K_{N}$, depending on $\left(
m_{0},h\right)  $, i.e., $\mathcal{H}\simeq L^{2}\left(  K_{N},\nu\right)  $,
in such a way that $\varphi$ is the constant function in $L^{2}\left(
K_{N},\nu\right)  $.
The construction of the measure
$\nu\left( m_{0},h\right) $ uses an inductive
limit procedure for subalgebras
of $\mathfrak{A}_{N}$ which is somewhat
analogous to
(but different from)
one used recently in \cite{Bre96},
\cite{Lac98}, \cite{Mur95}, and \cite{Sta93}.

In the proof of Theorem \ref{ThmWave.5}, we saw that $\mathfrak{A}_{N}$ contains
an abelian subalgebra which is generated by elements of the form
\begin{equation}
U^{-i}fU^{i},\qquad f\in C\left(  \mathbb{T}\right)  , \label{eqR-ha.1}%
\end{equation}
$i\in\left\{  0,1,2,\dots\right\}  $. Recall that $\mathfrak{A}_{N}$ is defined
from the relation%
\begin{equation}
UfU^{-1}=f\left(  z^{N}\right)  \label{eqR-ha.2}%
\end{equation}
on the generators $C\left(  \mathbb{T}\right)  $ and $U$.
Hence
\begin{equation}
U^{-i}fU^{i}=U^{-\left( i+1\right) }f\left(  z^{N}\right) U^{i+1}.
\label{eqR-ha.2bis}
\end{equation}

Let $\mathcal{A}%
_{N}\subset\mathfrak{A}_{N}$ be the (abelian) subalgebra which is generated by the
elements in (\ref{eqR-ha.1}), and let $K_{N}$ be the compact \emph{Gelfand
space} of $\mathcal{A}_{N}$, i.e., $\mathcal{A}_{N}\simeq C\left(
K_{N}\right)  $. In Chapter \ref{ax+b}, we introduced $\Lambda_{N}%
=\mathbb{Z}\left[  \frac{1}{N}\right]  $, and we consider $\Lambda_{N}$ as a
discrete abelian group. The corresponding dual compact group $\Lambda
_{N}\sphat$ (Pontryagin dual) consists of all characters on $\Lambda_{N}$,
i.e., all one-dimensional representations%
\begin{equation}
\chi\colon\Lambda_{N}\longrightarrow\mathbb{T}\text{\qquad such that }%
\begin{cases}
\displaystyle\chi\left( \lambda+\lambda^{\prime}\right) =\chi\left
( \lambda\right) \chi\left( \lambda^{\prime}\right) ,  \\
\displaystyle\chi\left( -\lambda\right) =\overline{\chi\left( \lambda\right
) },
\end{cases}%
\lambda,\lambda^{\prime}\in\Lambda_{N}, \label{eqR-ha.3}%
\end{equation}
with group operation%
\[
\left(  \chi\chi^{\prime}\right)  \left(  \lambda\right)  =\chi\left(
\lambda\right)  \chi^{\prime}\left(  \lambda\right)  .
\]

\begin{theorem}
\label{ThmR-ha.1}\ 

\begin{enumerate}
\item \label{ThmR-ha.1(1)}The Gelfand space $K_{N}$ of $\mathcal{A}_{N}$ is
the compact group $\Lambda_{N}\sphat$.

\item \label{ThmR-ha.1(2)}If $m_{0}$ is a non-singular wavelet filter, and
$h\in L^{1}\left(  \mathbb{T}\right)  $, $h\geq0$, solves $R_{m_{0}}\left(
h\right)  =h$, then there is a unique measure $\nu=\nu\left(  m_{0},h\right)
$ on $K_{N}\simeq\Lambda_{N}\sphat$ such that%
\begin{equation}
\int_{K_{N}}U^{-i}fU^{i}\,d\nu=\int_{\mathbb{T}}R^{i}\left(  fh\right)
\,d\mu, \label{eqR-ha.4}%
\end{equation}
with $\mu$ denoting the Haar measure on $\mathbb{T}$.

\item \label{ThmR-ha.1(3)}Let $e\left(  \frac{n}{N^{k}}\right)  $ be
identified with the $\left(  2\pi\right)  N^{k}$-periodic function
$x\mapsto\exp\left(  i\frac{nx}{N^{k}}\right)  $ on $\mathbb{R}$; then the
vectors%
\begin{equation}
\left\{  e\left(  \frac{n}{N^{k}}\right)  \biggm|n\in\mathbb{Z},\;k\in\left\{
0,1,2,\dots\right\}  \right\}  \label{eqR-ha.5}%
\end{equation}
span $L^{2}\left(  K_{N},\nu\right)  $, and the representation $\left(
\pi,U\right)  $ of $\mathfrak{A}_{N}$ which corresponds to $\left(  m_{0}%
,h\right)  $ is given by:%
\begin{equation}
\pi\left(  e_{n}\right)  e\left(  \lambda\right)  =e\left(  n+\lambda\right)
,\qquad\lambda\in\Lambda_{N}, \label{eqR-ha.6}%
\end{equation}%
\begin{equation}
Ue\left(  \frac{n}{N^{k+1}}\right)  =\sum_{j\in\mathbb{Z}}a_{j}e\left(
j+\frac{n}{N^{k}}\right)  , \label{eqR-ha.7}%
\end{equation}
where%
\[
m_{0}\left(  z\right)  =\sum_{j\in\mathbb{Z}}a_{j}e_{j}\left(  z\right)
=\sum_{j\in\mathbb{Z}}a_{j}z^{j}%
\]
and%
\begin{equation}
Ue_{0}=m_{0}=\sum_{n\in\mathbb{Z}}a_{n}e_{n}, \label{eqR-ha.8}%
\end{equation}
where $e_{0}=\openone$.
\end{enumerate}
\end{theorem}

\begin{proof}
\emph{Ad} (\ref{ThmR-ha.1(1)}): We introduced the function $e\left(  \frac
{n}{N^{k}}\right)  $ in (\ref{eqR-ha.5}) as a function on $\mathbb{R}$ of
period $\left(  2\pi\right)  N^{k}$, and by virtue of (\ref{eqR-ha.2}) it is
identified with $U^{-k}e_{n}U^{k}$. This corresponds to the case $f=e_{n}$ in
(\ref{eqR-ha.1}). From the theory of almost periodic functions (see, e.g.,
\cite{Rud90} or \cite{Bes55}), functions in $C\left(  K_{N}\right)  $ or
$L^{2}\left(  K_{N}\right)  $ may be identified with the corresponding
functions on $\mathbb{R}$ spanned by the frequencies from $\mathbb{Z}\left[
\frac{1}{N}\right]  $. For $C\left(  K_{N}\right)  $, the completion is in the
$\sup$-norm, and for $L^{2}\left(  K_{N}\right)  $ the completion is in the
norm which is defined as the limit, $T\rightarrow\infty$,%
\begin{equation}
\frac{1}{2T}\int_{-T}^{T}\left|  f\left(  x\right)  \right|  ^{2}\,dx.
\label{eqR-ha.9}%
\end{equation}
The dual of $K_{N}$ will be given by $\left\{  \frac{n}{N^{k}}\mid
n\in\mathbb{Z},\;k\in\left\{  0,1,2,\dots\right\}  \right\}  $ subject to the
following equivalence relation:%
\begin{equation}
\frac{n}{N^{k}}\sim\frac{l}{N^{i}}\text{\qquad if and only if }N^{i}n=N^{k}l.
\label{eqR-ha.10}%
\end{equation}
In defining operators on $e\left(  \frac{n}{N^{k}}\right)  $ we must then
check consistency with respect to the equivalence relation (\ref{eqR-ha.10}).
The conclusion in (\ref{ThmR-ha.1(1)}) follows if we check that $U^{-k}%
e_{n}U^{k}=U^{-i}e_{l}U^{i}$ when (\ref{eqR-ha.10}) holds. The proof of this
identity is based on (\ref{eqR-ha.1}),
(\ref{eqR-ha.2bis}),
and an induction. The first step is
consideration of $n$, $k$, $l$ such that $n=Nl$. But then $U^{-1}e_{n}%
U=U^{-1}e_{Nl}U=U^{-1}e_{l}\left(  z^{N}\right)  U=U^{-1}Ue_{l}U^{-1}U=e_{l}$,
where we used (\ref{eqR-ha.2}) in the second-to-last step. The induction is
left to the reader.

\emph{Ad} (\ref{ThmR-ha.1(2)}): We first check that the right-hand side in
formula (\ref{eqR-ha.4}) defines a linear functional $L$ on $\mathcal{A}_{N}$,
i.e.,%
\begin{equation}
L\left(  U^{-i}fU^{i}\right)  :=\int_{\mathbb{T}}R^{i}\left(  fh\right)
\,d\mu, \label{eqR-ha.11}%
\end{equation}
$i\in\left\{  0,1,2,\dots\right\}  $, $f\in C\left(  \mathbb{T}\right)  $.
Hence if $U^{-i}fU^{i}$ is represented in two ways,
as for example in (\ref{eqR-ha.2bis}),
we check that $L$ takes
the same value either way: that amounts to checking that%
\begin{equation}
\int_{\mathbb{T}}R^{i+1}\left(  f\left(  z^{N}\right)  h\right)  \,d\mu
=\int_{\mathbb{T}}R^{i}\left(  fh\right)  \,d\mu. \label{eqR-ha.12}%
\end{equation}
The other cases will then follow from this and induction. The proof of
(\ref{eqR-ha.12}) is based on the eigenvalue property, $Rh=h$. In fact,%
\begin{align*}
\int_{\mathbb{T}}R^{i+1}\left(  f\left(  z^{N}\right)  h\right)  \,d\mu &
=\int_{\mathbb{T}}R^{i}\left(  R\left(  f\left(  z^{N}\right)  h\right)
\right)  \,d\mu\\
&  =\int_{\mathbb{T}}R^{i}\left(  fRh\right)  \,d\mu\\
&  =\int_{\mathbb{T}}R^{i}\left(  fh\right)  \,d\mu,
\end{align*}
which is (\ref{eqR-ha.12}). Note that $Rh=h$ was used in the last step.

We now specialize to $f=e_{n}$, and define%
\begin{equation}
L\left(  \frac{n}{N^{k}}\right)  :=\int_{\mathbb{T}}R^{k}\left(
e_{n}h\right)  \,d\mu. \label{eqR-ha.13}%
\end{equation}
It follows that $L$ may be viewed as a function on $\Lambda_{N}=\mathbb{Z}%
\left[  \frac{1}{N}\right]  $, and we must check that it is positive definite.
We will then get the desired measure $\nu=\nu\left(  m_{0},h\right)  $ as a
solution to the corresponding \emph{moment problem.} We claim that, for all
finite sequences $\lambda\rightarrow A\left(  \lambda\right)  $ (i.e., at most
a finite number of nonzero scalar terms), we have%
\begin{equation}
\sum_{\lambda^{{}}\in\Lambda_{N}}\sum_{\lambda^{\prime}\in\Lambda_{N}%
}\overline{A\left(  \lambda\right)  }\,L\left(  \lambda^{\prime}%
-\lambda\right)  A\left(  \lambda^{\prime}\right)  \geq0. \label{eqR-ha.14}%
\end{equation}
Let $\lambda=l+\frac{n}{N^{k}}$, $\lambda^{\prime}=l^{\prime}+\frac{n^{\prime
}}{N^{k}}$, $l,l^{\prime},n,n^{\prime}\in\mathbb{Z}$, $k\in\left\{
0,1,2,\dots\right\}  $. (The general case may be reduced to this by
(\ref{eqR-ha.10}).) Then
\[
L\left(  \lambda^{\prime}-\lambda\right)  =\int_{\mathbb{T}}R^{k}\left(
e\left(  N^{k}\left(  l^{\prime}-l\right)  +n^{\prime}-n\right)  h\right)
\,d\mu,
\]
so%
\begin{align*}
&  \sum_{l^{{}},n}\sum_{l^{\prime},n^{\prime}}\overline{A\left(  l,n\right)
}\,A\left(  l^{\prime},n^{\prime}\right)  L\left(  \frac{N^{k}\left(
l^{\prime}-l\right)  +n^{\prime}-n}{N^{k}}\right) \\
&  \qquad=\sum_{l^{{}},n}\sum_{l^{\prime},n^{\prime}}\overline{A\left(
l,n\right)  }\,A\left(  l^{\prime},n^{\prime}\right)  \int_{\mathbb{T}}%
R^{k}\left(  e\left(  N^{k}\left(  l^{\prime}-l\right)  +n^{\prime}-n\right)
h\right)  \,d\mu\\
&  \qquad=\int_{\mathbb{T}}R^{k}\left(  \left|  \sum_{l,n}A\left(  l,n\right)
e\left(  N^{k}l+n\right)  \right|  ^{2}h\right)  \,d\mu,
\end{align*}
and this last term is $\geq0$ since $R^{k}$ takes nonnegative functions to
nonnegative functions, i.e., $R^{k}$ is positivity-preserving.

It then follows from a theorem of Akhiezer, \cite{Akh65},
\cite{Nel59},
and Kolmogorov,
\cite{Kol77} or
\cite{CFS82}, that there is a unique measure $\nu=\nu\left(  m_{0},h\right)  $
on $K_{N}$ such that%
\[
\int_{K_{N}}e\left(  \frac{n}{N^{k}}\right)  \,d\nu=L\left(  \frac{n}{N^{k}%
}\right)  ,
\]
and therefore%
\begin{equation}
\int_{K_{N}}e_{\lambda}\left(  x\right)  \,d\nu\left(  x\right)  =L\left(
\lambda\right)  \label{eqR-ha.16}%
\end{equation}
for all $\lambda\in\Lambda_{N}$ ($=\mathbb{Z}\left[  \frac{1}{N}\right]  $).
Details of the construction will be given below.

\emph{Ad} (\ref{ThmR-ha.1(3)}): Let $\nu=\nu\left(  m_{0},h\right)  $ be the
measure from (\ref{ThmR-ha.1(2)}). The proof that the operator $U$ on
$L^{2}\left(  K_{N},\nu\right)  $, and the representation $\pi$ of $L^{\infty
}\left(  \mathbb{T}\right)  $ on $L^{2}\left(  K_{N},\nu\right)  $, are given
by formulas (\ref{eqR-ha.6})--(\ref{eqR-ha.8}) follows from the corresponding
assertion in the proof of Theorem \ref{Thmax+b.3}; see Chapter \ref{Poof}
above. The correspondence which makes the connection to Chapter \ref{Poof} is
the identification%
\begin{equation}
U^{-k}e_{n}U^{k}\sim e\left(  \frac{n}{N^{k}}\right)  \label{eqR-ha.15}%
\end{equation}
from (\ref{ThmR-ha.1(1)}), and the basis property of $\left\{  e_{\lambda}%
\mid\lambda\in\Lambda_{N}\right\}  $ in $L^{2}\left(  K_{N},\nu\right)  $.

To show that the Kolmogorov construction applies, and to prove the uniqueness
part of the theorem, we must identify a projective system of measures, as
described, for example, in \cite[Proposition 27.8, page 124]{Par77}. Since,
for each $k$, we have the exponentials $\left\{  e\left(  \frac{n}{N^{k}%
}\right)  \mid n\in\mathbb{Z}\right\}  $ span the functions on $\mathbb{R}$
with period $\left(  2\pi\right)  N^{k}$, we will work with this scale of
periodic functions, $k=0,1,2,\dots$. Functions with period $\left(
2\pi\right)  N^{k}$ will be identified with functions on $\mathbf{X}%
_{k}=\mathbb{R}\diagup\left(  2\pi\right)  N^{k}\mathbb{Z}\simeq\left[  -\pi
N^{k},\pi N^{k}\right\rangle $ with the case $k=0$, $\mathbf{X}_{0}%
\simeq\mathbb{T}$. Restricting to the continuous case, we note that every
$f\in C\left(  \mathbf{X}_{k}\right)  $ has the representation $f\left(
z^{N^{k}}\right)  =F\left(  z\right)  $ for $F\in C\left(  \mathbb{T}\right)
$. The natural maps $\varphi_{k,k+l}$ defined by%
\[
\varphi_{k,k+l}\colon\mathbf{X}_{k+l}\ni z\longmapsto z^{N^{l}}\in
\mathbf{X}_{k}%
\]
then yield a commutative diagram of maps%
\[
\begin{minipage}{72pt}\mbox{\includegraphics
[bbllx=0bp,bblly=-36bp,bburx=72bp,bbury=72bp,width=69.8pt]
{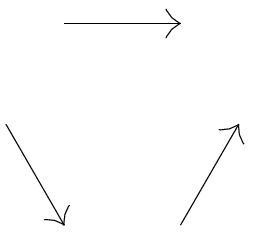}}\llap{\setlength{\unitlength}{32pt}\begin{picture}%
(2.25,3.375)(-1.125,-2.0625)
\put(1.5,0.866){\makebox(0,0){\rlap
{$\displaystyle \mathbf{X}_{k}\;.$}\phantom{$\displaystyle \mathbf{X}$}}}
\put(-1.5,0.866){\makebox(0,0){\rlap
{$\displaystyle \mathbf{X}_{k+l+m}$}\phantom{$\displaystyle \mathbf{X}$}}}
\put(0,-1.754){\makebox(0,0){\rlap
{$\displaystyle \mathbf{X}_{k+m}$}\phantom{$\displaystyle \mathbf{X}$}}}
\put(0,1.116){\makebox(0,0)[b]{$\scriptstyle \varphi_{k,k+l+m}^{{}}$}}
\put(-1,-0.444){\makebox(0,0)[r]{$\scriptstyle \varphi_{k+m,k+l+m}^{{}}$}}
\put(1,-0.444){\makebox(0,0)[l]{$\scriptstyle \varphi_{k,k+m}^{{}}$}}
\end{picture}}\end{minipage}
\]
Recalling the moment sequence in (\ref{eqR-ha.13}), i.e., $L\left(  \frac
{n}{N^{k}}\right)  =\int_{\mathbb{T}}R^{k}\left(  e_{n}h\right)  \,d\mu$, we
see that, for each $k$, a measure $\nu_{k}$ on $\mathbf{X}_{k}$ is determined
uniquely from the trigonometric moment problem (see \cite{Akh65}). The measure
$\nu_{k}$ is \emph{unique,} as it is known that the trigonometric moment
problem is \emph{determined,} i.e., the measure is determined uniquely from
its moments. The consistency condition which is required in the Kolmogorov
construction may be stated in the form, as an identity for $f\in C\left(
\mathbf{X}_{k}\right)  $:%
\[
\int_{\mathbf{X}_{k+l}}f\left(  z^{N^{l}}\right)  \,d\nu_{k+l}=\int
_{\mathbf{X}_{k}}f\,d\nu_{k}.
\]
Recalling that each of the measures $\nu_{k}$ derives from a moment problem,
this identity is equivalent to an identity on $F\in C\left(  \mathbb{T}%
\right)  $, \emph{viz.:}%
\[
\int_{\mathbb{T}}R^{k+l}\left(  F\left(  z^{N^{l}}\right)  h\left(
\,\cdot\,\right)  \right)  \left(  z\right)  \,d\mu\left(  z\right)
=\int_{\mathbb{T}}R^{k}\left(  Fh\right)  \left(  z\right)  \,d\mu\left(
z\right)  ,
\]
where $\mu$ is the Haar measure on $\mathbb{T}$. But this is the identity
which we derived above (by induction, starting with $l=1$ in (\ref{eqR-ha.12}%
).) Hence the Kolmogorov construction determines a measure $\nu$ uniquely. It
is a measure on the projective limit of the systems $\left(  \mathbf{X}%
_{k},\nu_{k}\right)  $. But the compact group $K_{N}$ was identified earlier
with this \emph{projective} limit, $\varprojlim\left\{  \mathbf{X}_{k}\right\}  $. Of
course, $C\left(  K_{N}\right)  $ will then be the \emph{injective} limit of
the system $\left\{  C\left(  \mathbf{X}_{k}\right)  \right\}  $. To see that
the measure $\nu$ on $K_{N}$ is unique, it only remains to invoke the
uniqueness part of the Kolmogorov construction for systems of probability
measures; see \cite[Prop.\ 27.8]{Par77} for the details on that point.

Let $\pi=\pi_{\nu}$ be the representation of $\mathfrak{A}_{N}$ on $L^{2}\left(
K_{N},\nu\right)  $ which is induced by the measure $\nu=\nu\left(
m_{0},h\right)  $ which we just constructed. Recall $m_{0}$ and $h$ are given.
If $R=R_{m_{0}}$ is the Ruelle operator, we have $R_{m_{0}}\left(  h\right)
=h$ at the outset. From the results in Chapter \ref{Cocy}, this means that%
\[
\left(  S_{0}f\right)  \left(  z\right)  :=m_{0}\left(  z\right)  f\left(
z^{N}\right)
\]
defines an \emph{isometry} in $L^{2}\left(  h\right)  $ ($=L^{2}\left(
\mathbb{T},h\,d\mu\right)  $ where $\mu$ is the Haar measure on $\mathbb{T}$).
We will show that the unitary operator $\pi\left(  U\right)  $ on
$L^{2}\left(  K_{N},\nu\right)  $ arises as an extension of $S_{0}$ when
$L^{2}\left(  h\right)  $ is identified (isometrically) with an invariant
subspace in $L^{2}\left(  K_{N},\nu\right)  $. Since clearly $S_{0}%
\openone=m_{0}$, it will follow immediately that $\pi\left(  U\right)
\openone=m_{0}$ ($\in L^{2}\left(  h\right)  \subset L^{2}\left(  K_{N}%
,\nu\right)  $).
\end{proof}

\begin{lemma}
\label{LemR-ha.2}Let $m_{0}$ and $h$ be as described in the statement of
Theorem \textup{\ref{ThmR-ha.1}},
and assume that $m_{0}$ is non-singular. Let $\nu=\nu\left(
m_{0},h\right)  $ be the corresponding measure on $K_{N}=\left(
\mathbb{Z}\left[  \frac{1}{N}\right]  \right)  \sphat{}$. Then $L^{2}\left(
h\right)  =L^{2}\left(  \mathbb{T},h\,d\mu\right)  $ embeds isometrically in
$L^{2}\left(  K_{N},\nu\right)  $ and the unitary operator $\pi_{\nu}\left(
U\right)  $ is an extension \textup{(}or power dilation\textup{)} to
$L^{2}\left(  K_{N},\nu\right)  $ of the isometry $S_{0}$ in $L^{2}\left(
h\right)  $.
\end{lemma}

\begin{proof}
Since $\mathbb{Z\hookrightarrow Z}\left[  \frac{1}{N}\right]  $ by the natural
inclusion, we have $K_{N}\hookrightarrow\mathbb{T}=\left(  \mathbb{Z}\right)
\sphat{}$ with the embedding from Pontryagin duality; see \cite{Rud90}.
Restricting functions on $\mathbb{T}$ to $K_{N}$, we then get the
identification of $L^{\infty}\left(  \mathbb{T}\right)  $ with a subspace of
$L^{\infty}\left(  K_{N}\right)  $. Since, for $f\in L^{\infty}\left(
\mathbb{T}\right)  $, we have%
\[
\left\|  f\right\|  _{L^{2}\left(  h\right)  }^{2}=\int_{\mathbb{T}}\left|
f\right|  ^{2}h\,d\mu=\left\|  f\right\|  _{L^{2}\left(  K_{N},\nu\right)
}^{2},
\]
it follows that $L^{2}\left(  h\right)  $ embeds isometrically into
$L^{2}\left(  K_{N},\nu\right)  $ as claimed. Since \linebreak $R_{m_{0}%
}\left(  h\right)  =h$, $S_{0}$ will be isometric on $L^{2}\left(  h\right)
$, and therefore identify with a partial isometry in $L^{2}\left(  K_{N}%
,\nu\right)  $.

Recall that the elements $\left\{  U^{-i}fU^{i}\mid f\in L^{\infty}\left(
\mathbb{T}\right)  ,\;i=0,1,2,\dots\right\}  $ in $\mathfrak{A}_{N}$ generate the
abelian algebra $\mathcal{A}_{N}$, and we showed that $K_{N}$ is the Gelfand
space of $\mathcal{A}_{N}$. Hence we may identify $\mathcal{A}_{N}$ also with
a linear subspace in $L^{2}\left(  K_{N},\nu\right)  $, and on this subspace
we set%
\[
\pi_{\nu}\left(  U\right)  \left(  U^{-\left(  i+1\right)  }fU^{i+1}\right)
:=\left(  U^{-i}fU^{i}\right)  m_{0}=U^{-i}\left(  f\left(  z\right)
m_{0}\left(  z^{N^{i}}\right)  \right)  U^{i}.
\]
To see that $\pi_{\nu}\left(  U\right)  $ is well defined and isometric (in
$L^{2}\left(  K_{N},\nu\right)  $) we must check that%
\[
\int_{K_{N}}\left|  U^{-\left(  i+1\right)  }fU^{i+1}\right|  ^{2}\,d\nu
=\int_{K_{N}}\left|  U^{-i}fU^{i}m_{0}\right|  ^{2}\,d\nu,
\]
which is equivalent to
\[
\int_{\mathbb{T}}R^{i+1}\left(  \left|  f\right|  ^{2}h\right)  \,d\mu
=\int_{\mathbb{T}}\left|  m_{0}\right|  ^{2}R^{i}\left(  \left|  f\right|
^{2}h\right)  \,d\mu.
\]
Since we already checked this identity, the result follows. It follows from the
formula for $\pi_{\nu}\left(  U\right)  $ that it maps \emph{onto}
$L^{2}\left(  K_{N},\nu\right)  $ if $m_{0}$ is given to be non-singular,
i.e., does not vanish on a subset of positive measure in $\mathbb{T}$. Hence
$\pi_{\nu}\left(  U\right)  $ is a unitary extension (or \emph{power
dilation\/}) as claimed in the lemma.
\end{proof}

\begin{remark}
\label{RemKean.3}
Alternatively, $L^{2}\left(  K_{N}\right)  $ may be defined relative to the
Haar measure $\mu_{N}$ on $K_{N}$. This Haar measure in turn is determined
uniquely by the \emph{ansatz} ($\lambda\in\Lambda_{N}$):%
\begin{equation}
\int_{K_{N}}e_{\lambda}\,d\mu_{N}=\delta_{\lambda}=%
\begin{cases}
0 &\text{if }\lambda\neq0, \\
1 &\text{if }\lambda=0.
\end{cases}%
\label{eqR-ha.17}%
\end{equation}
The fact that (\ref{eqR-ha.17}) determines a unique measure on $K_{N}$ follows
from the same argument which we used in (\ref{ThmR-ha.1(2)}) above. This
measure $\mu_{N}$ will be translation-invariant on $K_{N}$ by the following
calculation. (Hence it must be the Haar measure by the uniqueness theorem!) We
have for $\chi_{0}^{{}}\in K_{N}=\left(  \Lambda_{N}\right)  \sphat{}$:%
\[
\int_{K_{N}}\left(  \chi_{0}^{{}}\chi\right)  \left(  \lambda\right)
\,d\mu_{N}\left(  \chi\right)  =\chi_{0}^{{}}\left(  \lambda\right)  \int_{K_{N}%
}e_{\lambda}\,d\mu_{N}=\chi_{0}^{{}}\left(  \lambda\right)  \delta_{\lambda
}^{{}}=\int_{K_{N}}e_{\lambda}\,d\mu_{N}.
\]
\end{remark}
\section{\label{Sign}Signed solutions to $R\left(  f\right)  =f$}

In Example \ref{ExaWave.3(2)} we considered a wavelet filter $m_{0}$ of order
$2$ and a positive solution $g$ to $R_{m_{0}}\left(  g\right)  =g$ such that
$g$ had the following order-$3$ symmetry:%
\begin{equation}
\sum_{k=0}^{2}g\left(  e^{i\frac{2\pi k}{3}}z\right)  =1. \label{eqSign.1}%
\end{equation}
This means that $g$ satisfies a special case of the condition from Chapter
\ref{Kean} above, corresponding to $N=3$. Specifically, let%
\[
\mathbb{T}\ni z\underset{T_{3}}{\longmapsto}z^{3}\in\mathbb{T};
\]
then%
\[
\sum_{T_{3}w=z}g\left(  w\right)  =1,\qquad z\in\mathbb{T}.
\]
In this chapter, we will study a more general scaling duality which relates
scaling of order $N$ to that of order $p$, where $N$ and $p$ are positive
integers which are given and mutually prime, $\left(  N,p\right)  =1$, i.e.,
no common divisors other than $1$. We will then have a pair of Ruelle
operators and a specific duality between the eigenvalue problems for the
respective operators. It is also a concrete instance
of a case
when the dimension%
\begin{equation}
\dim\left\{  f\in L^{1}\left(  \mathbb{T}\right)  \mid R\left(  f\right)
=f\right\}  \label{eqSign.2}%
\end{equation}
can be calculated; in this case, it is shown to be equal to the number of
orbits for a certain finite dihedral group action.

The setting of this duality will be two given wavelet filters $m_{0}$ and
$m_{p}$ related as follows: It will be assumed that $\left(  N,p\right)  =1$,%
\begin{equation}
\sum_{w^{N}=z}\left|  m_{0}\left(  w\right)  \right|  ^{2}  
=N,\label{eqSign.3}
\end{equation}
and
\begin{equation}
m_{p}\left(  z\right)  =m_{0}\left(  z^{p}\right)  ,\qquad z\in\mathbb{T}.
\label{eqSign.4}
\end{equation}
Since $\left(  N,p\right)  =1$, it follows that $m_{p}$ will also satisfy the
scale-$N$ condition (\ref{eqSign.3}). (The simplest case of this is the one in
Example \ref{ExaWave.3(1)}, when $N=2$, $p=3$, $m_{0}\left(  z\right)
=\frac{1}{\sqrt{2}}\left(  1+z\right)  $, $m_{3}\left(  z\right)  =\frac
{1}{\sqrt{2}}\left(  1+z^{3}\right)  $.) We shall need both Ruelle operators
$R_{0}$ and $R_{p}$ constructed from $m_{0}$ and $m_{p}$ when the two filters
are related through (\ref{eqSign.4}), i.e.,%
\begin{align*}
R_{0}f\left(  z\right)   &  =\frac{1}{N}\sum_{w^{N}=z}\left|  m_{0}\left(
w\right)  \right|  ^{2}f\left(  w\right) \\%
\intertext{and}%
R_{p}f\left(  z\right)   &  =\frac{1}{N}\sum_{w^{N}=z}\left|  m_{p}\left(
w\right)  \right|  ^{2}f\left(  w\right)  .
\end{align*}

\begin{lemma}
\label{LemSign.1}Let $N$, $p$ be positive integers $\geq2$ such that $\left(
N,p\right)  =1$, and let $m_{0}$, $m_{p}$ be given wavelet filters of order
$N$ and related \emph{via} \textup{(\ref{eqSign.4}).} Let $R_{0}$ and $R_{p}$
be the respective Ruelle operators. Let $\alpha_{N}\in\operatorname*{Aut}%
\left(  \mathbb{Z}_{p}\right)  $ be the automorphism $i\mapsto Ni$, passed to
$\mathbb{Z}_{p}=\mathbb{Z}\diagup p\mathbb{Z}$, and let $\rho_{p}%
:=e^{i\frac{2\pi}{p}}$.

\begin{enumerate}
\item \label{LemSign.1(1)}Let $f\in L^{1}\left(  \mathbb{T}\right)  $ satisfy
$R_{p}f=f$, and set
\[
F_{0}\left(  z\right)  =\frac{1}{p}\sum_{j=0}^{p-1}f\left(  \rho_{p}%
^{j}\,z\right)  .
\]
Then $F_{0}$ is of the form $F_{0}\left(  z\right)  =H_{0}\left(
z^{p}\right)  $, $H_{0}\in L^{1}\left(  \mathbb{T}\right)  $, and
$R_{0}\left(  H_{0}\right)  =H_{0}$.

\item \label{LemSign.1(2)}Let $\left\langle j,k\right\rangle :=e^{i\frac{2\pi
jk}{p}}$, and
\begin{equation}
F_{k}\left(  z\right)  =\frac{1}{p}\sum_{j=0}^{p-1}\overline{\left\langle
j,k\right\rangle }\,f\left(  \rho_{p}^{j}\,z\right)  . \label{eqSign.5}%
\end{equation}
Then $F_{k}$ is of the form%
\begin{equation}
F_{k}\left(  z\right)  =z^{k}H_{k}\left(  z^{p}\right)  ,\qquad z\in
\mathbb{T},\;H_{k}\in L^{1}\left(  \mathbb{T}\right)  ,\label{eqSign.6}%
\end{equation}
with%
\begin{equation}
R_{0}\left(  H_{k}\right)  =H_{\alpha_{N}^{-1}\left(  k\right)  },\qquad
k\in\mathbb{Z}_{p}. \label{eqSign.7}%
\end{equation}
\end{enumerate}
\end{lemma}

\begin{proof}
Let $f_{j}\left(  z\right)  :=f\left(  \rho_{p}^{j}\,z\right)  $. Then
\begin{equation}
R_{p}\left(  f_{j}\right)  =f_{\alpha_{N}\left(  j\right)  }. \label{eqSign.8}%
\end{equation}
Indeed
\begin{align*}
R_{p}\left(  f_{j}\right)  \left(  z\right)   &  =\frac{1}{N}\sum_{w^{N}%
=z}\left|  m_{0}\left(  w^{p}\right)  \right|  ^{2}f\left(  \rho_{p}%
^{j}\,w\right) \\
&  =\frac{1}{N}\sum_{w^{N}=\rho_{p}^{Nj}z}\left|  m_{0}\left(  w^{p}\right)
\right|  ^{2}f\left(  w\right) \\
&  =R_{p}\left(  f\right)  \left(  \rho_{p}^{Nj}\,z\right) \\
&  =f\left(  \rho_{p}^{Nj}\,z\right)  =f_{\alpha_{N}\left(  j\right)  }\left(
z\right)  ,
\end{align*}
which is the assertion (\ref{eqSign.8}). Since
\[
\sum_{j=0}^{p-1}\left(  z\rho_{p}^{j}\right)  ^{n}=%
\begin{cases}
pz^{n} &\text{if }p|n,  \\
0 &\text{if }p\nmid n,
\end{cases}%
\]
there is an $H_{0}\in L^{1}\left(  \mathbb{T}\right)  $ such that $F_{0}\left(
z\right)  =H_{0}\left(  z^{p}\right)  $, and%
\begin{align*}
H_{0}\left(  z^{p}\right)   &  =F_{0}\left(  z\right)  =\sum_{j=0}^{p-1}%
f_{j}\left(  z\right) \\
&  =\sum_{j=0}^{p-1}f_{\alpha_{N}\left(  j\right)  }\left(  z\right)
=\sum_{j=0}^{p-1}R_{p}\left(  f_{j}\right)  \left(  z\right)  =\left(
R_{p}F_{0}\right)  \left(  z\right) \\
&  =\frac{1}{N}\sum_{w^{N}=z}\left|  m_{p}\left(  w\right)  \right|
^{2}F_{0}\left(  w\right) \\
&  =\frac{1}{N}\sum_{w^{N}=z}\left|  m_{0}\left(  w^{p}\right)  \right|
^{2}H_{0}\left(  w^{p}\right) \\
&  =\frac{1}{N}\sum_{w^{N}=z^{p}}\left|  m_{0}\left(  w\right)  \right|
^{2}H_{0}\left(  w\right) \\
&  =R_{0}\left(  H_{0}\right)  \left(  z^{p}\right)  ,
\end{align*}
which yields the desired identity $R_{0}\left(  H_{0}\right)  =H_{0}$ in
(\ref{LemSign.1(1)}).

The proof of (\ref{LemSign.1(2)}) is quite similar and will only be sketched.
We have%
\[
\frac{1}{p}\sum_{j=0}^{p-1}\rho_{p}^{-jk}\left(  \rho_{p}^{j}\,z\right)  ^{n}=%
\begin{cases}
z^{n} &\text{if }n\equiv k\mod{p},  \\
0 &\text{if }n\not\equiv k\mod{p}.
\end{cases}%
\]
Hence%
\[
F_{k}\left(  z\right)  =\sum_{l}c_{l}z^{k+lp}=z^{k}H_{k}\left(  z^{p}\right)
\]
for some $H_{k}\in L^{1}\left(  \mathbb{T}\right)  $. The argument from above
yields $F_{k}=R_{p}\left(  F_{\alpha_{N}\left(  k\right)  }\right)  $, and
therefore%
\begin{align*}
H_{k}\left(  z^{p}\right)   &  =z^{-k}F_{k}\left(  z\right)  =z^{-k}%
R_{p}\left(  F_{\alpha_{N}\left(  k\right)  }\right)  \left(  z\right) \\
&  =z^{-k}\frac{1}{N}\sum_{w^{N}=z}\left|  m_{0}\left(  w^{p}\right)  \right|
^{2}F_{\alpha_{N}\left(  k\right)  }\left(  w\right) \\
&  =\frac{1}{N}\sum_{w^{N}=z^{p}}\left|  m_{0}\left(  w\right)  \right|
^{2}H_{\alpha_{N}\left(  k\right)  }\left(  w\right) \\
&  =R_{0}\left(  H_{\alpha_{N}\left(  k\right)  }\right)  \left(
z^{p}\right)  ,
\end{align*}
which is the desired identity $R_{0}\left(  H_{\alpha_{N}\left(  k\right)
}\right)  =H_{k}$, or equivalently $R_{0}\left(  H_{k}\right)  =H_{\alpha
_{N}^{-1}\left(  k\right)  }$, $k\in\mathbb{Z}_{p}$.
\end{proof}

Consider the action $\tau_{p}$ on $L^{1}\left(  \mathbb{T}\right)  $ given by%
\begin{align*}
\left(  \tau_{p}f\right)  \left(  z\right)   &  =f\left(  \rho_{p}z\right)
,\qquad z\in\mathbb{T},\\%
\intertext{and}%
\tau_{p}^{j}f\left(  z\right)   &  =f\left(  \rho_{p}^{j}\,z\right)  ,\qquad
j\in\mathbb{Z}_{p}.
\end{align*}
If $V\subset L^{1}\left(  \mathbb{T}\right)  $ is a given subspace, we set%
\begin{align}
V^{\tau_{p}}  &  =\left\{  f\in V\mid\tau_{p}f=f\right\}  ,\label{eqSign.9}\\%
\intertext{and}%
V\left(  z^{p}\right)   &  =\left\{  f\in L^{1}\left(  \mathbb{T}\right)
\mid\exists h\in V\;\mathrm{s.t.}\;f\left(  z\right)  =h\left(  z^{p}\right)
\right\}  . \label{eqSign.10}%
\end{align}
Clearly then%
\[
V^{\tau_{p}}=\left\{  f\in V\biggm|\exists F\in L^{1}\left(  \mathbb{T}%
\right)  \;\mathrm{s.t.}\;f=\frac{1}{p}\sum_{j=0}^{p-1}\tau_{p}^{j}\left(
F\right)  \right\}  .
\]

Return to the setting of the lemma, i.e., two given wavelet filters $m_{0}$,
$m_{p}$ related \emph{via} (\ref{eqSign.4}) and corresponding Ruelle operators
$R_{0}$ and $R_{p}$.

We have the following reciprocity for the eigenspaces.

\begin{corollary}
\label{CorSign.2}Let $m_{p}\left(  z\right)  =m_{0}\left(  z^{p}\right)  $,
and let%
\begin{align*}
V_{0}  &  =\left\{  f\in L^{1}\left(  \mathbb{T}\right)  \mid R_{0}f=f\right\}
\\%
\intertext{and}%
V_{p}  &  =\left\{  f\in L^{1}\left(  \mathbb{T}\right)  \mid R_{p}%
f=f\right\}  .
\end{align*}
Then%
\begin{equation}
V_{p}^{\tau_{p}}=V_{0}\left(  z^{p}\right)  . \label{eqSign.11}%
\end{equation}
\end{corollary}

\begin{proof}
The inclusion $\subset$ already follows from Lemma \ref{LemSign.1}%
(\ref{LemSign.1(1)}). To prove $\supset$, let $f\left(  z\right)  =h\left(
z^{p}\right)  $ where $h\in L^{1}\left(  \mathbb{T}\right)  $ and $R_{0}h=h$.
Then%
\begin{align*}
\left(  R_{p}f\right)  \left(  z\right)   &  =\frac{1}{N}\sum_{w^{N}=z}\left|
m_{p}\left(  w\right)  \right|  ^{2}f\left(  w\right) \\
&  =\frac{1}{N}\sum_{w^{N}=z}\left|
m_{0}\left(  w^{p}\right)  \right|  ^{2}h\left(  w^{p}\right) \\
&  =\frac{1}{N}\sum_{w^{N}=z^{p}}\left|  m_{0}\left(  w\right)  \right|
^{2}h\left(  w\right) \\
&  =R_{0}\left(  h\right)  \left(  z^{p}\right)  =h\left(  z^{p}\right) \\
&  =f\left(  z\right)  ,
\end{align*}
proving $f\in V_{p}$. Since $\tau f=f$ is clear, the result follows.
\end{proof}

In many examples, scaling functions constructed from given wavelet filters may
not have orthogonal translates $\left\{  \varphi\left(  \,\cdot\,-k\right)
\mid k\in\mathbb{Z}\right\}  $ in $L^{2}\left(  \mathbb{R}\right)  $. It is
known \cite{Hor95} that if $\varphi$ is constructed from $m_{0}$, then the
orthogonality holds if and only if the eigenspace $\left\{  f\mid R_{m_{0}%
}f=f\right\}  $ is one-dimensional. In this case, we say that $m_{0}$ is
\emph{pure.} Let $N$ and $p$ be given, $\left(  N,p\right)  =1$, and let
$m_{p}\left(  z\right)  :=m_{0}\left(  z^{p}\right)  $ with $m_{0}$ pure. Then
it follows from Corollary \ref{CorSign.2} that $V_{p}^{\tau_{p}}$ ($=\left\{
f\in V_{p}\mid\tau_{p}f=f\right\}  $) must be one-dimensional. We saw in
Lemma \ref{LemSign.1} that every $f\in V_{p}$ decomposes uniquely as
\[
f=\sum_{j\in\mathbb{Z}_{p}}f_{j},
\text{\qquad where }\tau_{p}^{{}}f_{j}^{{}}=\rho_{p}%
^{j}\,f_{j}^{{}},\;j\in\mathbb{Z}_{p}. 
\]
The explicit form of this decomposition
is spelled out in Lemma \ref{LemSign.1}(\ref{LemSign.1(2)}). Specifically,
$f_{j}\left(  z\right)  =z^{j}\,h_{j}\left(  z^{p}\right)  $ with $R_{0}\left(
h_{j}\right)  =h_{\alpha_{N}^{-1}\left(  j\right)  }$, $j\in\mathbb{Z}_{p}$.
When $m_{0}$ is pure, it follows that a basis for $V_{p}$ may be labelled by
the finite orbits $j\mapsto Nj\mapsto N^{2}j\mapsto\dots\rightarrow
N^{k}j\rightarrow j$ in $\mathbb{Z}_{p}=\mathbb{Z}\diagup p\mathbb{Z}$.
Specifically, let $\mathcal{O}_{k}$ be such an orbit and let $f\in V_{p}$;
then $f_{\mathcal{O}_{k}}\left(  z\right)  =\sum_{s\in\mathcal{O}_{k}}f\left(
\rho_{p}^{s}\,z\right)  $ is in $V_{p}$, and the finite Fourier analysis on
$\mathbb{Z}_{p}$ of the lemma shows that the orbits parametrize a basis for
$V_{p}$ when $m_{0}$ is pure and $\left(  N,p\right)  =1$.

These same orbits were found in \cite{BrJo96b} to parametrize a class of
permutation representations of the Cuntz algebras. Specifically,
\cite[Proposition 8.2 and Remark 8.3]{BrJo96b} lists these orbits for $N=2$
and $p\in\mathbb{N}$ odd, for selected values of $p$. While there is a pattern
to this orbit counting, we do not have a general formula valid for all odd
values of $p$. This is just the case $N=2$, and more general pairs $N,p$ such
that $\left(  N,p\right)  =1$ may be more difficult. For $N=2$, it can be seen
that the period of an orbit starting at $j\in\mathbb{Z}_{p}$ is the order of
$2$ modulo $p\diagup\gcd\left(  j,p\right)  $.

\begin{acknowledgements}
Helpful discussions with Ola Bratteli,
Steen Pedersen, and Beth Peterson are gratefully
acknowledged, as are excellent typesetting,
diagram design and construction, and
\textsl{Mathematica}
programming by Brian Treadway.
Most of the original research was done
while the author recovered
in the hospital from an accident.
\end{acknowledgements}


\bibliographystyle{BFTALPHA}
\bibliography{jorgen}
\end{document}